# Jassologie
# Une vision originale sur les cartes orientables

*Catherine Mouet*

*Juillet 2013*

## Introduction

La *Jassologie* est une théorie générale sur les cartes issue d'expérimentations artistiques.

Le processus *Jassologique* consiste à associer un mot, appelé *Mot Jassologique*, à une carte planaire cubique (dont tous les sommets sont de valence 3) enracinée (qui possède une arête orientée distinguée) quelconque.
Les *Mots Jassologiques* font partie d'un monoïde libre dont la base, usuellement appelée Alphabet, est composée de 16 caractères spécifiques dont j'ai entièrement inventé le graphisme. Ils obéissent à un certain nombre de règles, dont la *Règle d'Emboîtement* (voir ci-dessous).
On montre que deux cartes planaires cubiques enracinées sont équivalentes si et seulement si les *Mots Jassologiques* auxquels ils correspondent sont les mêmes.
Inversement, à partir d'un *Mot Jassologique* quelconque, on construit une carte planaire cubique enracinée appelée C*arte Géométrique* car tous les sommets et arêtes ont des coordonnées spécifiques dans un Repère Cartésien.
On montre que le *Mot Jassologique* associé à cette carte est le même que le *Mot Jassologique* de départ.
On obtient alors une bijection entre l'ensemble des *Mots Jassologiques* et celui des classes d'équivalence des cartes planaires cubiques enracinées.

Du point de vue de la théorie informatique, le processus *Jassologique* est un algorithme qui n'utilise aucune manipulation d'image : il suffit de prendre la suite ordonnée des voisins qui entourent chaque face d'une carte de départ comme base de donnée, et le *Mot Jassologique* correspondant s'obtient uniquement par traitement de texte.

### La méthode *Jassologique*

Elle repose sur une méthode de décomposition en couches de cellules appelées *Jasses*, où « jas » signifie « dépôt, lie, stratification » en provençal. C'est une décomposition assez naturelle, à l'image des strates de sédimentation des roches, des bandes stratifiées de certaines pierres précieuses comme les agates, ou encore des courbes de niveau, du mouvement ondulatoire de l'eau ou du sable...
Les *Jasses* sont ensuite décomposées en chaînes de cellules appelées *Stratojasses*. Chaque *Stratojasse* a une hauteur bien définie sur l'*Échelle Jassologique* associée à la carte étudiée.
On peut alors visiter chacune des faces de la carte en partant de la face racine (définie par l'arête distinguée), comme une « visite guidée le long des chemins stratifiés » : le *Mot Jassologique* associé à la carte suit l'ordre des faces.
Une coloration spécifique, appelée *Coloration Stratojasse*, permet de mettre en évidence la structure. Bien qu'elle se réalise avec 4 couleurs différentes, elle ne répond pas du tout au célèbre problème des 4 couleurs. Elle est cependant très esthétique et facilite grandement la visualisation d'une carte.

### Les *Mots Jassologiques* : des symboles nouveaux pour la *Règle d'Emboîtement*

Rien n'empêche de remplacer les 16 caractères de l'*Alphabet Jassologique* par des lettres ou des chiffres.
Cependant, les caractères usuels ne sont pas très pratiques à utiliser pour composer les *Mots Jassologiques*, car ceux-ci ont besoin de caractères qui fonctionnent en binôme, comme les couples de parenthèses ouvrants et fermants : **(** et **)**. En particulier, ils doivent vérifier la *Règle d'Emboîtement* des parenthèses que chacun utilise intuitivement : par exemple « **( ( ) ( ) )** » est un *Mot Jassologique Simple*, mais « **( ( ( )** » ou « **( ( ) ) ) ( ( )** » n'en sont pas.
Je pense personnellement qu'il n'est pas non plus très judicieux d'utiliser les parenthèses ou les crochets pour représenter un mot mathématique, car ils ont une fonction syntaxique bien spécifique qu'il est bon de leur laisser. C'est pourquoi je préfère utiliser le binôme *Jassologique* [illegible], ce qui donne « [illegible] » à la place de « **( ( ) ( ) )** ».
En particulier, on montre qu'un seul type de binôme suffit pour représenter n'importe quel arbre ordonné par un unique *Mot Jassologique Simple*.
Mais pour ce qui concerne les cartes planaires cubiques, 6 binômes différents sont nécessaires, ainsi que quatre autres symboles qui ne fonctionnent pas en binôme. Cela fait au total 16 caractères :
Les 12 caractères ouvrants et fermants des 6 binômes : [illegible], [illegible], [illegible], [illegible], [illegible], [illegible].
Les 4 monômes : •, [illegible], [illegible], [illegible].
En particulier, le nombre de binômes qui apparaissent dans un *Mot Jassologique* est égal au nombre de faces de la carte à laquelle il correspond. Grâce à la Caractéristique d'Euler, on en déduit facilement le nombre d'arêtes et de sommets.

# *illustration graphique de la jassologie*

*carte enracinée initiale*

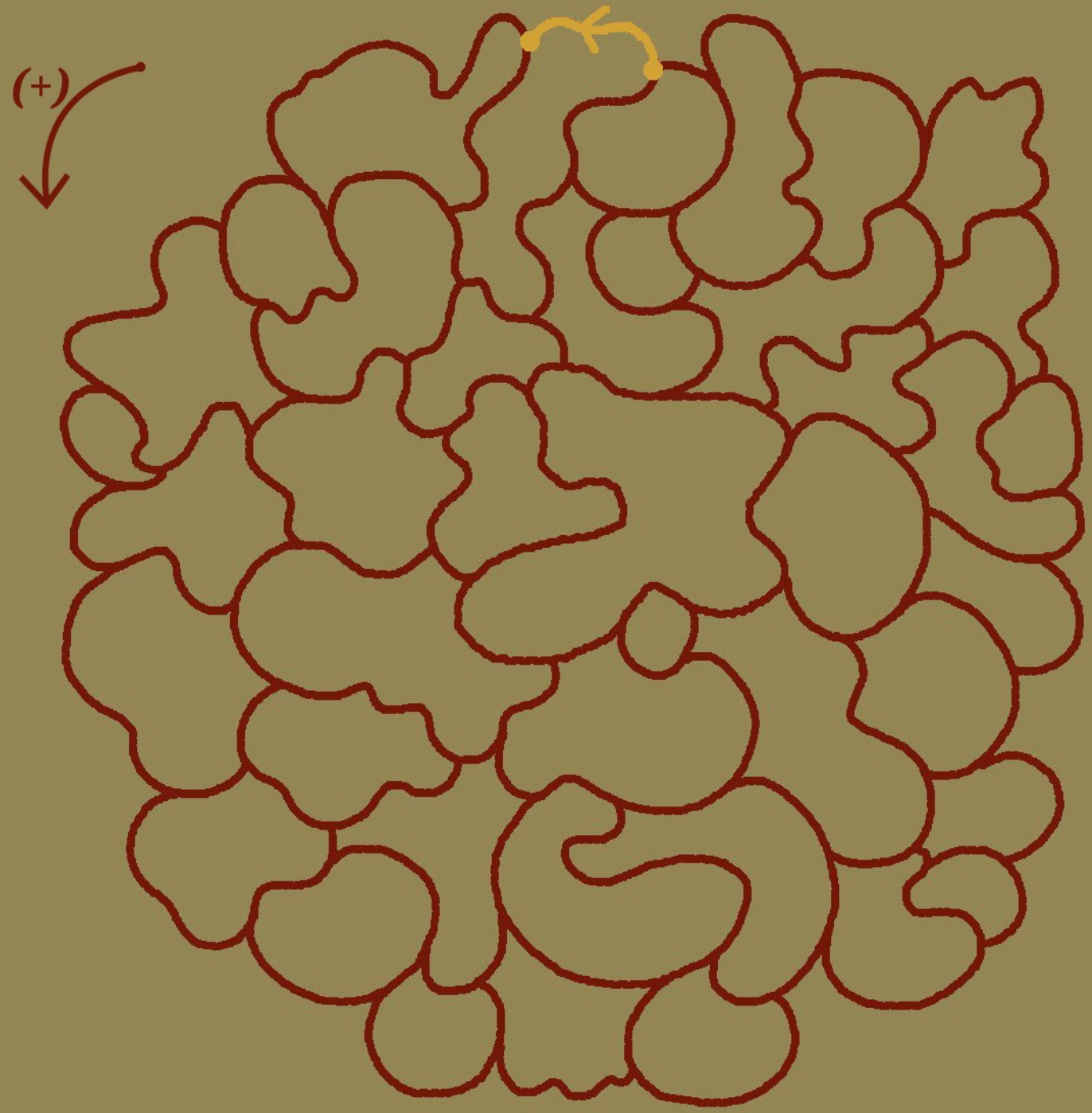

*coloration stratojasse*

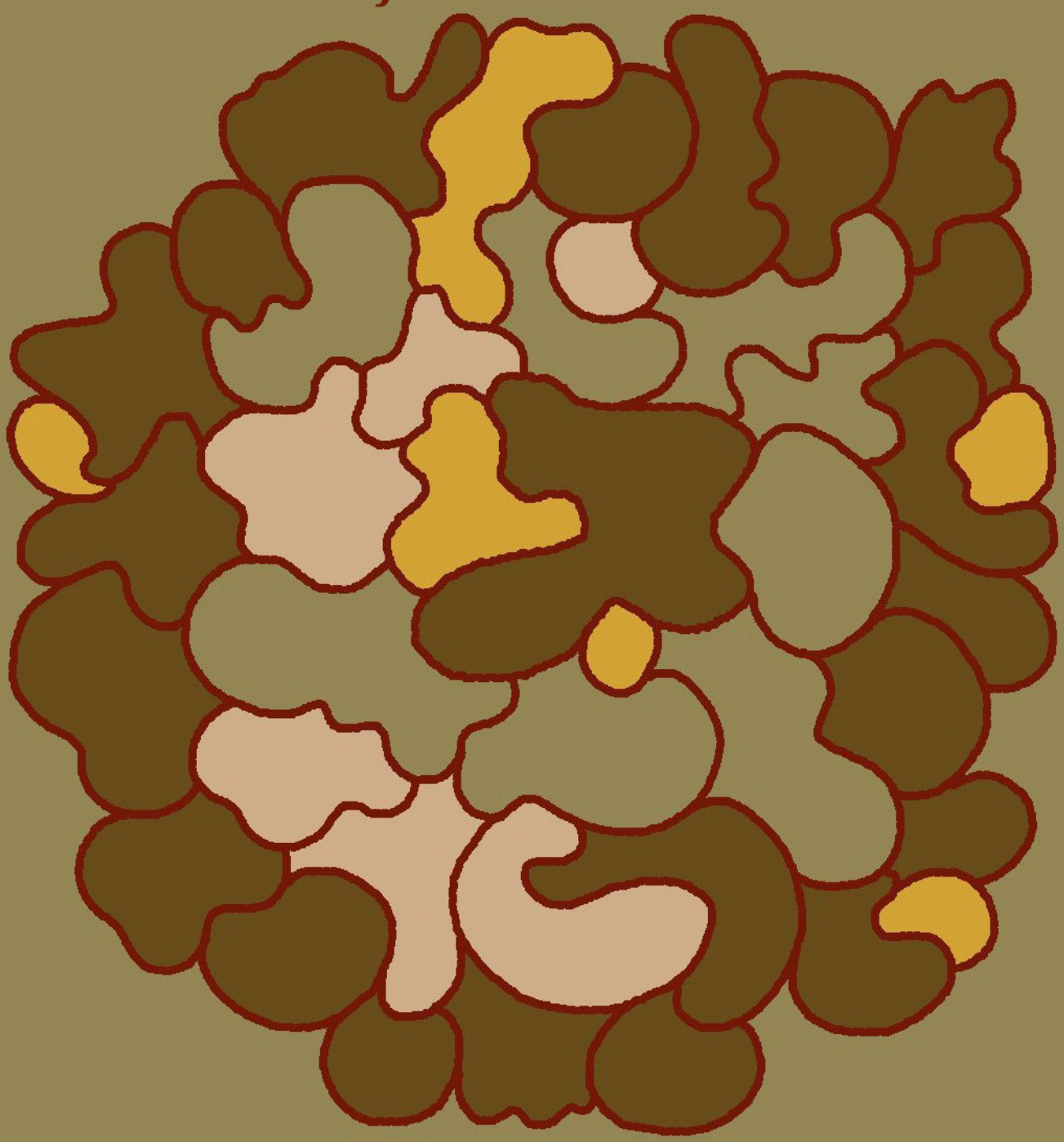

*mot jassologique correspondant*

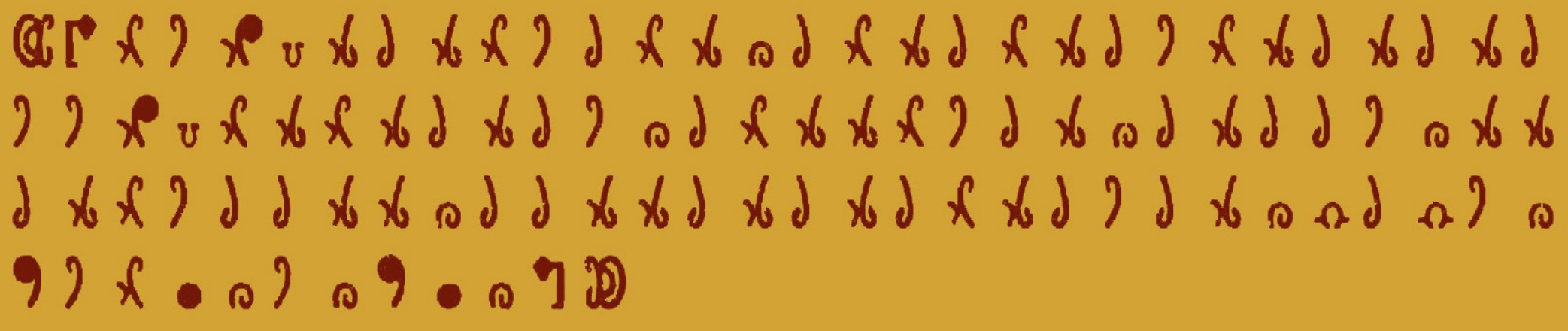

*échelle*

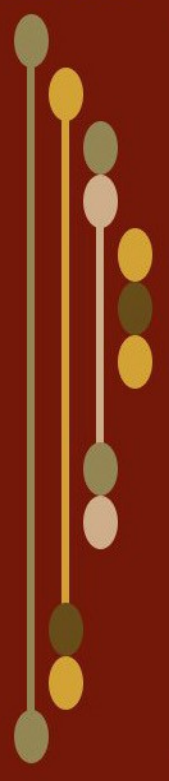

*carte géométrique équivalente à la carte initiale*

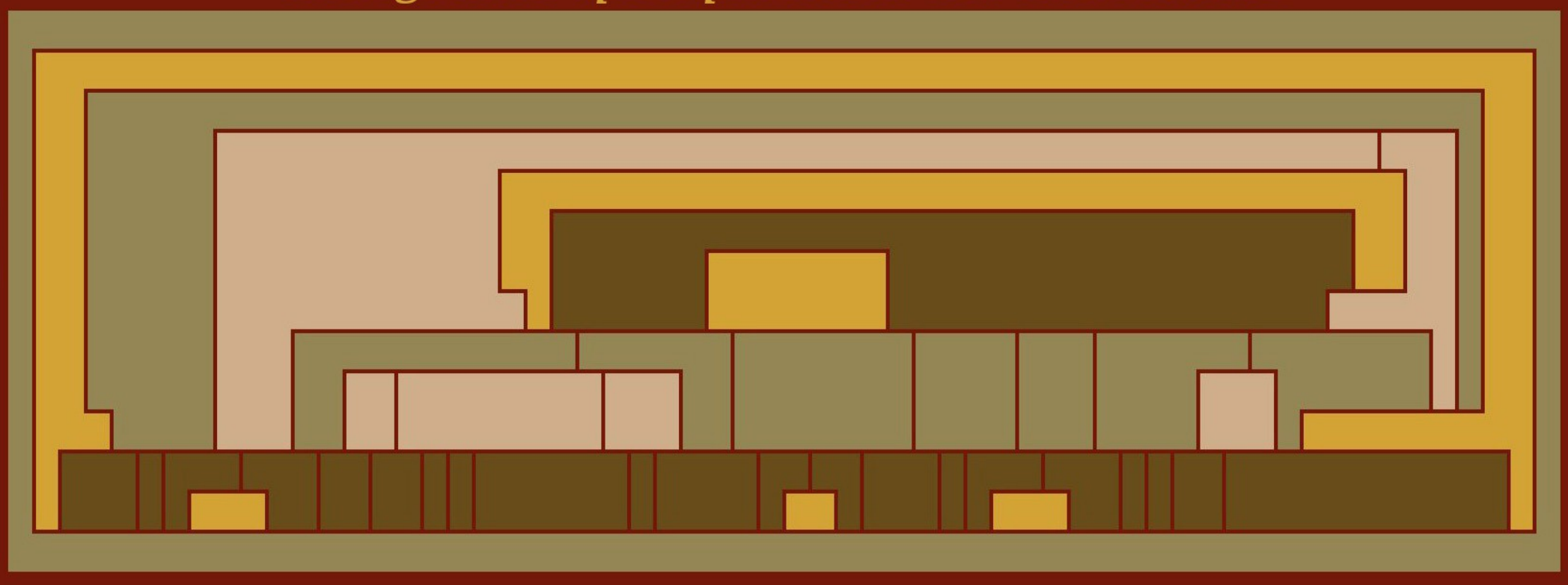

## *catherine mouet, novembre 2012*

# Table des matières

**N.B.** : L'objectif de cet article est de donner les outils de base de la *Jassologie.* Aussi, les démonstrations des propositions et théorèmes n'y figurent pas.

# 1
# Définitions et notations concernant les cartes

## 1.1. Cartes

### 1.1.1. Définition d'une carte :

Ma définition d'une carte est différente de celle définie par William T. Tutte dans les années 1960, car elle s'est basée sur un livre que j'avais trouvé à New York en 2004 (« Topology of Surfaces, Knots, and Manifolds » de Stephan C. Carlson) ainsi que le livre de Robin Wilson sur le théorème des 4 couleurs (« Four Colors Suffice - How the map problem was solved »). Tutte a été le premier mathématicien à tenter de résoudre le célèbre problème des 4 couleurs en répertoriant et énumérant le nombre de cartes en fonction du nombre d'arêtes par exemple. Ma ligne directive est toute autre, basée sur une recherche structurelle : comment trouver une structure intrinsèque commune à toutes les cartes.
William T. Tutte définit une carte comme plongement propre (sans croisement d'arêtes) d'un graphe connexe dans la sphère S (carte planaire) ou le g-Tore (carte de genre g) qui définit des faces simplement connexes.
Stephan C. Carlson définit une carte comme recouvrement de la Sphère ou d'un g-Tore par des polygones qui se rencontrent uniquement au niveau de leurs sommets et arêtes.
Ici, une carte sera définie comme un ensemble de cellules (espaces homéomorphes au disque D) qui recouvrent entièrement la sphère (carte planaire) ou le g-Tore (carte de genre g), et qui doivent se rencontrer uniquement le long de leur bord en un nombre fini de points (sommets) ou de segments de courbe (arêtes dont les extrémités sont des sommets). De plus, une carte est cubique si les sommets sont tous de valence 3, la valence d'un sommet étant le nombre d'arêtes qui ont pour extrémité ce sommet. J'ai trouvé cette définition dans le livre de Robin Wilson. Mais William T. Tutte utilise le terme *Trivalent*.
On vérifie facilement, et Robin Wilson l'explique très bien dans son livre, qu'il suffit de démontrer le théorème des 4 couleurs pour les cartes planaires cubiques.
Les cartes cubiques sont des plongements propres de graphes connexes, mais l'inverse n'est pas vrai.

### 1.1.2. Définition d'une carte enracinée :

Une carte est enracinée si elle possède une arête orientée distinguée (cette définition est utilisée de la même manière par William T. Tutte).

### 1.1.3. Orientation et faces racines :

Soit Ω une carte cubique enracinée de genre g.
Donc Ω est un ensemble de faces ou cellules, homéomorphes au disque D, dont le bord, qui est homéomorphe au cercle, peut être orienté dans le sens trigonométrique (on dira que c'est le sens (+)) ou dans le sens des aiguilles d'une montre (le sens (-)).
Notons λ l'arête orientée distinguée : elle est par définition à l'intersection de deux faces appelées faces racines. Pour l'une d'elle, appelée face racine (-), λ tourne dans le sens (-), alors que pour l'autre, appelée face racine (+), λ tourne dans le sens (+).

### 1.1.4. Connexité :

Soit Ω une carte cubique enracinée de genre g.
On dira qu'un sous-ensemble E de Ω est connexe si et seulement si l'union de ses éléments est connexe.
Pour tout sous-ensemble E de Ω, on définit alors les composantes connexes de E, qui sont uniques.

## 1.2. Equivalence

### 1.2.1. *Bordures* :

Soit Ω une carte cubique de genre g.
Pour toute cellule e de Ω, choisissons un sommet quelconque i(e) de référence sur le bord de e.
La bordure de e à partir de i(e) dans le sens (+), noté B(e), est la suite ordonnée des voisins qui entourent e dans le sens (+) à partir de i(e). Autrement dit, si on note $i(e)_1$, $i(e)_2$, …, $i(e)_n$, la suite des sommets qui se suivent dans le sens(+) sur le bord de e, avec $i(e)_1 = i(e)$ et $n \geq 2$, alors $B(e) = (g_1, g_2, \ldots, g_n)$, où pour tout k, $1 \leq k \leq n$, $g_k$ est l'unique cellule de Ω qui rencontre e le long de l'arête d'extrémités $i(e)_k$ et $i(e)_{k+1}$ (avec $i(e)_{n+1} = i(e)$).

On définit également la bordure étendue de e, notée $\overline{B}(e)$ : $\overline{B}(e) = (g_1, g_2, \ldots, g_n, g_1) = (g_1, g_2, \ldots, g_n, g_{n+1})$, avec $g_{n+1} = g_1$.
On vérifie facilement la propriété suivante : pour toutes cellules a, b et c, si (a, b) est un facteur de $\overline{B}(c)$, alors (b,c) est un facteur de $\overline{B}(a)$ et (c, a) de $\overline{B}(b)$.

**1.2.2. Cartes équivalentes :**

Soit $\Omega$ et $\Omega'$ deux cartes cubiques enracinées de genre g et g' respectivement.

Pour toute cellule e de $\Omega$, choisissons un sommet quelconque i(e) de référence sur le bord de e.

Nous définirons que $\Omega$ et $\Omega'$ sont équivalentes si et seulement s'il existe une bijection $\varphi : \Omega \rightarrow \Omega'$ et pour toute cellule e de $\Omega$ un sommet $i(\varphi(e))$ de référence sur le bord de $\varphi(e)$ tel que $B(\varphi(e))$ soit la suite image de B(e) par $\varphi$ : c'est-à-dire que si on note $B(e) = (g_1, g_2, \ldots, g_n)$, avec $n \geq 2$, alors $B(\varphi(e)) = (\varphi(g_1), \varphi(g_2), \ldots, \varphi(g_n))$.

# 2
# Définitions et notations concernant les Monoïdes Libres

## 2.1. Monoïde libre et *Nesiles*

Soit E un ensemble.
L'ensemble des suites finies de E et muni de l'opération de concaténation est un monoïde libre de base E, généralement noté E* et appelé monoïde libre sur E. L'élément neutre est la suite vide, qu'on notera toujours ε.
E est parfois appelé Alphabet et les éléments de E* des Mots.
Ces appellations sont pratiques quand E est un ensemble de caractères.
Dans le cas général, nous appellerons *Nesile* de E tout élément de E*, ε étant le *Nesile* vide.
Pour tout *Nesile* non vide de E, les éléments de E qui constituent le *Nesile* sont appelés les termes du N*esile*.
La longueur d'un *Nesile* ou d'un mot u de E est le nombre de termes du *Nesile*, notée généralement $|u|$ : ε est le seul *Nesile* de longueur 0 et les éléments de E les seuls *Nesiles* de longueur 1.
De même qu'on définit préfixe, suffixe, facteur, ou sous-mot d'un mot, on définit préfixe, suffixe, facteur ou *sous-nesile* d'un *Nesile*.
On dira qu'un mot ou un *Nesile* est *simple* si ses termes sont distincts deux à deux.
Une *Nesilpartition* de E sera un ensemble A de *Nesiles Simples* de E tels que pour tout élément x de E, il existe un et un seul *Nesile Simple* $N \in A$ tel que x soit un terme de N.
D'autre part, toute fonction $f : E \to F$ induit une fonction $f^* : E^* \to F^*$, définie par $f^*(\varepsilon) = \varepsilon$ et
$f^*(x_1, x_2, \ldots, x_n) = (f(x_1), f(x_2), \ldots, f(x_n))$ : on dira que $(f(x_1), f(x_2), \ldots, f(x_n))$ est le *Nesile* image (ou mot image) de $(x_1, x_2, \ldots, x_n)$ par f.

**Exemples :**
Notons AJ = {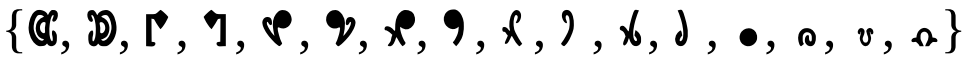} l'alphabet *Jassologique*.
Alors les *Mots Jassologiques* appartiennent à AJ*.
Pour toute carte cubique enracinée Ω de genre g, Ω* est l'ensemble des *Nesiles* de Ω.
On notera ΩAJ l'union de Ω et de AJ.
Alors ΩAJ* est l'ensemble des *Nesiles* de ΩAJ.

**Remarque :**
Soit Ω une carte cubique enracinée de genre g.
Pour toute cellule e de Ω, choisissons un sommet quelconque i(e) de référence sur le bord de e.
Alors B(e) est un *Nesile* de Ω de longueur supérieure ou égale à 2.

## 2.2. Relations d'ordre

Soit E un ensemble muni d'une relation d'ordre (<).
Celle-ci peut s'entendre à E* de multiples manières, dont deux relations d'ordre bien connues : l'ordre lexicographique ($<_{lex}$) et l'ordre des mots croisés ($<_{mc}$) (« radix order » ou « shortlex » en anglais : cf. cours de polytechnique disponible sur Internet).
L'ordre lexicographique, ou ordre du dictionnaire, est défini par $u <_{lex} v$ si et seulement si v est préfixe de u ou bien si u et v peuvent s'écrire u = pau', v = pbv', où p est un mot sur E, a et b des éléments de E, et a < b.

L'ordre des mots croisés est défini par $u <_{mc} v$ si et seulement si $|u| < |v|$ ou bien si $|u| = |v|$ et $u <_{lex} v$.
Par exemple : bar $<_{mc}$ car $<_{mc}$ barda $<_{mc}$ radar $<_{mc}$ abracadabra.

## 2.3. Les *Stratinos*

Les *Stratinos* sont les objets qui permettent de graduer l'*Echelle Jassologique* associée à une carte planaire cubique enracinée.
Soit NJ = {1, 2, 3, 4, …} ∪ {1#, 2#, 3#, 4#, …}, c'est-à-dire l'union des entiers non nuls et des entiers suivis du signe « # », qu'on appellera les entiers *décalés* (on pourrait remplacer # par n'importe quel signe ou par 1/2 ...).
Munissons NJ de la relation d'ordre suivante : 1 < 1# < 2 < 2# < 3 < 3# < 4 < 4# < 5 < 5# < 6 < 6# < 7 <...

L'ensemble des *Stratinos* sera alors l'ensemble NJ* des *Nesiles* de NJ muni de la relation d'ordre lexicographique induite par celle sur NJ.
Pour ne pas confondre avec l'écriture des entiers naturels, les *Stratinos* se présenteront sous forme de suite d'éléments de NJ séparés par des virgules.

**2.3.1. Les *Stratinos* unitaires :**
Un *Stratino* est dit *unitaire* si et seulement si son dernier terme est « 1 ».
Pour tout S*tratino* unitaire Y, il existe un et un seul S*tratino* X tel que Y = (X, 1).
On définit alors l'ensemble D(Y) de la façon suivante :
D(Y) = { (X, n) | n entier naturel ≥ 2} ∪ { (X, n#) | n entier naturel ≥ 1} ∪
{ (X, n, 1) | n entier naturel ≥ 1} ∪ { (X, n#, 1) | n entier naturel ≥ 1}.
Les éléments de D(Y) s'ordonnent (selon la relation d'ordre lexicographique) de la façon suivante :
(X, 1, 1) < (X, 1#) < (X, 1#, 1) < (X, 2) < (X, 2, 1) < (X, 2#) < (X, 2#, 1) < (X, 3) < (X, 3, 1) < (X, 3#) < (X, 3#, 1) < (X, 4) < (X, 4, 1) < (X, 4#) < (X, 4#, 1) < (X, 5) < ...
Pour tout p ≥ 0, notons DNJ(p) l'ensemble des *Stratinos* non unitaires de longueur p et des *Stratinos* unitaires de longueur p+1.
Pour p = 0, DNJ(0) = {ε, (1)}.
Alors pour tout p ≥ 1, l'ensemble des D(Y) pour tous les *Stratinos* unitaires Y de longueur p constitue une partition de DNJ(p).

**2.3.2. Les *Stratinos Naturels* et les *Stratinos Décalés* :**
Un *Stratino* est dit *naturel* (resp. *décalé*) si son dernier terme est un entier naturel (resp. un entier *décalé*).
On définit également le *naturel* et le *décalé* d'un *Stratino* quelconque :
Soit u ∈ NJ* : le naturel et le décalé de u seront notés respectivement (u)~ et (u)#.
Si u = ε , alors (ε)~ = (ε)# = ε.
Sinon u se décompose de façon unique : u = (v, x), où v ∈ NJ* et x ∈ NJ.
Si u est *naturel*, c'est-à-dire si x est un entier naturel, alors (u)~ = u et (u)# = (v, x#).
Si u est *décalé*, c'est-à-dire si x = n#, où n est un entier naturel, alors (u)~ = (v, n) et (u)# = u.

# 3
# Arbres Ordonnés, *Dallajascars* et *Mots Jassologiques Simples*

## 3.1. Arbres Ordonnés

### 3.1.1. Graphe :
Un graphe non orienté (resp. orienté) G = (S, A) est un couple formé d'un ensemble de nœuds S et d'un ensemble A d'arêtes (resp. d'arcs), qui est un ensemble de paires de nœuds (resp. une partie de S×S). Un chemin de s à t est une suite ($s = s_0, \ldots, s_n = t$) de nœuds tels que pour tout i, $1 \leq i \leq n$, ($s_{i-1}$, $s_i$) soit une arête (resp. un arc). Le nœud $s_0$ est l'origine du chemin et le nœud $s_n$ son extrémité. L'entier n est la longueur du chemin. C'est un entier positif ou nul. Un circuit est un chemin de longueur non nulle dont l'origine coïncide avec l'extrémité.
Un chemin est simple si tous ses nœuds sont distincts.
Un graphe est connexe si deux quelconques de ses nœuds sont reliés par un chemin.

### 3.1.2. Arbre libre :
Un arbre libre est un graphe non orienté non vide, connexe et sans circuit.
Soit G = (S, A) un graphe non orienté non vide. Les conditions suivantes sont équivalentes :

1) G est un arbre libre
2) Deux nœuds quelconques de S sont connectés par un chemin simple unique,
3) G est connexe mais ne l'est plus si l'on retire une arête quelconque
4) G est sans circuit, mais ne l'est plus si l'on ajoute une arête quelconque.
5) G est connexe, et Card(A) = Card(S) – 1
6) G est sans circuit, et Card(A) = Card (S) – 1.

### 3.1.3. Arbre enraciné :
Un arbre enraciné est un arbre libre muni d'un nœud distingué, appelé sa racine.
Soit T un arbre de racine r.
Pour tout nœud x, il existe un chemin simple unique de r à x. Tout nœud y sur ce chemin est un ancêtre de x, et x est un descendant de y.
Le sous-arbre de racine x est l'arbre contenant tous les descendants de x.
L'avant-dernier nœud y sur l'unique chemin reliant r à x est le parent (ou le père ou la mère) de x, et x est un enfant (ou un fils ou une fille) de y.
L'arité d'un nœud est le nombre de ses enfants. Un nœud sans enfant est une feuille, un nœud d'arité strictement positive est appelé nœud interne.
La hauteur d'un arbre T est la longueur maximale d'un chemin reliant sa racine à une feuille.
Un arbre réduit à un seul nœud est de hauteur 0.

### 3.1.4. Arbre ordonné :
Un arbre ordonné est un arbre enraciné dans lequel l'ensemble des enfants de chaque nœud est totalement ordonné.
Les enfants d'un nœud d'un arbre ordonné sont souvent représentés par une liste attachée au nœud, c'est-à-dire un *Nesile* de S.
Soit T = (S, A) un arbre ordonné et r sa racine.
On définit alors les fonctions $H : S \to S^*$ et $\varphi : S \to S \cup \{\varepsilon\}$ de la manière suivante : pour tout nœud x, H(x) est la suite ordonnée des enfants de x, avec $H(x) = \varepsilon$ si x est une feuille, et $\varphi(x)$ le parent de x, avec $\varphi(r) = \varepsilon$.
On dira que ce sont les *fonctions de ramification* de T, H étant la *fonction de ramification filiale* et φ la *fonction de ramification parentale*. En particulier, $\{H(x) \mid x \in S\}$ constitue une *Nesilpartition* de $S \backslash \{r\}$.

## 3.2. *Dallajascars*

### 3.2.1. Définition :
Soit S un ensemble fini non vide.
Un *Dallajascar* (φ,H) sur S est un couple de fonctions $\varphi : S \to S \cup \{\varepsilon\}$ et $H : S \to S^*$ tels qu'il existe un élément r de S, appelé racine, vérifiant les 4 propriétés suivantes :

i) $\varphi(r) = \varepsilon$,
ii) pour tout x différent de r, $\varphi(x) \neq \varepsilon$,

iii) pour tout x différent de r et tout entier n ≥ 1, $\varphi^n(x) \neq x$,

iv) pour tout x, H(x) est un *Nesile Simple* de S dont les termes constituent l'ensemble des antécédents de x par φ, avec H(x) = ε si x n'a pas d'antécédent.

Le couple T = (S, A), où A est l'ensemble des paires {x, φ(x)} pour tout x de S différent de r, constitue alors un arbre ordonné dont r est la racine. Les fonctions φ et H sont exactement les *fonctions de ramification* de T définies ci-dessus, H étant la *fonction de ramification filiale* et φ la *fonction de ramification parentale*.

**3.2.2. Les *Relations d'Emboîtement* :**

Soit (φ,H) un *Dallajascar* sur un ensemble fini non vide S, et notons r la racine.

1ère *Relation d'Emboîtement*, notée « ⌠ » :

Pour tout x et y, on dira que x *emboîte* y, et on écrira « x ⌠ y », si et seulement s'il existe un entier n ≥ 1 tel que $\varphi^n(y) = x$, avec x = y si et seulement si n = 0.

On vérifie que si x *emboîte*, y, alors y n'*emboîte* pas x.

2ème *Relation d'Emboîtement*, notée « ¬ » :

Pour tout x et y, on dira que x *précède* y, et on écrira « x ¬ y », si et seulement s'il existe deux entiers n ≥ 1 et m ≥ 1 et a ∈ S tels que $\varphi^n(x) = \varphi^n(y) = a$, et tels que $\varphi^{n-1}(x)$ soit situé avant $\varphi^{n-1}(y)$ dans H(a), auquel cas a et les entiers n et m sont nécessairement uniques.

On pourra dire éventuellement que « x ¬ y sous a ».

On vérifie que si x *précède* y, alors y ne *précède* pas x.

De plus, si x *précède* y, alors x n'*emboîte* pas y et y n'*emboîte* pas x, et symétriquement, si x *emboîte* y, alors x ne *précède* pas y et y ne *précède* pas x.

Quelques Propriétés :

1) Si x ⌠ y et y ⌠ z, alors x ⌠ z.
2) Si x ¬ y et y ¬ z, alors x ¬ z.
3) Si x ¬ y et y ⌠ z, alors x ¬ z.
4) Si x ⌠ y et y ¬ z, alors x ⌠ z ou bien x ¬ z.
5) Si x ⌠ z et y ⌠ z, alors x ⌠ y ou bien y ⌠ x.
6) Si x ¬ z et y ⌠ z, alors x ¬ y ou bien y ⌠ x.
7) Si x ⌠ y et x ¬ z, alors y ¬ z.
8) Si x ⌠ y, x' ⌠ y' et x ¬ x', alors y ¬ y'.
9) Si x ⌠ y, x' ⌠ y' et y ¬ y', alors x ⌠ x' ou bien x ¬ x'.

Relation d'ordre sur S induite par (φ,H), notée « $<_{(\varphi,H)}$ « : pour tout x et y, x $<_{(\varphi,H)}$ y si et seulement si x ⌠ y ou x ¬ y.

Les propriétés précédentes induisent que :

- pour tout x et y différents, x $<_{(\varphi,H)}$ y ou bien y $<_{(\varphi,H)}$ x.
- si x $<_{(\varphi,H)}$ y et y $<_{(\varphi,H)}$ z, alors x $<_{(\varphi,H)}$ z.

**3.2.3. Equivalence :**

Soit (φ,H) et (θ,K) deux *Dallajascars* sur des ensembles S et T respectifs.

On définit que (φ,H) et (θ,K) sont équivalents s'il existe une bijection f : S → T telle que pour tout élément a ∈ S, K(f(a)) soit le *Nesile* image de H(a) par f.

**3.3. Les *Mots Jassologiques Simples***

**3.3.1. La R*ègle d'Emboîtement* :**

Soit J = {⸢, ⸣}, où « ⸢ » est appelé caractère ouvrant et « ⸣ » caractère fermant.

Et soit M = $x_1\ x_2 \ldots x_r$ un mot sur J contenant exactement N caractères ouvrants, avec 1 ≤ N ≤ r, où r est sa longueur.

La R*ègle d'Emboîtement* se définit par récurrence sur N.

Si N = 1, alors M vérifie la *Règle d'Emboîtement* si et seulement si M = ⸢⸣.

Supposons que N ≥ 2 et notons $\alpha_N$ le plus grand entier, 1 ≤ $\alpha_N$ ≤ r, tel que $x_{\alpha N}$ = ⸢.

Alors M vérifie la *règle d'emboîtement* si et seulement si

i) $\alpha_N \neq r$, ce qui induit que $x_{\alpha N+1}$ = ⸣.

ii) le mot obtenu en supprimant le facteur $x_{\alpha N}\ x_{\alpha N+1}$, qui est un mot sur J contenant exactement N-1 caractères ouvrants, vérifie la *Règle d'Emboîtement*. On dira que $x_{\alpha N+1}$ est le caractère fermant correspondant à $x_{\alpha N}$ dans M.

Ce qui induit que si on note $\alpha_1, \alpha_2, \ldots, \alpha_N$, les entiers tels que $x_{\alpha k}$ = ⸢ pour tout k, $1 \le k \le N$, avec $1 \le \alpha_1 < \alpha_2 < \ldots < \alpha_N \le r$, alors pour tout k, il existe un et un seul entier $\beta_k$, tel que $x_{\beta k}$ soit le caractère fermant correspondant à $x_{\alpha k}$ dans M.
De plus, ils vérifient la propriété suivante : $\alpha_k < \beta_k$ et si $k < N$, $\beta_k < \alpha_{k+1}$ ou bien $\beta_{k+1} < \beta_k$.
Nécessairement, on aura toujours $\alpha_1 = 1$ et $r = 2N$.

### 3.3.2. Les *Mots Jassologiques Simples* :

Soit M un mot sur J = {⸢, ⸣}.
On dira que M est un *Mot Jassologique Simple* si et seulement si M vérifie la *Règle d'Emboîtement* et si le caractère fermant correspondant au premier caractère ouvrant est le dernier caractère de M.
Si on reprend les notations précédentes, cela signifie que $\alpha_1 = 1$ et $\beta_1 = r$.

### 3.3.3. Le *Mot Jassologique Simple* correspondant à un *Dallajascar* :

Soit S un ensemble fini non vide et $(\varphi,H)$ un *Dallajascar* sur S. Notons r sa racine.
Pour tout $x \in S$, définissons $LH(x) \in (S \cup J)^*$ de la façon suivante :

- si $H(x) = \varepsilon$, alors LH(x) = ⸢⸣,
- sinon notons $H(x) = (x_1, x_2, \ldots, x_r)$, avec $r \ge 1$, alors LH(x) = ⸢ $x_1\ x_2 \ldots x_r$ ⸣.

Le Mot sur J correspondant à $(\varphi,H)$, noté $M(\varphi,H)$ se construit alors de la façon suivante :
Pour $t = 0$, $M_0 = LH(r)$.
Si $M_0 \in J^*$, alors c'est terminé : $M(\varphi,H) = M_0$.
Sinon notons $x_1$ le premier terme de $M_0$ qui soit un élément de S.
$M_1$ s'obtient alors en remplaçant $x_1$ par $LH(x_1)$ dans $M_0$.
Si $M_1 \in J^*$, alors c'est terminé : $M(\varphi,H) = M_1$.
Sinon notons $x_2$ le premier terme de $M_1$ qui soit un élément de S.
$M_2$ s'obtient alors en remplaçant $x_2$ par $LH(x_2)$ dans $M_1$.
Si $M_2 \in J^*$, alors c'est terminé : $M(\varphi,H) = M_2$.
Sinon notons $x_3$ le premier terme de $M_2$ qui soit un élément de S.
$M_3$ s'obtient alors en remplaçant $x_3$ par $LH(x_3)$ dans $M_2$.
Et ainsi de suite.
Il existe nécessairement un et un seul entier $m \ge 1$ tel que $M(\varphi,H) = M_m$, vérifiant $M_m \in J^*$ et $M_{m-1} \notin J^*$. On vérifie alors les propriétés suivantes :

- $M(\varphi,H)$ est un *Mot Jassologique Simple*.
- m est égal au nombre d'éléments de S et au nombre de caractères ouvrants dans $M(\varphi,H)$.

### 3.3.4. Théorème fondamental :

Soit $(\varphi,H)$ et $(\theta,K)$ deux *Dallajascars*.
Alors $(\varphi,H)$ et $(\theta,K)$ sont équivalents si et seulement si $M(\varphi,H) = M(\theta,K)$.

### 3.3.5. Le *Dallajascar* correspondant à un *Mot Jassologique Simple* :

Soit $M = x_1\ x_2 \ldots x_r$ un *Mot Jassologique Simple* et notons m le nombre de caractères ouvrants, avec $1 \le m \le r$, où r est sa longueur.
Notons $\alpha_1, \alpha_2, \ldots, \alpha_m$, les entiers correspondant aux caractères ouvrants, avec $1 \le \alpha_1 < \alpha_2 < \ldots < \alpha_m \le r$,et pour tout k, $1 \le k \le m$, notons $\beta_k$ l'entier tel que $x_{\beta k}$ = ⸣ soit le caractère fermant correspondant à $x_{\alpha k}$ = ⸢ dans M.
Posons $a_k = (\alpha_k, \beta_k)$, $S = \{a_1, a_2, \ldots, a_m\}$, et définissons les fonctions $\varphi : S \to S \cup \{\varepsilon\}$ et $H : S \to S^*$ de la façon suivante :

- $\varphi(a_1) = \varepsilon$,
- si $m \ge 2$, pour tout k, $2 \le k \le m$, $\varphi(a_k) = a_q$, où q est le plus grand entier, $1 \le q \le p - 1$, tel que $\beta_k < \beta_q$. Comme $\beta_1 = r$, $\beta_q$ est bien défini.
- pour tout k, $1 \le k \le m$, si $a_k$ n'a pas d'antécédent pas $\varphi$, alors $H(a_k) = \varepsilon$, sinon $H(a_k) = (a_{k1}, a_{k2}, \ldots, a_{kp})$, où $\{a_{k1}, a_{k2}, \ldots, a_{kp}\}$ constitue l'ensemble des antécédents de $a_k$ par $\varphi$, avec $1 \le k_1 < k_2 < \ldots < \alpha_m \le m$.

On vérifie que $(\varphi,H)$ est bien un *Dallajascar* sur S, $a_1$ étant la racine.
De plus, on montre que $M(\varphi,H) = M$.

### 3.3.6. Corollaire fondamental :

On obtient une bijection entre l'ensemble des *Mots Jassologiques Simples* et celui des classes d'équivalence de *Dallajascars*.

# 4
# Le processus *Jassologique*

## 4.1. Base de données

Soit $\Omega$ une carte planaire cubique enracinée.
Notons $\lambda$ l'arête orientée distinguée : on a défini les faces racines correspondantes. Notons $w_0$ la face racine (-) et $w_1$ la face racine (+). Puis notons $i(w_1)$ et $i(w_0)$ les extrémités de $\lambda$, telles que $\lambda$ soit orientée de $i(w_1)$ vers $i(w_0)$.
Alors le premier terme de la bordure de $w_0$ à partir de $i(w_0)$, notée $B(w_0)$, est $w_1$, et le premier terme de la bordure de $w_1$ à partir de $i(w_1)$, notée $B(w_1)$, est $w_0$.
Pour toute cellule e de $\Omega$ différente de $w_0$ et $w_1$, on choisira au hasard un sommet de référence i(e) sur son bord, à partir duquel on définit la bordure B(e) et la bordure étendue $\overline{B}(e)$.
L'ensemble des bordures étendues ainsi fixées constituera la base de données du processus algorithmique que nous allons décrire ci-dessous.
L'objectif final est d'écrire le *Mot Jassologique* correspondant, noté $MJ(\Omega)$, qui est un mot sur l'*Alphabet Jassologique* AJ = {[illegible], [illegible], [illegible], [illegible], [illegible], [illegible], [illegible], [illegible], [illegible], [illegible], [illegible], [illegible], ●, [illegible], [illegible], [illegible]}.

D'autre par, nous définirons un certain nombre de propriétés caractéristiques des cartes planaires cubiques, qui ne sont plus valables pour les cartes non planaires, c'est-à-dire de genre $g \geq 1$.

**P1 : 1ère propriété caractéristique des cartes planaires cubiques :**
Soit e une cellule quelconque de $\Omega$, et notons $\overline{B}(e) = (g_1, g_2, \ldots, g_n, g_{n+1})$, avec $g_{n+1} = g_1$ et $n \geq 2$.
Alors pour tout k et r, $1 \leq k \leq n$ et $1 \leq r \leq n$, tels que $k \neq r$, si $g_k = g_r$, alors $g_{k+1} \neq g_{r+1}$.

## 4.2. *Jasses, Rovéjasses* et les fonctions relatives aux *Rovéjasses*

### 4.2.1. Les *Jasses* de $\Omega$ :
Posons $j_0 = \{w_0\}$ et $J_0 = \Omega \setminus j_0$.
La *jasse* $j_1$ sera alors l'ensemble des cellules de $J_0$ qui sont en contact avec $w_0$.
Puis posons $J_1 = J_0 \setminus j_1 = \Omega \setminus (j_0 \cup j_1)$.
Si $J_1 = \emptyset$, alors on posera $j_p = \emptyset$ pour tout $p \geq 2$.
Sinon notons $j_2$ l'ensemble des cellules de $J_1$ qui sont en contact avec au moins une cellule de $j_1$.
Puis posons $J_2 = J_1 \setminus j_2 = \Omega \setminus (j_0 \cup j_1 \cup j_2)$.
Si $J_2 = \emptyset$, alors on posera $j_p = \emptyset$ pour tout $p \geq 3$.
Sinon notons $j_3$ l'ensemble des cellules de $J_2$ qui sont en contact avec au moins une cellule de $j_2$.
Puis posons $J_3 = J_2 \setminus j_1 = \Omega \setminus (j_0 \cup j_1 \cup j_2 \cup j_3)$.
Et ainsi de suite.
Il existe un et un seul entier $m \geq 1$ tel que $j_{m+1} = \emptyset$ et $j_m \neq \emptyset$.
Les *jasses* $j_0, j_1, \ldots, j_m$, constituent une partition de $\Omega$.

### 4.2.2. Les *Rovéjasses* de $\Omega$ :
Pour tout p, $1 \leq p \leq m$, notons $V(\Omega,p)$ l'ensemble des composantes connexes de $j_p$, appelées *Rovéjasses* de $j_p$.
Elles constituent une partition de $j_p$.
Ce qui signifie que si $j_p$ est connexe, alors $V(\Omega,p) = \{j_p\}$.
C'est notamment toujours le cas pour $j_1$ : $V(\Omega,1) = \{j_1\}$.
Puis on notera $V(\Omega) = \cup_{1 \leq p \leq m} V(\Omega,p)$ l'ensemble des *Rovéjasses* de $\Omega$, qui constituent une partition de $\Omega$.

### 4.2.3. La fonction *Périmètre de contour*, notée $\pi$ :
La fonction *Périmètre de contour*, notée $\pi : V(\Omega) \rightarrow P(\Omega)$, est définie de la façon suivante :
Pour tout p, $1 \leq p \leq m$, et toute Rovéjasse $v \in V(\Omega,p)$, $\pi(v)$ sera l'ensemble des cellules de $j_{p-1}$ qui sont en contact avec au moins une cellule de v.
Il est évident que $\pi(j_1) = j_0 = \{w_0\}$. Si $p \geq 2$, on vérifie facilement que $\pi(v)$ contient au moins deux cellules.

**P2 : 2ème propriété caractéristique des cartes planaires cubiques :**
Pour tout p, $2 \le p \le m$, et toute *Rovéjasse* $v \in V(\Omega,p)$, $\pi(v)$ est connexe et est donc contenue dans une et une seule *Rovéjasse* de $j_{p-1}$ : on dit que v est *orthogonale* à cette *Rovéjasse* de $j_{p-1}$.

**P3 : 3ème propriété caractéristique des cartes planaires cubiques :**
Pour tout p, $2 \le p \le m$, toute *Rovéjasse* $v \in V(\Omega,p)$, et toute cellule $a \in \pi(v)$, les termes de $\overline{B}(a)$ qui appartiennent à v constituent un facteur de $\overline{B}(a)$ ou bien se divisent en deux facteurs de $\overline{B}(a)$ dont l'un est préfixe et l'autre suffixe.
Autrement dit, soit $\overline{B}(a) = (P, X, S)$, où $X \in v^*$ et P et $S \in (\Omega\backslash v)^*$,
soit $\overline{B}(a) = (X_1, M, X_2)$, où $X_1$ et $X_2 \in v^*$ et $M \in (\Omega\backslash v)^*$.
Il est évident que P, X et S (resp. $X_1$, M, $X_2$) sont tous différents de ε.
Ce qui nous permet de définir les fonctions $c(v)^+ : \pi(v) \rightarrow \pi(v)$ et $c(v)^- : \pi(v) \rightarrow \pi(v)$, de la façon suivante :
Si $\overline{B}(a) = (P, X, S)$, où $X \in v^*$ et P et $S \in (\Omega\backslash v)^*$, alors $c(v)^+(a)$ sera le dernier terme de P et $c(v)^-(a)$ le premier terme de S. Si $\overline{B}(a) = (X_1, M, X_2)$, où $X_1$ et $X_2 \in v^*$ et $M \in (\Omega\backslash v)^*$, alors $c(v)^+(a)$ sera le dernier terme de M et $c(v)^-(a)$ le premier terme de M.
Il est évident que dans tous les cas, $c(v)^+(a)$ et $c(v)^-(a)$ appartiennent à $\pi(v)$.
En particulier, il existe une et une seule cellule $e \in v$ et une et une seule cellule $g \in v$ telles que $(c(v)^+(a), a)$ soit un facteur de $\overline{B}(e)$ et $(a, c(v)^-(a))$ un facteur de $\overline{B}(g)$.
Remarquons que $c(v)^+(a)$ et $c(v)^-(a)$ ne dépendent pas du choix de i(a) pour $\overline{B}(a)$.
Inversement, pour toute cellule $e \in v$ et tout couple $(a, b) \in \pi(v)^*$ et facteur de $\overline{B}(e)$, alors $a = c(v)^+(b)$ et $b = c(v)^-(a)$.

**4.2.4. La fonction *baou* et la fonction *zouc* :**
La fonction baou : $\cup_{2\le p\le m}V(\Omega,p) \rightarrow \Omega^*$, et la fonction zouc : $V(\Omega) \rightarrow \Omega$, seront définies par récurrence (voir ci-dessous), et devront satisfaire les propriétés suivantes :
Pour tout p, $1 \le p \le m$, et toute *Rovéjasse* $v \in V(\Omega,p)$, $zouc(v) \in v$, avec $zouc(j_1) = w_1$ dans le cas où p = 1.
Si $p \ge 2$, baou(v) sera un couple $(a, b) \in \pi(v)^*$ tel que $c(v)^+(a) = b$ (et donc $c(v)^-(b) = a$).

**4.2.5. La fonction m*édiane* et la fonction *binôme* :**
La foncion médiane : $\Omega \rightarrow \{oui, non\}$, et la fonction binôme : $\Omega \rightarrow AJ^*$, seront définies par récurrence (voir ci-dessous), et devront satisfaire les propriétés suivantes :
Pour toute cellule $e \in \Omega$, binôme(e) = 𝕮𝕯 ou ⸢⸣ ou ❛❜ ou ❝❞ ou ⸌⸍ ou ⸜⸝ .
En particulier, binôme($w_0$) = 𝕮𝕯 et binôme($w_1$) = ⸢⸣, avec médiane($w_0$) = médiane($w_1$) = non.
Pour tout p, $2 \le p \le m$, et toute *Rovéjasse* $v \in V(\Omega,p)$, binôme(zouc(v)) = ❛❜ ou ❝❞ .
En particulier, si binôme(zouc(v)) = ❛❜, alors médiane(zouc(v)) = non.
Et si $\pi(v)$ ne contient que 2 cellules, alors également médiane(zouc(v)) = non.
Autrement dit, si médiane(zouc(v)) = oui, alors nécessairement binôme(zouc(v)) = ❝❞ et $\pi(v)$ contient au moins 3 cellules.
Pour toute cellule e de v autre que zouc(v), binôme(e) = ⸌⸍ ou ⸜⸝ , avec médiane(e) = oui ou non.

**4.2.6. La fonction c*aouly* et la fonction *fan* :**
La fonction caouly : $\cup_{2\le p\le m}V(\Omega,p) \rightarrow \Omega^*$, et la fonction fan : $\cup_{2\le p\le m}V(\Omega,p) \rightarrow \Omega^*$, seront définies automatiquement à partir des fonctions *baou* et *médiane* de la manière suivante :
Soit p, $2 \le p \le m$, et $v \in V(\Omega,p)$. Et notons baou(v) = (a, b).
Par définition, $(a, b) \in \pi(v)^*$, avec avec $c(v)^+(a) = b$.
La propriété P3 induit que $c(v)^+(b) = a$ si et seulement si $\pi(v) = \{a, b\}$.
Si $c(v)^+(b) = a$, alors caouly(v) = (b, a) et fan(v) = ε.
Sinon posons $g_1 = c(v)^+(b)$, $g_2 = c(v)^+(g_1)$, $g_3 = c(v)^+(g_2)$, et ainsi de suite pour tout $k \ge 2$, $g_k = c(v)^+(g_{k-1})$.
La propriété P3 induit qu'il existe un et un seul entier $r \ge 1$ tel que $g_{r+1} = a$ et $g_1, g_2, \ldots, g_r$, différentes de a.
De plus, les cellules $g_1, g_2, \ldots, g_r$ sont nécessairement distinctes et différentes de a.
Ce qui induit que $\pi(v) = \{a, b, g_1, g_2, \ldots, g_r\}$.
Autrement dit, $(a, b, g_1, g_2, \ldots, g_r)$ constitue un anneau autour de v.
Si médiane(zouc(v)) = non, alors caouly(v) = $(g_r, a)$ et fan(v) = $(g_1, g_2, \ldots, g_r)$.
Si médiane(zouc(v)) = oui, alors caouly(v) = $(g_{r-1}, g_r)$ et fan(v) = $(g_1, g_2, \ldots, g_{r-1})$, avec caouly(v) = $(b, g_1)$ et fan(v) = ε si r = 1.

### 4.3. Les fonctions de ramification H et φ

La fonction de ramification H : $\Omega \to \Omega AJ^*$, où $\Omega AJ = \Omega \cup AJ$, est l'étape ultime avant de pouvoir déterminer MJ(Ω).
Un certain nombre d'autres fonctions seront nécessaires à définir pour arriver à la fonction H.
Elles seront toutes déterminées par récurrence (voir ci-dessous).
Il en est de même pour l'autre fonction de ramification $\varphi : \Omega \to \Omega \cup \{\varepsilon\}$.

#### 4.3.1. Les fonctions gh et $\overline{\overline{B}}°$ :

La fonction gh : $\Omega \to \Omega^*$ : pour toute cellule $e \in \Omega$, gh(e) est un couple de cellules facteur de $\overline{B}(e)$ qui se détermine par récurrence (voir ci-dessous).
La fonction $\overline{\overline{B}}°$ : $\Omega \to \Omega^*$ se déduit alors des fonctions gh et $\overline{B}$ de la manière suivante :
Soit e une cellule de Ω et notons $\overline{B}(e) = (g_1, g_2, \ldots, g_r)$, avec $r \geq 3$ et $g_r = g_1$ par définition.
Alors il existe un et un seul entier k, $1 \leq k \leq r-1$, tel que $gh(e) = (g_k, g_{k+1})$.
Si k = 1, alors $\overline{\overline{B}}°(e) = \overline{B}(e)$.
Si $k \geq 2$, alors $\overline{\overline{B}}°(e) = (g_k, g_{k+1}, \ldots, g_{r-1}, g_1, g_2, \ldots, g_k)$ :
Autrement dit, $\overline{\overline{B}}°(e)$ est une rotation de $\overline{B}(e)$ afin que gh(e) soit un préfixe de $\overline{\overline{B}}°(e)$.

#### 4.3.2. Les fonctions dh et $\overline{\overline{Br}}$ :

La fonction dh : $\Omega \backslash \{w_0\} \to \Omega^*$ : pour toute cellule $e \in \Omega \backslash \{w_0\}$, dh(e) est un couple de cellules facteur de $\overline{B}(e)$, et donc aussi de $\overline{\overline{B}}°(e)$, qui se détermine par récurrence (voir ci-dessous) et qui sera différent de gh(e).

La fonction $\overline{\overline{Br}}$ : $\Omega \backslash \{w_0\} \to \Omega^*$ se déduit alors des fonctions dh et $\overline{\overline{B}}°$ de la manière suivante :
Soit e une cellule de $\Omega \backslash \{w_0\}$ et notons $\overline{\overline{B}}°(e) = (g_1, g_2, \ldots, g_r)$, avec $r \geq 3$ et $g_r = g_1$ par définition.
Alors il existe un et un seul entier k, $3 \leq k \leq r$, tel que $dh(e) = (g_{k-1}, g_k)$.
Définissons $\overline{\overline{Br}}(e) = (g_1, g_2, \ldots, g_k)$.
Autrement dit, $\overline{\overline{Br}}(e)$ est le préfixe de $\overline{\overline{B}}°(e)$ qui a pour propriété que gh(e) et dh(e) soient respectivement les couples préfixe et suffixe de $\overline{\overline{Br}}(e)$.

#### 4.3.3. Les fonctions ν, ξg, ξd, $\xi\Delta_k$ et $\xi R_k$ :

Les fonctions $\nu : \Omega \backslash \{w_0\} \to \{0, 1, 2, \ldots\}$ et les fonctions : ξg, ξd, $\xi\Delta_k$ et $\xi R_k : \Omega \backslash \{w_0\} \to \Omega^*$, pour tout entier $k \geq 0$, se déterminent par récurrence (voir ci-dessous) et vérifieront la propriété suivante : pour toute cellule e de $\Omega \backslash \{w_0\}$, ξg(e), ξd(e), $\xi\Delta_k(e)$ pour tout k, $0 \leq k \leq \nu(e)$, et $\xi R_k(e)$ pour tout k, $1 \leq k \leq \nu(e)$ et dans le cas où $\nu(e) \geq 1$, seront des facteurs de $\overline{\overline{Br}}(e)$ vérifiant :
$\overline{\overline{Br}}(e) = (\xi g(e), \xi\Delta_0(e), \xi R_1(e), \xi\Delta_1(e), \xi R_2(e), \xi\Delta_2(e), \ldots, \xi R_{\nu(e)}(e), \xi\Delta_{\nu(e)}(e), \xi d(e))$, avec
$\overline{\overline{Br}}(e) = (\xi g(e), \xi\Delta_0(e), \xi d(e))$ si $\nu(e) = 0$.
D'autre part, on posera $\xi R_0(e) = \varepsilon$ et pour tout $k \geq \nu(e)$, $\xi\Delta_k(e) = \xi R_k(e) = \varepsilon$.

#### 4.3.4. Les fonctions ξRg, ξRd, ξG et ξH :

Les fonctions ξRg, ξRd, ξG et ξH : $\Omega \backslash \{w_0\} \to \Omega^*$ se déterminent par récurrence (voir ci-dessous) et vérifieront la propriété suivante : pour toute cellule e de $\Omega \backslash \{w_0\}$, ξRg(e) est un suffixe de ξg(e), ξRd(e) est un préfixe de ξg(e) et ξG(e) un facteur de ξg(e) qui suit ξRd(e) dans ξg(e). ξH(e) sera alors un facteur de $\overline{\overline{Br}}(e)$ défini par :
$\xi H(e) = (\xi Rg(e), \xi\Delta_0(e), \xi R_1(e), \xi\Delta_1(e), \xi R_2(e), \xi\Delta_2(e), \ldots, \xi R_{\nu(e)}(e), \xi\Delta_{\nu(e)}(e), \xi Rd(e), \xi G(e))$, avec
$\xi H(e) = (\xi Rg(e), \xi\Delta_0(e), \xi Rd(e), \xi G(e))$ si $\nu(e) = 0$.

#### 4.3.5. Les fonctions Rg, Rd, G, $\Delta_k$, $R_k$, et H :

Les fonctions Rg, Rd, G, $\Delta_k$, $R_k$ : $\Omega \backslash \{w_0\} \to$ (Ω∪{•, ɷ, b, ɷ})* se déduisent des fonctions ξRg, ξRd, ξG, $L\Delta_k$, $\xi R_k$ par récurrence (voir ci-dessous).
La fonction H : $\Omega \to$ (Ω∪{•, ɷ, b, ɷ})* se définira alors de la manière suivante :
$H(w_0) = (w_1)$, et pour toute cellule e de $\Omega \backslash \{w_0\}$,
$H(e) = (Rg(e), \Delta_0(e), R_1(e), \Delta_1(e), R_2(e), \Delta_2(e), \ldots, R_{\nu(e)}(e), \Delta_{\nu(e)}(e), Rd(e), G(e))$, avec
$H(e) = (Rg(e), \Delta_0(e), Rd(e), G(e))$ si $\nu(e) = 0$.

### 4.4. La fonction *Stratojasse* Σ et la coloration *Stratojasse*

On a défini précédemment l'ensemble des *Stratinos*, noté NJ*, muni de la relation d'ordre lexicographique induite par celle sur NJ : 1 < 1# < 2 < 2# < 3 < 3# < 4 < 4# < 5 < 5# < 6 < 6# < 7 <...

#### 4.4.1. La fonction *Stratojasse* Σ et l'ensemble NJ(Ω) des *Stratinos* de Ω:

La fonction Stratojasse $\Sigma : NJ^* \to \Omega AJ^*$ se définira par récurrence (voir ci-desous).

L'ensemble $NJ(\Omega) = \{ X \in NJ^* \mid \Sigma(X) \neq \varepsilon \}$ constituera l'*Échelle Jassologique* associée à la carte.

La fonction Σ vérifiera les propriétés suivantes :

i) Pour toute cellule e de Ω, il existera un et un seul *Stratino* $X \in NJ(\Omega)$ et un et un seul entier $r \geq 1$ tel que e soit le $r^{ème}$ terme de Σ(X),

ii) $\Sigma(\varepsilon) = (w_0)$ et $\Sigma(1) = (w_1)$,

iii) Pour tout *Stratino* unitaire $Y \in NJ(\Omega)$ de longueur $p \geq 1$, et pour toute cellule e appartenant à Σ(Y), il existera une et une seule *Rovéjasse* v de $j_p$ telle que e = zouc(v). Ce qui induira que les cellules de Σ(Y) seront deux-à-deux disjointes (pas en contact). On notera X le *Stratino* tel que Y = (X, 1). Si X est *naturel*, alors binôme(e) = [illegible]. Si X est *décalé*, alors binôme(e) = [illegible].

iv) Pour tout *Stratino* non unitaire $X \in NJ(\Omega)$ de longueur $p \geq 1$, et pour toute cellule e appartenant à Σ(X), e appartiendra à une *Rovéjasse* v de $j_p$ et e ≠ zouc(v). Ce qui induira que binôme(e) = [illegible] ou [illegible] . De plus, e sera en contact avec au plus deux cellules de Σ(X), celle qui la précède dans Σ(X) et celle qui la suit. C'est pourquoi on dit que Σ(X) se factorise en chaînes de cellules. En particulier, si g est la première cellule qui suit e dans Σ(X), alors e et g seront en contact si et seulement si elles ne sont pas séparées par le caractère [illegible].

v) Pour tout *Stratino décalé* $X \in NJ(\Omega)$ et différent de ε, alors $\Sigma(X) \in AJ^*$ et même plus précisément, Σ(X) sera une puissance du mot ([illegible]) : il existera un entier $r \geq 1$ tel que $\Sigma(X) = ([illegible])^r$ = [illegible] [illegible] [illegible] … [illegible] (r fois).

#### 4.4.2. Les fonctions Tg, Td, $T_k$, et T :

Les fonctions Tg, Td, $T_k : \Omega\backslash\{w_0\} \to \Omega AJ^*$ se définissent parallèlement aux fonctions Rg, Rd, $R_k$ vues précédemment.

La fonction $T : \Omega\backslash\{w_0\} \to \Omega AJ^*$ se définira alors de la manière suivante :

Pour toute cellule e de $\Omega\backslash\{w_0\}$,

$T(e) = (Tg(e), T_1(e), T_2(e), ..., T_{\nu(e)}(e), Td(e))$, avec

$T(e) = (Tg(e), Td(e))$ si $\nu(e) = 0$.

#### 4.4.3. Les fonctions Sg, Sd, $S_k$, et S :

Les fonctions Sg, Sd, $S_k : \Omega\backslash\{w_0\} \to \Omega AJ^*$ se définissent parallèlement aux fonctions Rg, Rd, $R_k$ vues précédemment.

La fonction $S : \Omega\backslash\{w_0\} \to \Omega AJ^*$ se définira alors de la manière suivante :

Pour toute cellule e de $\Omega\backslash\{w_0\}$,

$S(e) = (Sg(e), S_1(e), S_2(e), ..., S_{\nu(e)}(e), Sd(e))$, avec

$S(e) = (Sg(e), Sd(e))$ si $\nu(e) = 0$.

#### 4.4.4. La fonction ΔG :

La fonction $\Delta G : \Omega\backslash\{w_0\} \to AJ^*$ se définit parallèlement à la fonction G vue précédemment.

Pour toute cellule e de $\Omega\backslash\{w_0\}$, notons r la longueur de G(e). Alors $\Delta G(e) = ([illegible])^r$, avec $\Delta G(e) = \varepsilon$ si r = 0.

#### 4.4.5. La fonction sΣ :

La fonction $s\Sigma : NJ^* \to \Omega^*$ se déduira de la fonction Σ de la façon suivante :

Pour tout *Stratino* $X \in NJ(\Omega)$, sΣ(X) sera par définition le sous-nesile de Σ(X) constitué des termes de Σ(X) qui sont des cellules de Ω. Ce qui induira que si Σ(X) ne contient que des *caractères Jassologiques*, alors $s\Sigma(X) = \varepsilon$.

Cette fonction sera utile pour déterminer Σ(X) par récurrence (voir ci-dessous) :

Soit $X \in NJ(\Omega)$, tel que $s\Sigma(X) \neq \varepsilon$, et notons $s\Sigma(X) = (e_1, e_2, \dots, e_r)$, avec $r \geq 1$.

Ce qui induira que X est un *stratino naturel* : X = (Z, n), où Z est un *Stratino* et n un entier.

Alors Σ(Z, n, 1), Σ(Z, n#), Σ(Z, n#, 1), et Σ(Z, n+1) se déduiront de la façon suivante :

$\Sigma(Z, n, 1) = (G(e_1), G(e_2), \dots, G(e_r))$,

$\Sigma(Z, n\#) = (\Delta G(e_1), \Delta G(e_2), \dots, \Delta G(e_r))$,

$\Sigma(Z, n\#, 1) = (T(e_1), T(e_2), \dots, T(e_r))$,

$\Sigma(Z, n+1) = (S(e_1), S(e_2), \ldots, S(e_r))$.

**4.4.6. Coloration Stratojasse :**

Soit $c_0$, $c_1$, $c_2$ et $c_3$ quatre couleurs différentes.
Il s'agit de colorer les cellules de chaque *Stratojasse* d'une des 4 couleurs.
Soit $X \in NJ(\Omega)$ de longueur $p \geq 0$, tel que $s\Sigma(X) \neq \varepsilon$ :
Si $X = \varepsilon$, alors $s\Sigma(X) = (w_0)$ par définition : $w_0$ sera colorée en $c_0$.
Sinon on notera n le dernier terme de X, qui est nécessairement un entier.
Si p est impair et n impair, alors toutes les cellules de $s\Sigma(X)$ seront colorées en $c_1$.
Si p est impair et n pair, alors les cellules de $s\Sigma(X)$ seront colorées en $c_3$.
Si p est pair et n impair, alors les cellules de $s\Sigma(X)$ seront colorées en $c_0$.
Si p est pair et n pair, alors les cellules de $s\Sigma(X)$ seront colorées en $c_2$.
En particulier, la coloration *Stratojasse* vérifiera la propriété suivante :
Pour toutes cellules e et g de même couleur, e et g seront en contact si et seulement si elles appartiennent à une même *Stratojasse* $\Sigma(X)$, l'une précédant l'autre dans $s\Sigma(X)$ et sans être séparées par le caractère ɷ dans $\Sigma(X)$.

## 4.5. Déroulement de l'algorithme

**4.5.1. Les ensembles DNJ(Ω, p) :**

$NJ(\Omega) = \{ X \in NJ^* \mid \Sigma(X) \neq \varepsilon \}$ est un ensemble ordonné qui permet de graduer l'*Echelle Jassologique* de Ω.
Pour tout entier $p \geq 0$, on a défini l'ensemble DNJ(p) des *Stratinos* non unitaires de longueur p et des Stratinos unitaires de longueur p+1. Définissons alors $DNJ(\Omega, p) = \{ X \in DNJ(p) \mid \Sigma(X) \neq \varepsilon \}$ et DNJU(Ω, p) le sous-ensemble des *Stratinos* Unitaires contenus dans DNJ(Ω, p), qui sont donc de longueur p+1.
L'ensemble des DNJ(Ω, p), $0 \leq p \leq m$, constitue alors une partition de NJ(Ω).

**4.5.2. Propriété $R_0$ :**

$DNJ(0) = \{\varepsilon, (1)\}$.
$\Sigma(\varepsilon) = (w_0)$, binôme($w_0$) = $\mathfrak{C}\mathfrak{D}$, médiane($w_0$) = non, $H(w_0) = (w_1)$ et $\varphi(w_0) = \varepsilon$.
$\Sigma(1) = (w_1)$, $w_1$ = zouc($j_1$), binôme($w_1$) = ⸢⸣, médiane($w_1$) = non, et $\varphi(w_1) = w_0$.
On en déduit que DNJ(Ω, 0) = DNJ(0), avec DNJU(Ω, 0) = {(1)}.

**4.5.3. Propriété $R_p$ pour tout p, $1 \leq p \leq m$ :**

Pour tout *Stratino* unitaire $Y = (X, 1) \in DNJU(p-1)$ et tout entier $n \geq 1$ :
i) on définit H(e), ν(e), T(e) et S(e) pour toute cellule $e \in \Sigma(X, n)$,
ii) on définit Σ(X, n, 1), Σ(X, n#), Σ(X, n#, 1), Σ(X, n+1),
iii) on définit binôme(e), médiane(e), et φ(e) pour toute cellule e de Σ(X, n, 1), Σ(X, n#), Σ(X, n#, 1) ou Σ(X, n+1),
iv) pour toute cellule e de Σ(X, n, 1) ou Σ(X, n#, 1), il existe une et une seule *Rovéjasse* v de $j_{p+1}$ telle que e = zouc(v), et on définit baou(v).
v) pour toute cellule e de Σ(X, n+1), on définit gh(e) et dh(e) : le premier terme de gh(e) et le dernier terme de dh(e) appartiennent à $j_p$ ou $j_{p+1}$ mais pas à $j_{p-1}$.

En particulier, nous verrons qu'il existera nécessairement un entier $n \geq 1$ unique tel que $\Sigma(X, n) \neq \varepsilon$ et $\Sigma(X, n+1) = \varepsilon$, ce qui induira automatiquement que $\Sigma(X, n\#, 1) = \varepsilon$ et que pour tout $r \geq n+1$ :
$\Sigma(X, r, 1) = \Sigma(X, r\#) = \Sigma(X, r\#, 1) = \Sigma(X, r+1) = \varepsilon$.

On a vu que l'ensemble des D(Y) pour tous les *Stratinos* unitaires $Y \in DNJU(p-1)$, de longueur p, constitue une partition de DNJ(p), avec
$D(Y) = \{ (X, n) \mid n \text{ entier naturel} \geq 2\} \cup \{ (X, n\#) \mid n \text{ entier naturel} \geq 1\} \cup$
$\{ (X, n, 1) \mid n \text{ entier naturel} \geq 1\} \cup \{ (X, n\#, 1) \mid n \text{ entier naturel} \geq 1\}$.
D'après la propriété $R_p$, Σ(Z) sera donc définie pour tous les *Stratinos* $Z \in DNJ(p)$, ce qui permettra de déduire DNJ(Ω, p) et DNJU(Ω, p).

Démontrons alors la propriété $R_p$ par récurrence sur p pour tout p, $0 \leq p \leq m$.
Pour p = 0, c'est déjà fait.
Soit $p \geq 0$.
Supposons que la propriété $R_p$ soit vérifiée et vérifions la propriété $R_{p+1}$.
Soit $Y = (X, 1) \in$ DNJU(p), de longueur p+1.
Par hypothèse de récurrence, $\Sigma(X, 1)$ est bien définie et les propriétés iii) et iv) sont vérifiées pour toute cellule éventuelle e de $\Sigma(X, 1)$.
Montrons alors les propriétés i), ii), iii), iv) et v) par récurrence sur $n \geq 1$.

**A) Pour n = 1 :**

Si $s\Sigma(X, 1) = \varepsilon$, alors on posera $\Sigma(X, 1, 1) = \Sigma(X, 1\#) = \Sigma(X, 1\#, 1) = \Sigma(X, 2) = \varepsilon$.
Sinon notons $s\Sigma(X, 1) = (e_1, e_2, \ldots, e_r)$, avec $r \geq 1$.
Pour tout t, $1 \leq t \leq r$, $e_t$ vérifie les propriétés iii) et iv) : binôme($e_t$), médiane($e_t$), $\varphi(e_t)$, sont bien définis, il existe une et une seule *Rovéjasse* $\varpi_t$ de $j_{p+1}$ telle que $e_t = \text{zouc}(\varpi_t)$, et baou($\varpi_t$) est bien défini.
On en déduit caouly($\varpi_t$) et fan($\varpi_t$) puisque ceux-ci se définissent à partir de baou($\varpi_t$) et médiane($e_t$).
Et d'après la propriété P1, $\pi(\varpi_t)$ est contenue dans une et une seule Rovéjasse $\omega_t$ de $j_p$.
Par définition, baou($\varpi_t$), caouly($\varpi_t$) et fan($\varpi_t$) sont des *Nesiles* de $\pi(\varpi_t)$, donc de $\omega_t$.

**A.1. Calcul de gh($e_t$) pour tout t, $1 \leq t \leq r$ :**
<u>Cas où $\varpi_t = j_1$ :</u>
Alors $e_t = w_1$.
Notons $\overline{B}(w_1) = (g_1, g_2, \ldots, g_r)$, avec $r \geq 3$ et $g_1 = w_0$ par hypothèse de départ.
Alors $gh(w_1) = (g_1, g_2)$.
<u>Cas où $\varpi_t \neq j_1$ :</u>
Si médiane($e_t$) = non, alors gh($e_t$) = baou($\varpi_t$).
Si médiane($e_t$) = oui, notons baou($\varpi_t$) = $(g_a, g_b)$ et caouly($\varpi_t$) = $(g'_b, g'_a)$, sachant que $g'_a \neq g_a$.
Alors gh($e_t$) = $(g'_a, g_a)$.
On a vu que le calcul de gh($e_t$) nous permet de définir $\overline{\overline{B}}°(e_t)$.
Dans tous les cas, les premier et dernier termes de $\overline{\overline{B}}°(e_t)$ appartiennent à $\omega_t$.

**A.2. Calcul de dh($e_t$) pour tout t, $1 \leq t \leq r$ :**
Dans tous les cas, dh($e_t$) sera le suffixe de $\overline{\overline{B}}°(e_t)$ constitué des deux derniers termes de $\overline{\overline{B}}°(e_t)$.
Ce qui induit que $\overline{\overline{Br}}(e_t) = \overline{\overline{B}}°(e_t)$.

**A.3. Calcul de $\nu(e_t)$, $\xi g(e_t)$, $\xi d(e_t)$, $\xi\Delta_k(e_t)$ et $\xi R_k(e_t)$ pour tout t, $1 \leq t \leq r$ :**
Comme les premier et dernier termes de $\overline{\overline{Br}}(e_t)$ appartiennent à $\omega_t$, $\overline{\overline{Br}}(e_t)$ se décompose de façon unique $\overline{\overline{Br}}(e_t) = (d_0, r_1, d_1, r_2, d_2, \ldots, r_\rho, d_\rho)$, avec $\rho \geq 0$, où pour tout k, $0 \leq k \leq \rho$ :

- $d_k$ est un *Nesile* non vide de $\omega_t$,
- si $1 \leq k \leq \rho$, $r_k$ est un *Nesile* non vide de $\Omega \backslash \omega_t$.

Posons : $\nu(e_t) = \rho$, $\xi g(e_t) = \xi d(e_t) = \varepsilon$, et si $\rho \geq 1$, pour tout k, $1 \leq k \leq \rho$, $\xi\Delta_\kappa(e_t) = d_k$ et $\xi R_k(e_t) = r_k$.
Si médiane($e_t$) = non, alors $\xi\Delta_0(e_t) = d_0$. D'où :
$\overline{\overline{Br}}(e_t) = (\xi\Delta_0(e_t), \xi R_1(e_t), \xi\Delta_1(e_t), \xi R_2(e_t), \xi\Delta_2(e_t), ..., \xi R_{\nu(et)}(e_t), \xi\Delta_{\nu(et)}(e_t))$, avec
$\overline{\overline{Br}}(e_t) = (\xi\Delta_0(e_t))$ si $\nu(e_t) = 0$.
Remarquons que dans le cas où p = 0, avec $e_t = w_1$, comme $\omega_t = j_0 = \{w_0\}$, alors pour tout k, $0 \leq k \leq \rho$, $\xi\Delta_k(w_1) = (w_0)$.

Si médiane($e_t$) = oui, alors $\varpi_t \neq j_1$ et gh($e_t$) est nécessairement un préfixe de $d_0$. Notons à nouveau $g'_a$ le premier terme de gh($e_t$). Alors $\xi\Delta_0(e_t)$ est le suffixe de $d_0$ tel que $d_0 = (g'_a, \xi\Delta_0(e_t))$. Remarquons que $\xi\Delta_0(e_t)$ est non vide. D'où :
$\overline{\overline{Br}}(e_t) = (g'_a, \xi\Delta_0(e_t), \xi R_1(e_t), \xi\Delta_1(e_t), \xi R_2(e_t), \xi\Delta_2(e_t), ..., \xi R_{\nu(et)}(e_t), \xi\Delta_{\nu(et)}(e_t))$, avec
$\overline{\overline{Br}}(e_t) = (g'_a, \xi\Delta_0(e_t))$ si $\nu(e_t) = 0$.

**A.4. Calcul de $\xi Rg(e_t)$, $\xi Rd(e_t)$, $\xi G(e_t)$ et $\xi H(e_t)$ pour tout t, $1 \leq t \leq r$ :**
Comme $\xi g(e_t) = \xi d(e_t) = \varepsilon$, alors on posera $\xi Rg(e_t) = \xi Rd(e_t) = \xi G(e_t) = \varepsilon$. D'où :

$\xi H(e_t) = (\xi\Delta_0(e_t), \xi R_1(e_t), \xi\Delta_1(e_t), \xi R_2(e_t), \xi\Delta_2(e_t), ..., \xi R_{\nu(et)}(e_t), \xi\Delta_{\nu(et)}(e_t))$, avec
$\xi H(e_t) = (\xi\Delta_0(e_t))$ si $\nu(e_t) = 0$.

**A.5. Calcul de $G(e_t)$, $\Delta G(e_t)$, $Rg(e_t)$ et $Rd(e_t)$, pour tout t, $1 \leq t \leq r$ :**
Comme $\xi Rg(e_t) = \xi Rd(e_t) = \xi G(e_t) = \varepsilon$, alors on posera $Rg(e_t) = Rd(e_t) = G(e_t) = \Delta G(e_t) = \varepsilon$.

**A.6. Calcul des $\Delta_k(e_t)$ pour tout t, $1 \leq t \leq r$ :**
Soit k, $0 \leq k \leq \nu(e_t)$.
Dans le cas où p = 0, on a vu que $e_t = w_1$ et $\xi\Delta_k(w_1) = (w_0)$.
On posera alors $\Delta_k(w_1) = \varepsilon$.
Supposons que $p \geq 1$ et notons $\xi\Delta_k(e_t) = (g_0, g_1, g_2, \ldots, g_u)$, avec $u \geq 0$.
Il y a 4 cas à distinguer :
1er cas : $(g_0, g_1) = $ baou$(\varpi_t)$ et $(g_{u-1}, g_u) = $ caouly$(\varpi_t)$, auquel cas $u \geq 2$.
Alors $\Delta_k(e_t) = ($ʋ$, g_2, g_3, \ldots, g_{u-1},$ ᴖ$)$, avec $\Delta_k(e_t) = ($ʋ, ᴖ$)$ si u = 2.

2ème cas : $(g_0, g_1) = $ baou$(\varpi_t)$ et $(g_{u-1}, g_u) \neq $ caouly$(\varpi_t)$, auquel cas $u \geq 1$.
Alors $\Delta_k(e_t) = ($ʋ$, g_2, g_3, \ldots, g_u)$, avec $\Delta_k(e_t) = ($ʋ$)$ si u = 1.

3ème cas : $(g_0, g_1) \neq $ baou$(\varpi_t)$ et $(g_{u-1}, g_u) = $ caouly$(\varpi_t)$, auquel cas $u \geq 1$.
Alors $\Delta_k(e_t) = (g_1, g_2, \ldots, g_{u-1},$ ᴖ$)$, avec $\Delta_k(e_t) = ($ᴖ$)$ si u = 1.

4ème cas : $(g_0, g_1) \neq $ baou$(\varpi_t)$ et $(g_{u-1}, g_u) \neq $ caouly$(\varpi_t)$.
Alors $\Delta_k(e_t) = (g_1, g_2, \ldots, g_u)$, avec $\Delta_k(e_t) = \varepsilon$ si u = 0.

**A.7. Calcul des $R_k(e_t)$, $T_k(e_t)$, et $S_k(e_t)$ pour tout t, $1 \leq t \leq r$ :**
Supposons que $\nu(e_t) \geq 1$ et soit k, $1 \leq k \leq \nu(e_t)$.
Par construction, $\xi R_k(e_t)$ est un *Nesile* non vide de $\Omega\backslash\omega_t$.
Donc nécessairement, les cellules de $\xi R_k(e_t)$ appartiennent à $\varpi_t$ ou à une *Rovéjasse* de $j_{p+2}$ orthogonale à $\varpi_t$.
$\xi R_k(e_t)$ se décompose alors de façon unique $\xi R_k(e_t) = (s_0, t_1, s_1, t_2, s_2, \ldots, t_u, s_u)$, avec $u \geq 0$, où pour tout q, $0 \leq q \leq u$ :
- $s_q$ est un *Nesile* non vide de $\varpi_t$,
- si $q \geq 1$, $t_q$ est un *Nesile* non vide de $\Omega\backslash\varpi_t$, donc d'une *Rovéjasse* $w_q$ de $j_{p+2}$.

Si u = 0, alors $T_k(e_t) = \varepsilon$, et $R_k(e_t) = S_k(e_t) = (s_0,$ ɷ$)$.
Supposons que $u \geq 1$.
Alors d'après la propriété P3, $w_1, w_2, \ldots, w_u$ sont nécessairement des *Rovéjasses* de $j_{p+2}$ distinctes.
Pour tout q, $1 \leq q \leq u$, zouc$(w_q)$ sera par définition le premier terme de $t_q$. Alors $c(w_q)^+(e_t)$ est le dernier terme de $s_{q-1}$ et $c(w_q)^-(e_t)$ le premier terme de $s_q$.
On posera baou$(w_q) = (e_t, c(w_q)^+(e_t))$, $\varphi($zouc$(w_q)) = e_t$, binôme(zouc$(w_q)$) = ɾ9 et médiane(zouc$(w_q)$) = non.
On en déduit caouly$(w_q) = (c(w_q)^-(e_t), e_t)$ et fan$(w_q)$, dont le dernier terme est $c(w_q)^-(e_t)$.
Pour toute cellule e de fan$(w_q)$, on posera binôme(e) = ʎ∂ , médiane(e) = non, et $\varphi(e)$ sera par définition la cellule de $w_q$ telle que $(c(w_q)^-(e), e)$ soit un facteur de $\overline{B}(\varphi(e))$.
Puis notons $s'_q$ le suffixe de $s_q$ tel que $s_q = (c(w_q)^-(e_t), s'_q)$.
Pour toute cellule e de $s'_q$, on posera binôme(e) = ɾ) , médiane(e) = non, et $\varphi(e) = e_t$.
De même, pour toute cellule e de $s_0$, on posera binôme(e) = ɾ) , médiane(e) = non, et $\varphi(e) = e_t$.
Alors $R_k(e_t) = (s_0,$ zouc$(w_1)$, •, $s'_1$, zouc$(w_2)$, •, $s'_2, \ldots,$ zouc$(w_u)$, •, $s'_u$, ɷ$)$
$T_k(e_t) = ($zouc$(w_1)$, •, zouc$(w_2)$, •, …, zouc$(w_u)$, •$)$, et
$S_k(e_t) = (s_0,$ fan$(w_1)$, $s'_1$, fan$(w_2)$, $s'_2, \ldots,$ fan$(w_u)$, $s'_u$, ɷ$)$.

Calcul de gh(c) et dh(c) pour toute cellule $c \in S_k(e_t)$ :
Notons $S_k(e_t) = (c_1, c_2, \ldots, c_\gamma,$ ɷ$)$ avec $\gamma \geq 1$, a le dernier terme de $\xi\Delta_{k-1}(e_t)$ et b le premier terme de $\xi\Delta_k(e_t)$ .
Alors $gh(c_1) = (e_t, a) = (\varphi(c_1) , a)$, et $dh(c_\gamma) = (b, e_t)$.
Si $\gamma \geq 2$, alors pour tout s, $2 \leq s \leq \gamma$, $gh(c_s) = (\varphi(c_s), c_{s-1})$, et pour tout s, $1 \leq s \leq \gamma-1$, $dh(c_s) = (c_{s+1}, \varphi(c_{s+1}))$.
Dans tous les cas, le premier terme de $gh(c_s)$ et le dernier terme de $dh(c_s)$ appartiennent à $j_{p+1}$ ou $j_{p+2}$, mais pas à $j_p$.

**A.8. Calcul de $H(e_t)$, $T(e_t)$, et $S(e_t)$ pour tout t, $1 \leq t \leq r$ :**

On en déduit que $H(e_t) = (\Delta_0(e_t), R_1(e_t), \Delta_1(e_t), R_2(e_t), \Delta_2(e_t), ..., R_{\nu(e)}(e_t), \Delta_{\nu(e)}(e_t))$, avec $H(e_t) = (\Delta_0(e_t))$ si $\nu(e) = 0$.
De plus, $T(e_t) = (T_1(e_t), T_2(e_t), ..., T_{\nu(e)}(e_t))$, avec $T(e_t) = \varepsilon$ si $\nu(e) = 0$, et $S(e_t) = (S_1(e_t), S_2(e_t), ..., S_{\nu(e)}(e_t))$, avec $S(e_t) = \varepsilon$ si $\nu(e) = 0$.

**A.9. Conclusion :**

Pour tout t, $1 \leq t \leq r$, $H(e_t)$, $\nu(e_t)$, $T(e_t)$ et $S(e_t)$ sont bien définis donc i) est vérifiée.
Par définition :
$\Sigma(X, 1, 1) = \varepsilon$,
$\Sigma(X, 1\#) = \varepsilon$,
$\Sigma(X, 1\#, 1) = (T(e_1), T(e_2), \ldots, T(e_u))$,
$\Sigma(X, 2) = (S(e_1), S(e_2), \ldots, S(e_u))$.
Donc ii) est vérifiée.
Pour toute cellule e de $\Sigma(X, 1\#, 1)$ ou $\Sigma(X, 2)$, on a défini binôme(e), médiane(e), et $\varphi(e)$, donc iii) est vérifiée.
Pour toute cellule e de $\Sigma(X, 1\#, 1)$, on a vu qu'il existait une et une seule *Rovéjasse* v de $j_{p+2}$ telle que e = zouc(v), et on a défini baou(v). Donc la propriété iv) est vérifiée.
Pour toute cellule e de $\Sigma(X, 2)$, on a défini gh(e) et dh(e) : le premier terme de gh(e) et le dernier terme de dh(e) appartiennent à $j_{p+1}$ ou $j_{p+2}$ mais pas à $j_p$. Donc la propriété v) est vérifiée.

**B) Pour tout $n \geq 2$ :**

Soit $n \geq 2$ et supposons que les propriétés i), ii), iii), iv) et v) soient vérifiées pour n-1.
Donc $\Sigma(X, n)$ est bien définie.
Si $s\Sigma(X, n) = \varepsilon$, alors on posera $\Sigma(X, n, 1) = \Sigma(X, n\#) = \Sigma(X, n\#, 1) = \Sigma(X, n+1) = \varepsilon$.
Sinon notons $s\Sigma(X, n) = (e_1, e_2, \ldots, e_r)$, avec $r \geq 1$.
Pour tout t, $1 \leq t \leq r$, $e_t$ vérifie les propriétés iii) et v) : binôme($e_t$), médiane($e_t$), $\varphi(e_t)$, gh($e_t$) et dh($e_t$) sont bien définis : le premier terme de gh($e_t$) et le dernier terme de dh($e_t$) appartiennent à $j_{p+1}$ ou $j_{p+2}$, mais pas à $j_p$.
On en déduit $\overline{B}°(e_t)$ et $\overline{Br}(e_t)$, dont le préfixe est gh($e_t$) et le suffixe dh($e_t$) : les premier et dernier termes de $\overline{Br}(e_t)$ appartiennent donc à $j_{p+1}$ ou $j_{p+2}$, mais pas à $j_p$.
D'autre part, on notera $\varpi_t$ la *Rovéjasse* de $j_{p+1}$ à laquelle $e_t$ appartient.
Par définition, baou($\varpi_t$), caouly($\varpi_t$) et fan($\varpi_t$) sont des nesiles de $\pi(\varpi_t)$, donc de $j_p$.

**P4 : 4ème propriété caractéristique des cartes planaires cubiques :**

Les cellules de $s\Sigma(X, n) = (e_1, e_2, \ldots, e_r)$, avec $r \geq 1$, vérifient les propriétés suivantes :
1) $e_1, e_2, \ldots, e_r$ sont des cellules distinctes de $j_{p+1}$.
2) pour tout t, $1 \leq t \leq r$, $e_t$ n'est pas un terme de $\Sigma(X, n-1)$.
3) si $r \geq 2$, pour tout s, $1 \leq s \leq r-2$, et tout t, $s+2 \leq t \leq r$, $e_s$ n'est pas en contact avec $e_t$.
4) si $r \geq 2$,pour tout t, $1 \leq t \leq r-1$, $e_t$ est en contact avec $e_{t+1}$ si et seulement si le terme « ɷ »ne se trouve pas entre $e_t$ et $e_{t+1}$ dans $\Sigma(X, n)$.

**B.1. Calcul de $\nu(e_t)$, $\xi g(e_t)$, $\xi d(e_t)$, $\xi\Delta_k(e_t)$ et $\xi R_k(e_t)$ pour tout t, $1 \leq t \leq r$ :**

Comme les premier et dernier terme de $\overline{Br}(e_t)$ n'appartiennent pas à $j_p$, $\overline{Br}(e_t)$ se décompose de façon unique
$\overline{Br}(e_t) = (r_0, d_0, r_1, d_1, r_2, d_2, \ldots, r_\rho, d_\rho, r_{\rho+1})$ avec $\rho \geq 0$, où :

- pour tout k, $0 \leq k \leq \rho$, $d_k$ est un *Nesile* non vide de $j_p$,
- pour tout k, $0 \leq k \leq \rho+1$, $r_k$ est un *Nesile* non vide de $\Omega \backslash j_p$.

Posons : $\nu(e_t) = \rho$, $\xi g(e_t) = r_0$, $\xi d(e_t) = r_{\rho+1}$ , pour tout k, $0 \leq k \leq \rho$, $\xi\Delta_k(e_t) = d_k$, et si $\rho \geq 1$, pour tout k, $1 \leq k \leq \rho$, $\xi R_k(e_t) = r_k$. D'où :
$\overline{Br}(e_t) = (\xi g(e_t), \xi\Delta_0(e_t), \xi R_1(e_t), \xi\Delta_1(e_t), \xi R_2(e_t), \xi\Delta_2(e_t), ..., \xi R_{\nu(et)}(e_t), \xi\Delta_{\nu(et)}(e_t), \xi d(e_t))$, avec
$\overline{Br}(e_t) = (\xi g(e_t), \xi\Delta_0(e_t), \xi d(e_t))$ si $\nu(e_t) = 0$.
Remarquons que dans le cas où p = 0, comme $j_0 = \{w_0\}$, alors pour tout k, $0 \leq k \leq \rho$, $\xi\Delta_\kappa(e_t) = (w_0)$.
Dans le cas où $p \geq 1$, pour tout k, $0 \leq k \leq \rho$, $\xi\Delta_k(e_t)$ est un *Nesile* de $\pi(\varpi_t)$.

**B.2. Calcul de $\xi Rg(e_t)$, $\xi Rd(e_t)$, $\xi G(e_t)$ et $\xi H(e_t)$ pour tout t, $1 \leq t \leq r$ :**

$\xi Rg(e_t)$ se déduit de $\xi g(e_t)$ de la façon suivante :
Si $\xi g(e_t)$ est constitué d'une seul terme, alors $\xi Rg(e_t) = \varepsilon$.
Sinon notons $\xi g(e_t) = (a_1, a_2, \ldots, a_v)$, avec $v \geq 2$, vérifiant $(a_1, a_2) = gh(e_t)$.
Par construction, on vérifie que nécessairement $t \geq 2$ et $a_2 = e_{t-1}$.
Il existe u le plus grand entier, $2 \leq u \leq v$, tel que $a_u = a_2$.
Alors $\xi Rg(e_t)$ est le suffixe de $\xi g(e_t)$ vérifiant $\xi g(e_t) = (a_1, a_2, \ldots, a_u, \xi Rg(e_t))$.

$\xi Rd(e_t)$ et $\xi G(e_t)$ se déduisent de $\xi d(e_t)$ de la façon suivante :
Si $\xi d(e_t)$ est constitué d'une seul terme, alors $\xi Rd(e_t) = \xi G(e_t) = \varepsilon$.
Sinon notons $\xi d(e_t) = (b_1, b_2, \ldots, b_s)$, avec $s \geq 2$, vérifiant $(b_{s-1}, b_s) = dh(e_t)$.
Par construction, on vérifie que nécessairement $t < r$ et $b_{s-1} = e_{t+1}$.
Il existe r le plus grand entier, $1 \leq r \leq s-1$, tel que $b_r = b_{s-1}$.
Alors $\xi G(e_t) = (b_r, b_{r+1}, \ldots, b_{s-1})$ et $\xi Rd(e_t)$ est le préfixe de $\xi d(e_t)$ vérifiant $\xi d(e_t) = (\xi Rd(e_t), \xi G(e_t), b_s)$.

On en déduit $\xi H(e_t)$, qui est un facteur de $\overline{\overline{Br}}(e_t)$ :
$\xi H(e_t) = (\xi Rg(e_t), \xi\Delta_0(e_t), \xi R_1(e_t), \xi\Delta_1(e_t), \xi R_2(e_t), \xi\Delta_2(e_t), \ldots, \xi R_{\nu(et)}(e_t), \xi\Delta_{\nu(et)}(e_t), \xi Rd(e_t), \xi G(e_t))$, avec
$\xi H(e_t) = (\xi Rg(e_t), \xi\Delta_0(e_t), \xi Rd(e_t), \xi G(e_t))$ si $\nu(e_t) = 0$.

**B.3. Calcul des $\Delta_k(e_t)$ pour tout t, $1 \leq t \leq r$ :**

Soit k, $0 \leq k \leq \nu(e_t)$.
Dans le cas où $p = 0$, on a vu que $\xi\Delta_k(e_t) = (w_0)$.
On posera alors $\Delta_k(e_t) = \varepsilon$.
Supposons que $p \geq 1$ et notons $\xi\Delta_k(e_t) = (g_0, g_1, g_2, \ldots, g_u)$, qui est un nesile de $\pi(\varpi_t)$, avec $u \geq 0$.
Il y a 4 cas à distinguer :
<u>1er cas :</u> $(g_0, g_1) = baou(\varpi_t)$ et $(g_{u-1}, g_u) = caouly(\varpi_t)$, auquel cas $u \geq 2$.
Alors $\Delta_k(e_t) = ($ʋ$, g_2, g_3, \ldots, g_{u-1},$ ɷ$)$, avec $\Delta_k(e_t) = ($ʋ, ɷ$)$ si $u = 2$.

<u>2ème cas :</u> $(g_0, g_1) = baou(\varpi_t)$ et $(g_{u-1}, g_u) \neq caouly(\varpi_t)$, auquel cas $u \geq 1$.
Alors $\Delta_k(e_t) = ($ʋ$, g_2, g_3, \ldots, g_u)$, avec $\Delta_k(e_t) = ($ʋ$)$ si $u = 1$.

<u>3ème cas :</u> $(g_0, g_1) \neq baou(\varpi_t)$ et $(g_{u-1}, g_u) = caouly(\varpi_t)$, auquel cas $u \geq 1$.
Alors $\Delta_k(e_t) = (g_1, g_2, \ldots, g_{u-1},$ ɷ$)$, avec $\Delta_k(e_t) = ($ɷ$)$ si $u = 1$.

<u>4ème cas :</u> $(g_0, g_1) \neq baou(\varpi_t)$ et $(g_{u-1}, g_u) \neq caouly(\varpi_t)$.
Alors $\Delta_k(e_t) = (g_1, g_2, \ldots, g_u)$, avec $\Delta_k(e_t) = \varepsilon$ si $u = 0$.

**B.4. Calcul de $G(e_t)$ et $\Delta G(e_t)$, pour tout t, $1 \leq t \leq r$ :**

Si $\xi G(e_t) = \varepsilon$, alors $G(e_t) = \Delta G(e_t) = \varepsilon$.
Sinon notons a le dernier terme de $\xi G(e_t)$, qui est également le premier par construction.
On a vu qu'alors nécessairement $t < r$ et $a = e_{t+1}$.
Si $\xi G(e_t) = (e_{t+1})$, alors $G(e_t) = \Delta G(e_t) = \varepsilon$.
Sinon $\xi G(e_t)$ se décompose de façon unique $\xi G(e_t) = (e_{t+1}, G_1, e_{t+1}, G_2, \ldots, e_{t+1}, G_u, e_{t+1})$, avec $u \geq 1$, où pour tout s, $1 \leq s \leq u$, $G_s$ est un *Nesile* non vide de $\Omega \backslash \{e_{t+1}\}$.
La propriété P4 caractéristique des cartes planaires cubiques (voir ci-dessous) induit par récurrence que nécessairement pour tout s, $1 \leq s \leq u$, $G_s$ est un *Nesile* $j_{p+2}$, et donc d'une *Rovéjasse* $w_s$ de $j_{p+2}$.
En particulier, $\pi(w_s) = \{e_t, e_{t+1}\}$ : on posera $baou(w_s) = (e_t, e_{t+1})$ et $zouc(w_s)$ sera le premier terme de $G_s$.
De plus, $binôme(zouc(w_s)) =$ ❛❜, $médiane(zouc(w_s)) = non$, et $\varphi(zouc(w_s)) = e_t$.
On en déduit que $caouly(w_s) = (e_{t+1}, e_t)$ et $fan(w_s) = \varepsilon$.
Alors $G(e_t) = (zouc(w_1), zouc(w_2), \ldots, zouc(w_u))$ et $\Delta G(e_t) = ($ʋ, ɷ$)^u$.

**B.5. Calcul de $Rd(e_t)$, $Td(e_t)$, et $Sd(e_t)$ pour tout t, $1 \leq t \leq r$ :**

Si $\xi Rd(e_t) = \varepsilon$, alors $Rd(e_t) = Td(e_t) = Sd(e_t) = \varepsilon$.
Sinon, puisque $\xi Rd(e_t)$ est un facteur de $\xi d(e_t)$, les cellules de $\xi Rd(e_t)$ appartiennent à $\varpi_t$ ou à $j_{p+2}$.
Rappelons qu'alors $t < r$ et que le terme qui suit $\xi Rd(e_t)$ dans $\xi d(e_t)$ est $e_{t+1}$, le premier terme de $\xi G(e_t)$.
En particulier, le premier terme appartient nécessairement à $\varpi_t$.
$\xi Rd(e_t)$ se décompose alors de façon unique $\xi Rd(e_t) = (s_0, t_1, s_1, t_2, s_2, \ldots, t_u, s_u)$, avec $u \geq 0$, où pour tout q, $0 \leq q \leq u$ :

- $s_q$ est un *Nesile* non vide de $\varpi_t$, sauf éventuellement pour q = u,
- si $q \geq 1$, $t_q$ est un *Nesile* non vide de $j_{p+2}$, donc d'une *Rovéjasse* $w_q$ de $j_{p+2}$.

Si u = 0, alors $Td(e_t) = \varepsilon$, et $Rd(e_t) = Sd(e_t) = \xi Rd(e_t) = (s_0)$.
Supposons que $u \geq 1$.
Alors d'après la propriété P3, $w_1, w_2, \ldots, w_u$ sont nécessairement des *Rovéjasses* de $j_{p+2}$ toutes différentes.
Notons v l'entier défini par v = u si $s_u \neq \varepsilon$ et $v = u - 1$ si $s_u = \varepsilon$, auquel cas $u \geq 1$.
Pour tout q, $1 \leq q \leq v$, $zouc(w_q)$ sera par définition le premier terme de $t_q$. Alors $c(w_q)^+(e_t)$ est le dernier terme de $s_{q-1}$ et $c(w_q)^-(e_t)$ le premier terme de $s_q$.
On posera $baou(w_q) = (e_t, c(w_q)^+(e_t))$, $\varphi(zouc(w_q)) = e_t$, $binôme(zouc(w_q)) =$ [symbol] et $médiane(zouc(w_q)) =$ non.
On en déduit $caouly(w_q) = (c(w_q)^-(e_t), e_t)$ et $fan(w_q)$, dont le dernier terme est $c(w_q)^-(e_t)$.
Pour toute cellule e de $fan(w_q)$, on posera binôme(e) = [symbol] , et $\varphi(e)$ sera par définition la cellule de $w_q$ telle que $(c(w_q)^-(e), e)$ soit un facteur de $\overline{B}(\varphi(e))$.
Puis notons $s'_q$ le suffixe de $s_q$ tel que $s_q = (c(w_q)^-(e_t), s'_q)$.
Pour toute cellule e de $s'_q$, on posera binôme(e) = [symbol] , et $\varphi(e) = e_t$.
De même, pour toute cellule e de $s_0$, on posera binôme(e) = [symbol] , et $\varphi(e) = e_t$.
Il y a alors 2 cas à considérer :
<u>1er cas :</u> v = u.
Alors $Rd(e_t) = (s_0, zouc(w_1), \bullet, s'_1, zouc(w_2), \bullet, s'_2, \ldots, zouc(w_u), \bullet, s'_u)$
$Td(e_t) = (zouc(w_1), \bullet, zouc(w_2), \bullet, \ldots, zouc(w_u), \bullet)$, et
$Sd(e_t) = (s_0, fan(w_1), s'_1, fan(w_2), s'_2, \ldots, fan(w_u), s'_u)$.
Pour toute cellule e de $Sd(e_t)$, si e est le dernier terme de $Sd(e_t)$, alors médiane(e) = oui, sinon médiane(e) = non.
<u>2ème cas :</u> $v = u - 1$, c'est-à-dire $s_u = \varepsilon$.
Alors $e_t$ et $e_{t+1}$ appartiennent à $\pi(w_u)$, avec $c(w_q)^+(e_{t+1}) = e_t$.
$zouc(w_u)$ sera par définition le dernier terme de $t_q$, qui est en contact avec $e_t$ et $e_{t+1}$.
On posera $baou(w_u) = (e_t, c(w_q)^+(e_t))$, $\varphi(zouc(w_u)) = e_t$, $binôme(zouc(w_u)) =$ [symbol] et $médiane(zouc(w_u)) =$ oui.
On en déduit $caouly(w_u) = (c(w_q)^-(e_{t+1}), e_{t+1})$ et $fan(w_u)$, dont le dernier terme est $c(w_u)^-(e_{t+1})$.
Pour toute cellule e de $fan(w_u)$, on posera binôme(e) = [symbol] , et $\varphi(e)$ sera par définition la cellule de $w_u$ telle que $(c(w_u)^-(e), e)$ soit un facteur de $\overline{B}(\varphi(e))$.
Alors $Rd(e_t) = (s_0, zouc(w_1), \bullet, s'_1, zouc(w_2), \bullet, s'_2, \ldots, zouc(w_{u-1}), \bullet, s'_{u-1}, zouc(w_u))$
$Td(e_t) = (zouc(w_1), \bullet, zouc(w_2), \bullet, \ldots, zouc(w_{u-1}), \bullet, zouc(w_u))$, et
$Sd(e_t) = (s_0, fan(w_1), s'_1, fan(w_2), s'_2, \ldots, fan(w_{u-1}), s'_{u-1}, fan(w_u))$.
Pour toute cellule e de $Sd(e_t)$, médiane(e) = non.

<u>Calcul de gh(c) et dh(c) pour toute cellule $c \in Sd(e_t)$ :</u>
Dans tous les cas, notons $Sd(e_t) = (c_1, c_2, \ldots, c_\gamma)$ avec $\gamma \geq 1$, et a le dernier terme de $\xi\Delta_{v(et)}(e_t)$.
Alors $gh(c_1) = (e_t, a) = (\varphi(c_1), a)$, et $dh(c_\gamma) = (b, e_{t+1})$, où b dépend de $Sg(e_{t+1})$ : voir ci-dessous pour son calcul.
Si $Sg(e_{t+1})$ ne contient aucune cellule, alors b est le premier terme de $\xi\Delta_0(e_{t+1})$.
Sinon b est le premier terme de $Sg(e_{t+1})$.
Si $\gamma \geq 2$, alors pour tout s, $2 \leq s \leq \gamma$, $gh(c_s) = (\varphi(c_s), c_{s-1})$, et pour tout s, $1 \leq s \leq \gamma-1$, $dh(c_s) = (c_{s+1}, \varphi(c_{s+1}))$.
Dans tous les cas, le premier terme de $gh(c_s)$ et le dernier terme de $dh(c_s)$ appartiennent à $j_{p+1}$ ou $j_{p+2}$, mais pas à $j_p$.

**B.6. Calcul de $Rg(e_t)$, $Tg(e_t)$, et $Sg(e_t)$ pour tout t, $1 \leq t \leq r$ :**

Si $\xi Rg(e_t) = \varepsilon$, alors $Rg(e_t) = Tg(e_t) = Sg(e_t) = \varepsilon$.
Sinon, puisque $\xi Rg(e_t)$ est un facteur de $\xi g(e_t)$, les cellules de $\xi Rg(e_t)$ appartiennent à $\varpi_t$ ou à $j_{p+2}$.
Rappelons qu'alors $t \geq 2$ et le terme qui précède $\xi Rg(e_t)$ dans $\xi g(e_t)$ est $e_{t-1}$.
De plus, le dernier terme appartient nécessairement à $\varpi_t$.
$\xi Rg(e_t)$ se décompose alors de façon unique $\xi Rg(e_t) = (s_0, t_1, s_1, t_2, s_2, \ldots, t_u, s_u)$, avec $u \geq 0$, où pour tout q, $0 \leq q \leq u$ :

- $s_q$ est un *Nesile* non vide de $\varpi_t$, sauf éventuellement pour q = 0,
- si $q \geq 1$, $t_q$ est un *Nesile* non vide de $j_{p+2}$, donc d'une *Rovéjasse* $w_q$ de $j_{p+2}$.

Alors d'après la propriété P3, $w_1, w_2, \ldots, w_u$ sont nécessairement des *Rovéjasses* de $j_{p+2}$ distinctes.
Notons v l'entier défini par v = 1 si $s_0 \neq \varepsilon$ et v = 2 si $s_0 = \varepsilon$, auquel cas si u = 1, alors u < v.
Dans le cas où $v \leq u$, alors pour tout q, $v \leq q \leq u$, zouc($w_q$) sera par définition le premier terme de $t_q$. Alors $c(w_q)^+(e_t)$ est le dernier terme de $s_{q-1}$ et $c(w_q)^-(e_t)$ le premier terme de $s_q$.
On posera baou($w_q$) = ($e_t$, $c(w_q)^+(e_t)$), $\varphi$(zouc($w_q$)) = $e_t$, binôme(zouc($w_q$)) = ✗❜ et médiane(zouc($w_q$)) = non.
On en déduit caouly($w_q$) = ($c(w_q)^-(e_t)$, $e_t$) et fan($w_q$), dont le dernier terme est $c(w_q)^-(e_t)$.
Pour toute cellule e de fan($w_q$), on posera binôme(e) = ✗❛❜, médiane(e) = non et $\varphi$(e) sera par définition la cellule de $w_q$ telle que ($c(w_q)^-(e)$, e) soit un facteur de $\overline{B}(\varphi(e))$.
D'autre part, pour tout q, $v-1 \leq q \leq u$, notons $a_q$ le premier terme de $s_q$ et $s'_q$ le suffixe de $s_q$ tel que $s_q = (a_q, s'_q)$.
Pour toute cellule e de $s'_q$, on posera binôme(e) = ✗⁾, médiane(e) = non et $\varphi(e) = e_t$.
Il y a alors 2 cas à considérer :
<u>1er cas :</u> v = 1.
Alors $Rg(e_t)$ = ($s'_0$, zouc($w_1$), •, $s'_1$, zouc($w_2$), •, $s'_2$, …, zouc($w_u$), •, $s'_u$, ɷ)
$Tg(e_t)$ = (zouc($w_1$), •, zouc($w_2$), •, …, zouc($w_u$), •), et
$Sg(e_t)$ = ($s'_0$, fan($w_1$), $s'_1$, fan($w_2$), $s'_2$, …, fan($w_u$), $s'_u$, ɷ).
En particulier, si u = 0, alors $Rg(e_t) = Sg(e_t)$ = ($s'_0$, ɷ) et $Tg(e_t) = \varepsilon$.
<u>2ème cas :</u> v = 2, c'est-à-dire $s_0 = \varepsilon$.
Alors $Rg(e_t)$ = (•, $s'_1$, zouc($w_2$), •, $s'_2$, …, zouc($w_u$), •, $s'_u$, ɷ)
$Tg(e_t)$ = (•, zouc($w_2$), •, zouc($w_3$), •, …, zouc($w_u$), •), et
$Sg(e_t)$ = ($s'_1$, fan($w_2$), $s'_2$, …, fan($w_u$), $s'_u$, ɷ).
En particulier, si u = 1, alors $Rg(e_t)$ = (•, $s'_1$, ɷ), $Tg(e_t)$ = (•), et $Sg(e_t)$ = ($s'_1$, ɷ).
<u>Calcul de gh(c) et dh(c) pour toute cellule $c \in Sd(e_t)$ :</u>
Dans le cas où $Sg(e_t)$ contient au moins une cellule, notons $Sg(e_t)$ = ($c_1, c_2, \ldots, c_\gamma$, ɷ) avec $\gamma \geq 1$, et b le premier terme de $\xi\Delta_0(e_t)$.
Alors $dh(c_\gamma) = (b, e_t)$ et $gh(c_1) = (e_{t-1}, a)$, où a est le dernier terme de $Sd(e_{t-1})$ : voir ci-dessus pour son calcul.
Si $\gamma \geq 2$, alors pour tout s, $2 \leq s \leq \gamma$, $gh(c_s) = (\varphi(c_s), c_{s-1})$, et pour tout s, $1 \leq s \leq \gamma-1$, $dh(c_s) = (c_{s+1}, \varphi(c_{s+1}))$.
Dans tous les cas, le premier terme de $gh(c_s)$ et le dernier terme de $dh(c_s)$ appartiennent à $j_{p+1}$ ou $j_{p+2}$, mais pas à $j_p$.

**B.7. Calcul des $R_k(e_t)$, $T_k(e_t)$, et $S_k(e_t)$ pour tout t, $1 \leq t \leq r$ :**
Supposons que $\nu(e_t) \geq 1$ et soit k, $1 \leq k \leq \nu(e_t)$.
Par construction, $\xi R_k(e_t)$ est un *Nesile* non vide de $\Omega \backslash j_p$.
Donc nécessairement, les cellules de $\xi R_k(e_t)$ appartiennent à $\varpi_t$ ou à une Rovéjasse de $j_{p+2}$ orthogonale à $\varpi_t$.
$\xi R_k(e_t)$ se décompose alors de façon unique $\xi R_k(e_t) = (s_0, t_1, s_1, t_2, s_2, \ldots, t_u, s_u)$, avec $u \geq 0$, où pour tout q, $0 \leq q \leq u$ :

- $s_q$ est un *Nesile* non vide de $\varpi_t$,
- si $q \geq 1$, $t_q$ est un *Nesile* non vide de $\Omega \backslash j_p$, donc d'une *Rovéjasse* $w_q$ de $j_{p+2}$.

Si u = 0, alors $T_k(e_t) = \varepsilon$, et $R_k(e_t) = S_k(e_t)$ = ($s_0$, ɷ).
Supposons que $u \geq 1$.
Alors d'après la propriété P3, $w_1, w_2, \ldots, w_u$ sont nécessairement des *Rovéjasses* de $j_{p+2}$ toutes différentes.
Pour tout q, $1 \leq q \leq u$, zouc($w_q$) sera par définition le premier terme de $t_q$. Alors $c(w_q)^+(e_t)$ est le dernier terme de $s_{q-1}$ et $c(w_q)^-(e_t)$ le premier terme de $s_q$.
On posera baou($w_q$) = ($e_t$, $c(w_q)^+(e_t)$), $\varphi$(zouc($w_q$)) = $e_t$, binôme(zouc($w_q$)) = ✗❜ et médiane(zouc($w_q$)) = non.
On en déduit caouly($w_q$) = ($c(w_q)^-(e_t)$, $e_t$) et fan($w_q$), dont le dernier terme est $c(w_q)^-(e_t)$.
Pour toute cellule e de fan($w_q$), on posera binôme(e) = ✗❛❜, médiane(e) = non, et $\varphi$(e) sera par définition la cellule de $w_q$ telle que ($c(w_q)^-(e)$, e) soit un facteur de $\overline{B}(\varphi(e))$.
Puis notons $s'_q$ le suffixe de $s_q$ tel que $s_q = (c(w_q)^-(e_t), s'_q)$.
Pour toute cellule e de $s'_q$, on posera binôme(e) = ✗⁾, médiane(e) = non, et $\varphi(e) = e_t$.
De même, pour toute cellule e de $s_0$, on posera binôme(e) = ✗⁾, médiane(e) = non, et $\varphi(e) = e_t$.
Alors $R_k(e_t)$ = ($s_0$, zouc($w_1$), •, $s'_1$, zouc($w_2$), •, $s'_2$, …, zouc($w_u$), •, $s'_u$, ɷ)
$T_k(e_t)$ = (zouc($w_1$), •, zouc($w_2$), •, …, zouc($w_u$), •), et
$S_k(e_t)$ = ($s_0$, fan($w_1$), $s'_1$, fan($w_2$), $s'_2$, …, fan($w_u$), $s'_u$, ɷ).

Calcul de gh(c) et dh(c) pour toute cellule $c \in S_k(e_t)$ :

Notons $S_k(e_t) = (c_1, c_2, \ldots, c_\gamma,$ ɷ) avec $\gamma \geq 1$, a le dernier terme de $\xi\Delta_{k-1}(e_t)$ et b le premier terme de $\xi\Delta_k(e_t)$ .
Alors $gh(c_1) = (e_t, a) = (\varphi(c_1), a)$, et $dh(c_\gamma) = (b, e_t)$.
Si $\gamma \geq 2$, alors pour tout s, $2 \leq s \leq \gamma$, $gh(c_s) = (\varphi(c_s), c_{s-1})$, et pour tout s, $1 \leq s \leq \gamma-1$, $dh(c_s) = (c_{s+1}, \varphi(c_{s+1}))$.
Dans tous les cas, le premier terme de $gh(c_s)$ et le dernier terme de $dh(c_s)$ appartiennent à $j_{p+1}$ ou $j_{p+2}$, mais pas à $j_p$.

**B.8. Calcul de $H(e_t)$, $T(e_t)$, et $S(e_t)$ pour tout t, $1 \leq t \leq r$ :**

On en déduit que $H(e_t) = (Rg(e_t), \Delta_0(e_t), R_1(e_t), \Delta_1(e_t), R_2(e_t), \Delta_2(e_t), ..., R_{\nu(e)}(e_t), \Delta_{\nu(e)}(e_t), Rd(e_t), G(e_t))$, avec
$H(e_t) = (Rg(e_t), \Delta_0(e_t), Rd(e_t), G(e_t))$ si $\nu(e) = 0$.
De plus, $T(e_t) = (Tg(e_t), T_1(e_t), T_2(e_t), ..., T_{\nu(e)}(e_t), Td(e_t))$, avec $T(e_t) = (Tg(e_t), Td(e_t))$ si $\nu(e) = 0$, et
$S(e_t) = (Sg(e_t), S_1(e_t), S_2(e_t), ..., S_{\nu(e)}(e_t), Sd(e_t))$, avec $S(e_t) = (Sg(e_t), Sd(e_t))$ si $\nu(e) = 0$.

**B.9. Conclusion :**

Pour tout t, $1 \leq t \leq r$, $H(e_t)$, $\nu(e_t)$, $T(e_t)$ et $S(e_t)$ sont bien définis donc i) est vérifiée.
Par définition :
$\Sigma(X, n, 1) = (G(e_1), G(e_2), \ldots, G(e_r))$,
$\Sigma(X, n\#) = (\Delta G(e_1), \Delta G(e_2), \ldots, \Delta G(e_r))$,
$\Sigma(X, n\#, 1) = (T(e_1), T(e_2), \ldots, T(e_r))$,
$\Sigma(X, n+1) = (S(e_1), S(e_2), \ldots, S(e_r))$.
Donc ii) est vérifiée.
Pour toute cellule e de $\Sigma(X, n, 1)$, $\Sigma(X, n\#)$, $\Sigma(X, n\#, 1)$ ou $\Sigma(X, n+1)$, on a défini binôme(e), médiane(e), et $\varphi(e)$, donc iii) est vérifiée.
Pour toute cellule e de $\Sigma(X, n, 1)$ ou $\Sigma(X, n\#, 1)$, on a vu qu'il existait une et une seule *Rovéjasse* v de $j_{p+2}$ telle que e = zouc(v), et on a défini baou(v). Donc la propriété iv) est vérifiée.
Pour toute cellule e de $\Sigma(X, n+1)$, on a défini gh(e) et dh(e) : le premier terme de gh(e) et le dernier terme de dh(e) appartiennent à $j_{p+1}$ ou $j_{p+2}$ mais pas à $j_p$. Donc la propriété v) est vérifiée.

**C) Remarque générale :**

Notons $s\Sigma(X, 1) = (zouc(\varpi_1), zouc(\varpi_2), \ldots, zouc(\varpi_r))$, avec $r \geq 1$, où $\varpi_1, \varpi_2, \ldots, \varpi_r$ sont des *Rovéjasses* de $j_{p+1}$ toutes différentes.
Notons $E = \varpi_1 \cup \varpi_2 \cup \ldots \cup \varpi_r$, qui est donc une union disjointe, et W l'ensemble des *Rovéjasses* de $j_{p+2}$ qui sont orthogonales à l'une des *Rovéjasses* $\varpi_1, \varpi_2, \ldots, \varpi_r$.
Pour tout $n \geq 1$, notons $E_n$ l'ensemble des cellules de $\Sigma(X, n)$ et $E_n^* = E \setminus (E_1 \cup E_2 \cup \ldots \cup E_n)$.
Puis notons $V_n$ l'ensemble des Rovéjasses correspondant aux cellules de $\Sigma(X, n, 1)$, $W_n$ l'ensemble de celles correspondant aux cellules de $\Sigma(X, n\#, 1)$, ainsi que
$V_n^* = W \setminus (V_1 \cup W_1 \cup V_2 \cup W_2 \cup \ldots \cup V_n)$, et $W_n^* = W \setminus (V_1 \cup W_1 \cup V_2 \cup W_2 \cup \ldots \cup V_n \cup W_n)$.
En particulier, $E_1 = \{zouc(\varpi_1), zouc(\varpi_2), \ldots, zouc(\varpi_r)\}$, $E_1^* = E \setminus E_1$, et $V_1 = \emptyset$, donc $V_1^* = W$.
Alors pour tout $n \geq 2$, :
i) $W_{n-1}$ constitue l'ensemble des *Rovéjasses* de $V_{n-1}^*$ qui sont en contact avec au moins une cellule de $E_{n-1}$,
ii) $E_n$ constitue l'ensemble des cellules de $E_n^*$ qui sont en contact avec au moins une cellule de $E_n$ ou une cellule d'une *Rovéjasse* appartenant à $W_{n-1}$,
iii) $V_n$ constitue l'ensemble des *Rovéjasses* de $W_{n-1}^*$ qui sont en contact avec exactement deux cellules consécutives de $s\Sigma(X, n)$.
Comme E est un ensemble fini, il existe nécessairement un et un seul entier $n_0 \geq 1$ tel que $E_{n0} \neq \emptyset$ et $E_r = \emptyset$ pour tout $r \geq n_0 + 1$. Ce qui induit que $V_r = W_r = \emptyset$ pour tout $r \geq n_0 + 1$, avec aussi $W_{n0} = \emptyset$.
Alors $E_1, E_2, \ldots, E_{n0}$, constituent une partition de E, avec $E_k \neq \emptyset$ pour tout k, $1 \leq k \leq n_0$.
Et $V_1, W_1, V_2, W_2, \ldots, W_{n0-1}, V_{n0}$ constituent une partition de W, à la différence que pour tout k, $1 \leq k \leq n_0$, $V_k$ ou $W_k$ sont éventuellement vides.

## 4.6. Calcul de MJ(Ω)

### 4.6.1. Le *Dallajascar* de Ω :

On a défini les fonctions de ramification $H : \Omega \rightarrow \Omega AJ^*$ et $\varphi : \Omega \rightarrow \Omega \cup \{\varepsilon\}$.
Soit $H' : \Omega \rightarrow \Omega^*$, défini de la façon suivante : pour toute cellule e de Ω, H'(e) est le sous-nesile de H(e) constitué des éléments de H(e) qui sont des cellules. En particulier, $H'(e) = \varepsilon$ si et seulement si $H(e) \in AJ^*$.
On vérifie alors que (φ,H') est un *Dallajascar* sur Ω, $w_0$ étant la racine.
En particulier, $\{H'(e) \mid x \in \Omega\}$ constitue une *Nesilpartition* de $\Omega \backslash \{w_0\}$.

### 4.6.2. La fonction LH :

La fonction $LH : \Omega \rightarrow \Omega AJ^*$ se définit à partir des fonctions $H : \Omega \rightarrow \Omega AJ^*$ et binôme : $\Omega \rightarrow AJ^*$ de la façon suivante : soit e une cellule de Ω et notons o(e) le caractère ouvrant de binôme(e) et f(e) le caractère fermant.
Si $H(e) = \varepsilon$, alors $LH(e) = o(e)\ f(e)$.
Sinon notons $H(e) = (\gamma_1, \gamma_2, \ldots, \gamma_u)$, avec $u \geq 1$.
Alors $LH(e) = o(e)\ \gamma_1\ \gamma_2 \ldots \gamma_u\ f(e)$.
En particulier, $LH(w_0) = \mathfrak{C}\ w_1\ \mathfrak{D}$.

### 4.6.3. Construction de MJ(Ω)

MJ(Ω) se construit alors de la façon suivante :
Pour t = 0, $MJ_0 = LH(w_0) = \mathfrak{C}\ w_1\ \mathfrak{D}$.
Pour t = 1 : $MJ_1 = \mathfrak{C}\ LH(w_1)\ \mathfrak{D}$.
Si $MJ_1 \in AJ^*$, alors c'est terminé : $MJ(\Omega) = M_1$.
Sinon notons $w_2$ le premier terme de $M_1$ qui soit une cellule de Ω.
$M_2$ s'obtient alors en remplaçant $w_2$ par $LH(w_2)$ dans $M_1$.
Si $MJ_2 \in AJ^*$, alors c'est terminé : $MJ(\Omega) = M_2$.
Sinon notons $w_3$ le premier terme de $M_2$ qui soit une cellule de Ω.
$M_3$ s'obtient alors en remplaçant $w_3$ par $LH(w_3)$ dans $M_2$.
Et ainsi de suite.
Il existe nécessairement un et un seul entier $N \geq 1$ tel que $MJ(\Omega) = M_N$, vérifiant $M_N \in AJ^*$ et $M_{N-1} \notin AJ^*$. D'où $\Omega = \{w_0, w_1, \ldots, w_N\}$.
On vérifie que pour tout p, $0 \leq p \leq N\text{-}1$, $w_p <_{(\varphi,H)} w_{p+1}$, où « $<_{(\varphi,H)}$ » est la relation d'ordre définie par le *Dallajascar* (φ,H).

### 4.6.4. Théorème fondamental :

Soit Ω et Ω' deux cartes planaires cubiques enracinées.
Alors Ω et Ω' sont équivalentes si et seulement si MJ(Ω) = MJ(Ω').

## 4.7. Exemple (Figure 1) :

Soit Ω = {a, b, c, d, e, f, g, h, i, j, k}, avec $w_0$ = a et $w_1$ = b, $i(w_1)$ = I et $i(w_0)$ = J (voir Figure 1).
Base de données : $\overline{B}(a)$ = (b, f, j, g, c, b), $\overline{B}(b)$ = (a, c, d, e, f, a), $\overline{B}(c)$ = (a, g, k, g, d, b, a), $\overline{B}(d)$ = (c, g, e, b, c), $\overline{B}(e)$ = (g, f, b, d, g), $\overline{B}(f)$ = (a, b, e, g, i, j, a), $\overline{B}(g)$ = (c, a, j, h, i, f, e, d, c, k, c), $\overline{B}(h)$ = (g, j, i, g), $\overline{B}(i)$ = (f, g, h, j, f), $\overline{B}(j)$ = (i, h, g, a, f, i), et $\overline{B}(k)$ = (c, g, c).

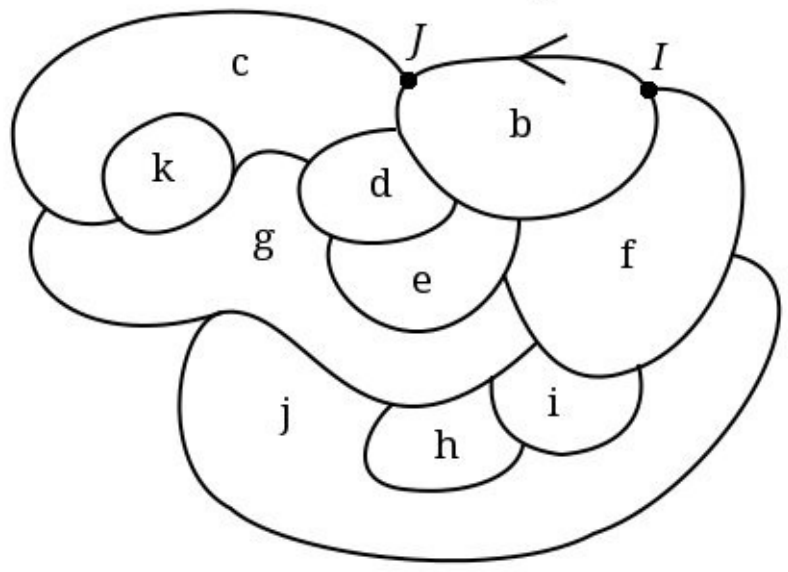

Figure 1. Carte initiale Ω.

**4.7.1. Les *Jasses* de Ω (Figure 2) :**
$j_0$ = {a} et comme $\overline{B}$(a) = (b, f, j, g, c, b), alors $j_1$ = {b, c, f, g, j}.
D'où $J_1 = \Omega \setminus (j_0 \cup j_1)$ = {d, e, h, i, k}.
Toutes les faces de $J_1$ sont accolées à au-moins une face de $j_1$ puisque d et e appartiennent à $\overline{B}$(b), h et i appartiennent à $\overline{B}$(g), et k appartient à $\overline{B}$(c). Donc $j_2 = J_1$ et $J_2 = \emptyset$, et donc $j_p = \emptyset$ pour tout $p \geq 3$.
D'où m = 2 : les *Jasses* $j_0$, $j_1$, $j_2$, constituent une partition de Ω. (voir Figure 2)

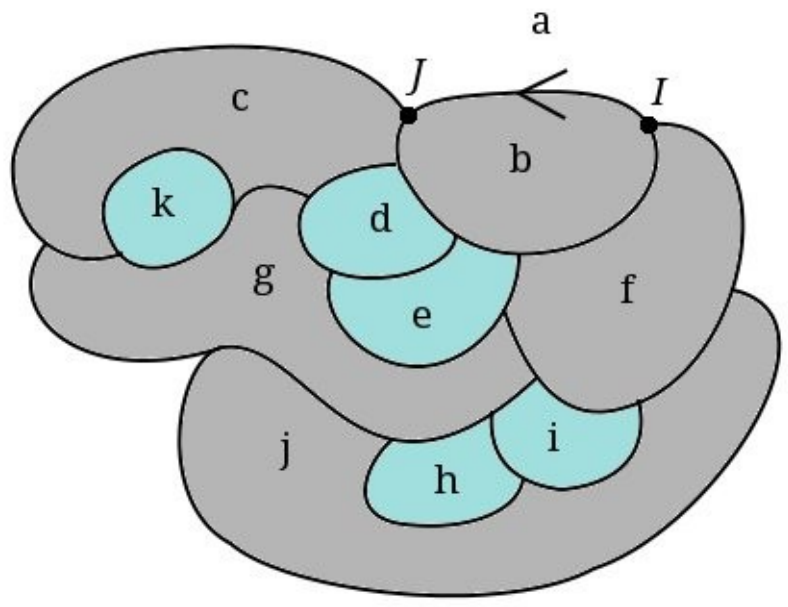

Figure 2.
Coloration en gris des faces de $j_1$, et en bleu celles de $j_2$.

**4.7.2. Les *Rovéjasses* de Ω :**
On constate bien que $j_1$ est connexe, ce qui induit que V(Ω,1) = {$j_1$}.
Les composantes connexes de $j_2$ = {d, e, h, i, k} sont u = {k}, v = {d, e} et w = {h, i}.
Donc V(Ω,2) = {u, v, w}.
Par définition, π(u) (resp. π(v) et π(w)) est l'ensemble des cellules de $j_1$ qui sont en contact avec au moins une cellule de u (resp. v et w). Donc π(u) = {c, g}, π(v) = {b, c, g, f}, et π(w) = {f, g, j}.
π(u), π(v) et π(w) sont connexes, donc la propriété P2 est bien vérifiée.

Vérifions la propriété P3 pour u :
$\overline{B}$(c) = (a, g, k, g, d, b, a) = (P, X, S), avec X = (k) ∈ u*, P = (a, g) ∈ (Ω\u)*, et S = (g, d, b, a) ∈ (Ω\u)*.
D'où $c(u)^+(c)$ = g et $c(u)^-(c)$ = g.
$\overline{B}$(g) = (c, a, j, h, i, f, e, d, c, k, c) = (P, X, S), avec X = (k) ∈ u*, P = (c, a, j, h, i, f, e, d, c) ∈ (Ω\u)*, et S = (c) ∈ (Ω\u)*.
D'où $c(u)^+(g)$ = c et $c(u)^-(g)$ = c.
La propriété P3 est donc bien vérifiée pour u.
Nécessairement zouc(u) = k. Supposons que baou(u) = (c, g), qui est bien défini puisque $c(u)^+(c)$ = g et $c(u)^-(g)$ = c.
Or $c(u)^+(g)$ = c, donc caouly(u) = (g, c) et fan(u) = ε.

Vérifions la propriété P3 pour v :
$\overline{B}$(b) = (a, c, d, e, f, a) = (P, X, S), avec X = (d, e) ∈ v*, P = (a, c) ∈ (Ω\v)*, et S = (f, a) ∈ (Ω\v)*.
D'où $c(v)^+(b)$ = c et $c(v)^-(b)$ = f.
$\overline{B}$(c) = (a, g, k, g, d, b, a) = (P, X, S), avec X = (d) ∈ v*, P = (a, g, k, g) ∈ (Ω\v)*, et S = (b, a) ∈ (Ω\v)*.
D'où $c(v)^+(c)$ = g et $c(v)^-(c)$ = b.
$\overline{B}$(g) = (c, a, j, h, i, f, e, d, c, k, c) = (P, X, S), avec X = (e, d) ∈ v*, P = (c, a, j, h, i, f) ∈ (Ω\v)*, et S = (c, k, c) ∈ (Ω\v)*.
D'où $c(v)^+(g)$ = f et $c(v)^-(c)$ = c.
$\overline{B}$(f) = (a, b, e, g, i, j, a) = (P, X, S), avec X = (e) ∈ v*, P = (a,b) ∈ (Ω\v)*, et S = (g, i, j, a) ∈ (Ω\v)*.
D'où $c(v)^+(f)$ = b et $c(v)^-(f)$ = g.
La propriété P3 est donc bien vérifiée pour v.
D'autre part, supposons que zouc(v) = d et baou(v) = (b, c), qui est bien défini puisque $c(v)^+(b)$ = c et $c(v)^-(c)$ = b.
$c(v)^+(c)$ = g, $c(v)^+(g)$ = f, et $c(v)^+(f)$ = b : (b, c, g, f) constitue un anneau autour de v.
De plus, supposons que médiane(d) = non.
Dans ce cas, caouly(v) = (f, b) et fan(v) = (g, f).

Vérifions la propriété P3 pour w :
$\overline{B}$(f) = (a, b, e, g, i, j, a) = (P, X, S), avec X = (i) ∈ w*, P = (a,b, e, g) ∈ (Ω\w)*, et S = (j, a) ∈ (Ω\w)*.
D'où $c(w)^+(f)$ = g et $c(w)^-(f)$ = j.

$\overline{B}(g) = (c, a, j, h, i, f, e, d, c, k, c) = (P, X, S)$, avec $X = (h, i) \in w^*$, $P = (c, a, j) \in (\Omega\backslash w)^*$, et $S = (f, e, d, c, k, c) \in (\Omega\backslash w)^*$. D'où $c(w)^+(g) = j$ et $c(w)^-(g) = f$.
$\overline{B}(j) = (i, h, g, a, f, i) = (X_1, M, X_2)$, avec $X_1 = (i, h) \in w^*$, $X_2 = (i) \in w^*$, et $M = (g, a, f) \in (\Omega\backslash w)^*$.
D'où $c(w)^+(j) = f$ et $c(w)^-(j) = g$.
La propriété P3 est donc bien vérifiée pour w.
D'autre part, supposons que zouc(w) = i et baou(w) = (g, j), qui est bien défini puisque $c(w)^+(g) = j$ et $c(w)^-(j) = g$.
$c(w)^+(j) = f$ et $c(w)^+(f) = g$ : (g, j, f) constitue un anneau autour de w.
De plus, supposons que médiane(i) = oui.
Dans ce cas, caouly(w) = (j, f) et fan(w) = ε.

### 4.7.3. Propriété $R_0$ :

$w_0 = a$ et $w_1 = b$, donc :
Σ(ε) = (a), binôme(a) = [symbole], médiane(a) = non, H(a) = (b) et φ(a) = ε.
Σ(1) = (b), b = zouc($j_1$), binôme(b) = [symbole], médiane(b) = non, et φ(b) = a.
DNJ(Ω, 0) = {ε, (1)}, avec DNJU(Ω, 0) = {(1)}.

### 4.7.4. Propriété $R_1$ :

Y = (1) est le seul *Stratino* unitaire de DNJU(Ω, 0).
Pour n = 1 : b est la seule cellule de Σ(1).
<u>A1 à A5</u> : $\overline{B}(b) = (a, c, d, e, f, a)$
Comme b = zouc($j_1$), alors gh(b) = (a, c), et donc $\overline{Br}(b) = \overline{B}°(b) = \overline{B}(b)$.
$\overline{Br}(b) = (d_0, r_1, d_1)$, avec $d_0 = d_1 = (a) \in j_0^*$, et $r_1 = (c, d, e, f) \in (\Omega\backslash j_0)^*$.
D'où ν(b) = 1, ξg(b) = ξd(b) = ε, $\xi\Delta_0(b) = \xi\Delta_1(b) = (a)$, et $\xi R_1(b) = (c, d, e, f)$.
Et donc $\xi H(b) = (\xi\Delta_0(b), \xi R_1(b), \xi\Delta_1(b))$.
<u>A6 à A7</u> :
$\Delta_0(b) = \Delta_1(b) = \varepsilon$.
Calculons $R_1(b)$, $T_1(b)$, et $S_1(b)$ :
$\xi R_1(b) = (c, d, e, f) = (s_0, t_1, s_1)$, avec $s_0 = (c) \in j_1^*$, $s_1 = (f) \in j_1^*$, et $t_1 = (d, e) \in (\Omega\backslash j_1)^*$.
En particulier, $t_1$ est un *Nesile* de la *Rovéjasse* v de $j_2$ orthogonale à $j_1$.
Alors zouc(v) = (d), baou(v) = (b, c), φ(d) = (b), binôme(d) = [symbole] et médiane(d) = non.
On a vu alors que caouly(v) = (f, b) et fan(v) = (g, f).
D'où binôme(g) = binôme(f) = [symbole] , médiane(g) = médiane(f) = non, et φ(g) et φ(f) appartiennent à v = {d, e}.
Comme (c, g) est un facteur de $\overline{B}(d) = (c, g, e, b, c)$, on en déduit que φ(g) = d.
Et comme (g, f)) est un facteur de $\overline{B}(e) = (g, f, b, d, g)$, on en déduit que φ(f) = e.
D'autre part, $s'_1 = \varepsilon$ puisque $s_1 = (f)$, et $s_0 = (c)$ induit que binôme(c) = [symbole] , médiane(c) = non, et φ(c) = b.
On en déduit que $R_1(b) = (s_0$, zouc(v), •, $s'_1$, ɷ) = (c, d, •, ɷ), $T_1(b)$ = (zouc(v), •) = (d, •), et
$S_1(b) = (s_0$, fan(v), $s'_1$, ɷ) = (c, g, f, ɷ).
Comme a est le dernier terme de $\xi\Delta_0(b)$ et le dernier de $\xi\Delta_1(b)$, alors gh(c) = (φ(c), a) = (b, a), et dh(f) = (a, b).
De plus, dh(c) = (g, φ(g) = (g, d), gh(g) = (φ(g), c) = (d, c), dh(g) = (f, φ(f) = (f, e), et gh(f) = (φ(f), g) = (e, g).
<u>A8 à A9</u> :
$H(b) = (\Delta_0(b), R_1(b), \Delta_1(b))$ = (c, d, •, ɷ),
$T(b) = T_1(b)$ = (d, •), et
$S(b) = S_1(b)$ = (c, g, f, ɷ).
Σ(1, 1) = ε,
Σ(1#) = ε,
Σ(1#, 1) = (T(b)) = (d, •), et
Σ(2) = (S(b)) = (c, g, f, ɷ).
binôme(d) = [symbole], médiane(d) = non, et φ(d) = b, avec d = zouc(v) et baou(v) = (b, c).
binôme(c) = [symbole] , médiane(c) = non, et φ(c) = b, avec gh(c) = (b, a) et dh(c) = (g, d).
binôme(g) = [symbole] , médiane(g) = non, et φ(g) = d, avec gh(g) = (d, c) et dh(g) = (f, e).
binôme(f) = [symbole] , médiane(f) = non, et φ(f) = e, avec gh(f) = (e, g) et dh(f) = (a, b).

Pour n = 2 : $s\Sigma(2)$ = (c, g, f).
En particulier, $s\Sigma(2)$ vérifie la propriété P4 puisque c, g et f sont des cellules distinctes, et n'appartiennent pas à $\Sigma(1)$, f n'est pas en contact avec c, g est en contact avec c et f est en contact avec g : $s\Sigma(2)$ constitue bien une chaîne de cellules.
$\overline{B}(c)$ = (a, g, k, g, d, b, a), avec gh(c) = (b, a) et dh(c) = (g, d).
Donc $\overline{\overline{B}}°(c)$ = (b, a, g, k, g, d, b) et $\overline{\overline{Br}}(c)$ = (b, a, g, k, g, d).
$\overline{B}(g)$ = (c, a, j, h, i, f, e, d, c, k, c), avec gh(g) = (d, c) et dh(g) = (f, e).
Donc $\overline{\overline{B}}°(g)$ = (d, c, k, c, a, j, h, i, f, e, d) et $\overline{\overline{Br}}(g)$ = (d, c, k, c, a, j, h, i, f, e).
$\overline{B}(f)$ = (a, b, e, g, i, j, a), avec gh(f) = (e, g) et dh(f) = (a, b).
Donc $\overline{\overline{B}}°(f)$ = (e, g, i, j, a, b, e) et $\overline{\overline{Br}}(f)$ = (e, g, i, j, a, b).

<u>B1 à B2</u> :
$\overline{\overline{Br}}(c) = (r_0, d_0, r_1)$, avec $d_0$ = (a) $\in j_0^*$, $r_0$ = (b) $\in (\Omega\backslash j_0)^*$, et $r_1$ = (g, k, g, d) $\in (\Omega\backslash j_0)^*$.
D'où ν(c) = 0, ξg(c) = (b), ξd(c) =(g, k, g, d), et $\xi\Delta_0(c)$ = (a).
ξg(c) est constitué d'un seul terme, donc ξRg(c) = ε.
ξd(c) = (ξRd(c), ξG(c), d), avec ξG(c) = (g, k, g) et ξRd(c) = ε.
D'où ξH(c) = ($\xi\Delta_0(c)$, ξG(c)).

$\overline{\overline{Br}}(g) = (r_0, d_0, r_1)$, avec $d_0$ = (a) $\in j_0^*$, $r_0$ = (d, c, k, c) $\in (\Omega\backslash j_0)^*$, et $r_1$ = (j, h, i, f, e) $\in (\Omega\backslash j_0)^*$.
D'où ν(g) = 0, ξg(g) = (d, c, k, c), ξd(g) =((j, h, i, f, e), et $\xi\Delta_0(g)$ = (a).
ξg(g) = (d, c, k, c, ξRg(g)), donc ξRg(g) = ε.
ξd(g) = (ξRd(g), ξG(g), e), avec ξG(g) = (f) et ξRd(g) = (j, h, i).
D'où ξH(g) = ($\xi\Delta_0(g)$, ξRd(g), ξG(g)).

$\overline{\overline{Br}}(f) = (r_0, d_0, r_1)$, avec $d_0$ = (a) $\in j_0^*$, $r_0$ = (e, g, i, j) $\in (\Omega\backslash j_0)^*$, et $r_1$ = (b) $\in (\Omega\backslash j_0)^*$.
D'où ν(f) = 0, ξg(f) = (e, g, i, j), ξd(f) =(b), et $\xi\Delta_0(f)$ = (a).
ξg(f) = (e, g, ξRg(f)), donc ξRg(f) = (i, j).
ξd(f) est constitué d'une seul terme, alors ξRd(f) = ξG(f) = ε.
D'où ξH(f) = (ξRg(f), $\xi\Delta_0(f)$).

<u>B3 à B4</u> :
$\Delta_0(c) = \Delta_0(g) = \Delta_0(f)$ = ε.
ξG(c) = (g, G, g), où G = (k) est un nesile non vide de Ω\{g}.
En particulier, G est un *nesile* de la *Rovéjasse* u de $j_2$ orthogonale à $j_1$.
Alors zouc(u) = (k), baou(v) = (c, g), φ(k) = (c), binôme(k) = ❛❜ et médiane(k) = non.
Ce qui induit que caouly(u) = (g, c) et fan(u) = ε.
Alors G(c) = (k) et ΔG(c) = ( ʋ , ɷ ).
ξG(g) estt constitué d'un seul terme, donc G(g) = ΔG(g) = ε.
ξG(f) = ε, donc G(f) = ΔG(f) = ε.

<u>B5 à B6</u> :
ξRd(g) = (j, h, i) = $(s_0, t_1, s_1)$, avec $s_0$ = (j) $\in j_1^*$, $t_1$ = (h, i) $\in (\Omega\backslash j_1)^*$, et $s_1$ = ε.
En particulier, $t_1$ est un *nesile* de la *Rovéjasse* w de $j_2$ orthogonale à $j_1$.
Comme $s_1$ = ε, on est dans le 2$^{ème}$ cas : g et f appartiennent à π(w), avec $c(w)^+(f)$ = g.
Alors zouc(w) = i, qui est le dernier terme de $t_1$ et en contact aved g et f.
Puis baou(w) = (g, j), φ(i) = (g), binôme(i) = ❛❜ et médiane(i) = oui.
On a vu alors que caouly(w) = (j, f) et fan(w) = ε.
$s_0$ = (j) induit que binôme(j) = ❛❜ , médiane(j) = non, et φ(j) = g.
On en déduit que Rd(g) = ($s_0$, zouc(w)) = (j, i), Td(g) = (zouc(w)) = (i), et Sd(g) = ($s_0$, fan(w)) = (j).
Comme a est le dernier terme de $\xi\Delta_0(g)$, alors gh(j) = (φ(j), a) = (g, a).
Nous allons voir que Sg(f) ne contient aucune cellule, et comme a est le dernier terme de $\xi\Delta_0(f)$, on en déduit que dh(j) = (a, f).

ξRg(f) = (i, j) = $(s_0, t_1, s_1)$, avec $s_1$ = (j) $\in j_1^*$, $t_1$ = (i) $\in (\Omega\backslash j_1)^*$, et $s_0$ = ε.
En particulier, $t_1$ est un *nesile* de la *Rovéjasse* w de $j_2$ orthogonale à $j_1$.

Comme $s_0 = \varepsilon$, on est dans le 2ème cas : g et f appartiennent à π(w), avec $c(w)^+(f) = g$.
$s_1 = (j, s'_1)$, donc $s'_1 = \varepsilon$.
On en déduit que Rg(f) = (•, ɷ), Tg(f) = (•), et Sg(f) = (ɷ).

<u>B8 à B9</u> :
H(c) = ($\Delta_0$(c), G(c)) = (k), avec T(c) = S(c) = ε.
H(g) = ($\Delta_0$(g), Rd(g), G(g)) = (j, i), avec T(g) = Td(g) = (i) et S(g) = Sd(g) = (j).
H(f) = (Rg(f), $\Delta_0$(f)) = (•, ɷ), avec T(f) = Tg(f) = (•), et S(f) = Sg(f) = (ɷ).
Σ(2, 1) = (G(c), G(g), G(f)) = (k),
Σ(2#) = (ΔG(c), ΔG(g), ΔG(f)) = (ʊ, ᴒ),
Σ(2#, 1) = (T(c), T(g), T(f)) = (i, •)
Σ(3) = (S(c), S(g), S(f)) = (j, ɷ).
binôme(k) = [symbole], médiane(k) = non, et φ(k) = c, avec k = zouc(u) et baou(u) = (c, g).
binôme(i) = [symbole], médiane(i) = oui, et φ(i) = g, avec i = zouc(w) et baou(w) = (g, j).
binôme(j) = [symbole], médiane(j) = non, et φ(j) = g, avec gh(j) = (g, a) et dh(j) = (a, f).

Pour n = 3 : sΣ(3) = (j).
Il est évident que sΣ(3) vérifie la propriété P4.
$\overline{B}$(j) = (i, h, g, a, f, i), avec gh(j) = (g, a) et dh(j) = (a, f).
Donc $\overline{B}$°(j) = (g, a, f, i, h, g) et $\overline{Br}$(j) = (g, a, f).

Alors $\overline{\overline{Br}}$(j) = ($r_0$, $d_0$, $r_1$), avec $d_0 = (a) \in j_0^*$, $r_0 = (g) \in (\Omega \backslash j_0)^*$, et $r_1 = (f) \in (\Omega \backslash j_0)^*$.
D'où ν(j) = 0, ξg(j) = (g), ξd(j) =(f), et ξ$\Delta_0$(j) = (a).
ξg(j) et ξd(j) sont constitués d'un seul terme, donc ξRg(j) = ξRd(j) = ξG(j) = ε.
D'où ξH(j) = (ξ$\Delta_0$(j)).
Or $\Delta_0$(j) = ε.
Donc H(j) = ε.
Ce qui induit que Σ(3, 1) = Σ(3#) = Σ(3#, 1) = Σ(4) = ε.
Et donc pour tout r ≥ 4, Σ(r, 1) = Σ(r#) = Σ(r#, 1) = Σ(r+1) = ε.
La propriété $R_1$ est terminée.
En particulier, DNJ(Ω, 1) = {(1#, 1), (2), (2, 1), (2#, 1), (3)}, avec DNJU(Ω, 1) = {(1#, 1), (2, 1), (2#, 1)}.

**4.7.5. Propriété $R_2$ :**
$Y_1$ = (1#, 1), $Y_2$ = (2, 1), et $Y_3$ = (2#, 1) sont les trois *Stratinos* unitaires de DNJU(Ω, 1).
Pour n = 1 : sΣ(1#, 1) = (d), sΣ(2, 1) = (k) et sΣ(2#, 1) = (i).
<u>A1 à A2</u> :
$\overline{B}$(d) = (c, g, e, b, c), avec d = zouc(v), médiane(d) = non, baou(v) = (b, c) et caouly(v) = (f, b).
Donc gh(d) = baou(v) = (b, c).
D'où $\overline{Br}$(d) = $\overline{B}$°(d) = (b, c, g, e, b).

$\overline{B}$(k) = (c, g, c), avec k = zouc(u), médiane (k) = non, baou(u) = (c, g) et caouly(u) = (g, c).
Donc gh(k) = baou(u) = (c, g).
D'où $\overline{Br}$(k) = $\overline{B}$°(k) = (c, g, c).

$\overline{B}$(i) = (f, g, h, j, f), avec i = zouc(w), médiane(i) = oui, baou(i) = (g, j) et caouly(i) = (j, f).
Donc gh(i) = (f, g).
D'où $\overline{Br}$(i) = $\overline{B}$°(i) = (f, g, h, j, f).

<u>A3 à A5</u> :
π(v) = {b, c, f, g} donc $\overline{\overline{Br}}$(d) = ($d_0$, $r_1$, $d_1$), avec $d_0 = (b, c, g) \in \pi(v)^*$, $d_1 = (b) \in \pi(v)^*$, et $r_1 = (e) \in (\Omega \backslash \pi(v))^*$.
D'où ν(d) = 1, ξg(b) = ξd(b) = ε, ξ$\Delta_1$(d) = (b), et ξ$R_1$(d) = (e).
Comme médiane(d) = non, ξ$\Delta_0$(d) = $d_0$ = (b, c, g).
Et donc ξH(d) = (ξ$\Delta_0$(d), ξ$R_1$(d), ξ$\Delta_1$(d)).

$\pi(u) = \{c, g\}$ donc $\overline{\overline{\text{Br}}}(k) = (d_0)$, avec $d_0 = (c, g, c) \in \pi(u)^*$.
D'où $\nu(k) = 0$, $\xi g(k) = \xi d(b) = \varepsilon$, et comme médiane(k) = non, $\xi\Delta_0(k) = (c, g, c)$.
Et donc $\xi H(k) = (\xi\Delta_0(d))$.

$\pi(w) = \{f, g, j\}$ donc $\overline{\overline{\text{Br}}}(i) = (d_0, r_1, d_1)$, avec $d_0 = (f, g) \in \pi(w)^*$, $d_1 = (j, f) \in \pi(w)^*$, et $r_1 = (h) \in (\Omega\backslash\pi(w))^*$.
D'où $\nu(i) = 1$, $\xi g(i) = \xi d(i) = \varepsilon$, $\xi\Delta_1(i) = (j, f)$, et $\xi R_1(i) = (h)$.
Comme médiane(i) = oui, alors $\xi\Delta_0(i)$ est le préfixe de $d_0$ vérifiant $d_0 = (f, \xi\Delta_0(i))$. Donc $\xi\Delta_0(i) = (g)$.
Et donc $\xi H(i) = (\xi\Delta_0(i), \xi R_1(i), \xi\Delta_1(i))$.

<u>A6</u> :
baou(v) = (b, c) et caouly(v) = (f, b).
$\xi\Delta_0(d) = (b, c, g)$, donc on est dans le 2$^{\text{ème}}$ cas : $\Delta_0(d)$ = (ʋ, g).
$\xi\Delta_1(d) = (b)$, donc on est dans le 4$^{\text{ème}}$ cas : $\Delta_1(d) = \varepsilon$.

baou(u) = (c, g) et caouly(u) = (g, c).
$\xi\Delta_0(k) = (c, g, c)$, donc on est dans le 1$^{\text{er}}$ cas : $\Delta_0(k)$ = (ʋ, ᴖ).

baou(w) = (g, j) et caouly(w) = (j, f).
$\xi\Delta_0(i) = (g)$, donc on est dans le 4$^{\text{ème}}$ cas : $\Delta_0(i) = \varepsilon$.
$\xi\Delta_1(i) = (j, f)$, donc on est dans le 4$^{\text{ème}}$ cas : $\Delta_1(i)$ = (ᴖ).

<u>A7</u> :
$\xi R_1(d) = (s_0)$, avec $s_0 = (e) \in v^*$.
Donc $T_1(d) = \varepsilon$, et $R_1(d) = S_1(d) = (s_0$, ɷ) = (e, ɷ).
De plus, binôme(e) = ʎʔ , médiane(e) = non, et $\varphi(e) = d$.
Comme g est le dernier terme de $\xi\Delta_0(d)$ et b le premier terme de $\xi\Delta_1(d)$, on en déduit que gh(e) = (d, g) et dh(e) = (b, d).

$\xi R_1(i) = (s_0)$, avec $s_0 = (h) \in w^*$.
Donc $T_1(i) = \varepsilon$, et $R_1(i) = S_1(i) = (s_0$, ɷ) = (h, ɷ).
De plus, binôme(h) = ʎʔ , médiane(h) = non, et $\varphi(h) = i$.
Comme g est le dernier terme de $\xi\Delta_0(i)$ et j le premier terme de $\xi\Delta_1(i)$, on en déduit que gh(h) = (i, g) et dh(h) = (j, i).

<u>A8 à A9</u> :
$H(d) = (\Delta_0(d), R_1(d), \Delta_1(d))$ = (ʋ, g, e, ɷ), $T(d) = T_1(d) = \varepsilon$, $S(d) = S_1(d)$ = (e, ɷ). D'où :
$\Sigma(1\#, 1, 1) = \varepsilon$,
$\Sigma(1\#, 1\#) = \varepsilon$,
$\Sigma(1\#, 1\#, 1) = (T(d)) = \varepsilon$, et
$\Sigma(1\#, 2) = (S(d))$ = (e, ɷ).
binôme(e) = ʎʔ , médiane(e) = non, et $\varphi(e) = d$, avec gh(e) = (d, g) et dh(e) = (b, d).

$H(k) = (\Delta_0(k))$ = (ʋ, ᴖ), $T(h) = S(h) = \varepsilon$. D'où :
$\Sigma(2, 1, 1) = \Sigma(2, 1\#) = \Sigma(2, 1\#, 1) = \Sigma(2, 2) = \varepsilon$.
Et donc pour tout $r \geq 2$, $\Sigma(2, r, 1) = \Sigma(2, r\#) = \Sigma(2, r\#, 1) = \Sigma(2, r+1) = \varepsilon$.

$H(i) = (\Delta_0(i), R_1(i), \Delta_1(i))$ = (h, ɷ, ᴖ), $T(i) = T_1(\iota) = \varepsilon$, $S(i) = S_1(\iota)$ = (h, ɷ). D'où :
$\Sigma(2\#, 1, 1) = \varepsilon$,
$\Sigma(2\#, 1\#) = \varepsilon$,
$\Sigma(2\#, 1\#, 1) = (T(i)) = \varepsilon$, et
$\Sigma(2\#, 2) = (S(i))$ = (h, ɷ).
binôme(h) = ʎʔ , médiane(h) = non, et $\varphi(h) = i$, avec gh(h) = (i, g) et dh(h) = (f, i).

Pour n = 2 : sΣ(1#, 2) = (e) et sΣ(2#, 2) = (h).
Il est évident que sΣ(1#, 2) = (e) et sΣ(2#, 2) vérifient la propriété P4.
$\overline{B}$(e) = (g, f, b, d, g), avec gh(e) = (d, g) et dh(e) = (b, d).
Donc $\overline{\overline{B°}}$(e) = (d, g, f, b, d) et $\overline{\overline{Br}}$(e) = (d, g, f, b, d).
$\overline{B}$(h) = (g, j, i, g), avec gh(h) = (i, g) et dh(h) = (j, i).
Donc $\overline{\overline{B°}}$(h) = (i, g, j, i) et $\overline{\overline{Br}}$(h) = (i, g, j, i).

<u>B1 à B2</u> :
π(v) = {b, c, f, g} donc $\overline{\overline{Br}}$(e) = ($r_0$, $d_0$, $r_1$), avec $d_0$ = (g, f, b) ∈ π(v)*, et $r_0$ = $r_1$ = (d) ∈ (Ω\π(v))*.
D'où ν(e) = 0, ξg(e) = ξd(e) =(d), et ξ$\Delta_0$(e) = (g, f, b).
ξg(e) et ξd(e) sont constitués d'un seul terme, donc ξRg(e) = ξRd(e) = ξG(e) = ε.
D'où ξH(e) = (ξ$\Delta_0$(e)).

π(w) = {f, g, j} donc $\overline{\overline{Br}}$(h) = ($r_0$, $d_0$, $r_1$), avec $d_0$ = (g, j) ∈ π(v)*, et $r_0$ = $r_1$ = (i) ∈ (Ω\π(v))*.
D'où ν(h) = 0, ξg(h) = ξd(h) =(i), et ξ$\Delta_0$(h) = (g, j).
ξg(eh et ξd(h) sont constitués d'un seul terme, donc ξRg(h) = ξRd(h) = ξG(h) = ε.
D'où ξH(h) = (ξ$\Delta_0$(h)).

<u>B3</u> :
baou(v) = (b, c) et caouly(v) = (f, b).
ξ$\Delta_0$(e) = (g, f, b), donc on est dans le 3$^{ème}$ cas : $\Delta_0$(e) = (f, ɷ).

baou(w) = (g, j) et caouly(w) = (j, f).
ξ$\Delta_0$(h) = (g, j), donc on est dans le 2$^{ème}$ cas : $\Delta_0$(h) = (ʋ).

<u>B8 à B9</u> :
H(e) = $\Delta_0$(e) = (f, ɷ), avec T(e) = S(e) = G(e) = ΔG(e) = ε. D'où :
Σ(1#, 2, 1) = Σ(1#, 2, 1#) = Σ(1#, 2, 1#, 1) = Σ(1#, 2, 2) = ε.
Et donc pour tout r ≥ 2, Σ(1#, 2, r, 1) = Σ(1#, 2, r#) = Σ(1#, 2, r#, 1) = Σ(1#, 2, r+1) = ε.

H(h) = $\Delta_0$(h) = (ʋ), avec T(h) = S(h) = G(h) = ΔG(h) = ε. D'où :
Σ(2#, 2, 1) = Σ(2#, 2, 1#) = Σ(2#, 2, 1#, 1) = Σ(2#, 2, 2) = ε.
Et donc pour tout r ≥ 2, Σ(2#, 2, r, 1) = Σ(2#, 2, r#) = Σ(2#, 2, r#, 1) = Σ(2#, 2, r+1) = ε.

La propriété $R_2$ est terminée.
En particulier, DNJ(Ω, 2) = {(1#, 2), (2#, 2)}, avec DNJU(Ω∪, 2) = ∅.

**4.7.6. Coloration *Stratojasse* (Figure 3) :**
NJ(Ω) = DNJ(Ω, 0) ∪ DNJ(Ω, 1) ∪ DNJ(Ω, 2) = {ε, (1), (1#, 1), (2), (2, 1), (2#, 1), (3), (1#, 2), (2#, 2)}.
Soit $c_0$, $c_1$, $c_2$ et $c_3$ quatre couleurs différentes : $c_0$ = blanc, $c_2$ = bleu clair, $c_1$ = gris clair et $c_3$ = gris foncé.

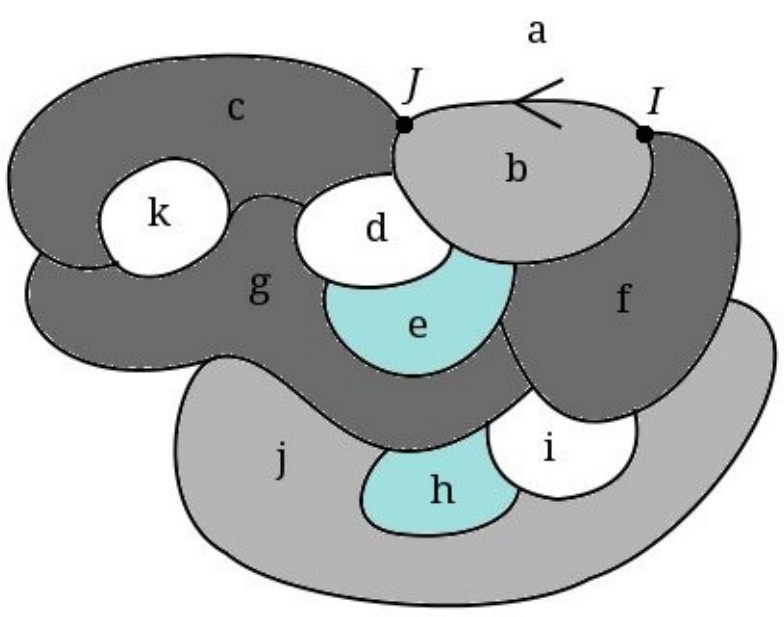

Figure 3. Coloration *Stratojasse*.

sΣ(1) = (b) et sΣ(3) = (j) : b et j sont colorées en $c_1$.
sΣ(2) = (c, g, f) : c, g et f sont colorées en $c_3$.
sΣ(ε) = (a), sΣ(1#, 1) = (d), sΣ(2, 1) = (k), et sΣ(2#, 1) = (i) : d, k et i sont colorées en $c_0$.
sΣ(1#, 2) = (e) et sΣ(2#, 2) = (h) : e et h sont colorées en $c_2$.

### 4.7.7. Construction de MJ(Ω)

Pour t = 0, $MJ_0$ = LH($w_0$) = ⟪ b ⟫, avec $w_0$ = a.

Pour t = 1 : $MJ_1$ = ⟪ LH(b) ⟫, avec $w_1$ = b.
binôme(b) = ⸢⸣ et H(b) = (c, d, •, ɷ), donc LH(b) = ⸢ c d • ɷ ⸣.
D'où $MJ_1$ = ⟪ ⸢ c d • ɷ ⸣ ⟫.

Pour t = 2 : $MJ_2$ = ⟪ ⸢ LH(c) d • ɷ ⸣ ⟫, avec $w_2$ = c.
binôme(c) = ⸤⸥ et H(c) = (k), donc LH(c) = ⸤ k ⸥ .
D'où $MJ_2$ = ⟪ ⸢ ⸤ k ⸥ d • ɷ ⸣ ⟫.

Pour t = 3 : $MJ_3$ = ⟪ ⸢ ⸤ LH(k) ⸥ d • ɷ ⸣ ⟫, avec $w_3$ = k.
binôme(k) = ⦃⦄ et H(k) = (ʊ, ᴖ), donc LH(k) = ⦃ ʊ ᴖ ⦄.
D'où $MJ_3$ = ⟪ ⸢ ⸤ ⦃ ʊ ᴖ ⦄ ⸥ d • ɷ ⸣ ⟫.

Pour t = 4 : $MJ_4$ = ⟪ ⸢ ⸤ ⦃ ʊ ᴖ ⦄ ⸥ LH(d) • ɷ ⸣ ⟫, avec $w_4$ = d.
binôme(d) = ⦑⦒ et H(d) = (ʊ, g, e, ɷ), donc LH(d) = ⦑ ʊ g e ɷ ⦒.
D'où $MJ_4$ = ⟪ ⸢ ⸤ ⦃ ʊ ᴖ ⦄ ⸥ ⦑ ʊ g e ɷ ⦒ • ɷ ⸣ ⟫.

Pour t = 5 : $MJ_5$ = ⟪ ⸢ ⸤ ⦃ ʊ ᴖ ⦄ ⸥ ⦑ ʊ LH(g) e ɷ ⦒ • ɷ ⸣ ⟫, avec $w_5$ = g.
binôme(g) = ⦗⦘ et H(g) = (j, i), donc LH(g) = ⦗ j i ⦘ .
D'où $MJ_5$ = ⟪ ⸢ ⸤ ⦃ ʊ ᴖ ⦄ ⸥ ⦑ ʊ ⦗ j i ⦘ e ɷ ⦒ • ɷ ⸣ ⟫.

Pour t = 6 : $MJ_6$ = ⟪ ⸢ ⸤ ⦃ ʊ ᴖ ⦄ ⸥ ⦑ ʊ ⦗ LH(j) i ⦘ e ɷ ⦒ • ɷ ⸣ ⟫, avec $w_6$ = j.
binôme(j) = ⸤⸥ et H(j) = ε, donc LH(j) = ⸤⸥ .
D'où $MJ_6$ = ⟪ ⸢ ⸤ ⦃ ʊ ᴖ ⦄ ⸥ ⦑ ʊ ⦗ ⸤⸥ i ⦘ e ɷ ⦒ • ɷ ⸣ ⟫.

Pour t = 7 : $MJ_7$ = ⟪ ⸢ ⸤ ⦃ ʊ ᴖ ⦄ ⸥ ⦑ ʊ ⦗ ⸤⸥ LH(i) ⦘ e ɷ ⦒ • ɷ ⸣ ⟫, avec $w_7$ = i.
binôme(i) = ⦑⦒ et H(i) = (h, ɷ, ᴖ), donc LH(i) = ⦑ h ɷ ᴖ ⦒.
D'où $MJ_7$ = ⟪ ⸢ ⸤ ⦃ ʊ ᴖ ⦄ ⸥ ⦑ ʊ ⦗ ⸤⸥ ⦑ h ɷ ᴖ ⦒ ⦘ e ɷ ⦒ • ɷ ⸣ ⟫.

Pour t = 8 : $MJ_8$ = ⟪ ⸢ ⸤ ⦃ ʊ ᴖ ⦄ ⸥ ⦑ ʊ ⦗ ⸤⸥ ⦑ LH(h) ɷ ᴖ ⦒ ⦘ e ɷ ⦒ • ɷ ⸣ ⟫, avec $w_8$ = h.
binôme(h) = ⸤⸥ et H(h) = (ʊ), donc LH(h) = ⸤ ʊ ⸥ .
D'où $MJ_8$ = ⟪ ⸢ ⸤ ⦃ ʊ ᴖ ⦄ ⸥ ⦑ ʊ ⦗ ⸤⸥ ⦑ ⸤ ʊ ⸥ ɷ ᴖ ⦒ ⦘ e ɷ ⦒ • ɷ ⸣ ⟫.

Pour t = 9 : $MJ_9$ = ⟪ ⸢ ⸤ ⦃ ʊ ᴖ ⦄ ⸥ ⦑ ʊ ⦗ ⸤⸥ ⦑ ⸤ ʊ ⸥ ɷ ᴖ ⦒ ⦘ LH(e) ɷ ⦒ • ɷ ⸣ ⟫, avec $w_9$ = e.
binôme(e) = ⸤⸥ et H(e) = (f, ᴖ), donc LH(e) = ⸤ f ᴖ ⸥ .
D'où $MJ_9$ = ⟪ ⸢ ⸤ ⦃ ʊ ᴖ ⦄ ⸥ ⦑ ʊ ⦗ ⸤⸥ ⦑ ⸤ ʊ ⸥ ɷ ᴖ ⦒ ⦘ ⸤ f ᴖ ⸥ ɷ ⦒ • ɷ ⸣ ⟫.

Pour t = 10 : $MJ_{10}$ = ⟪ ⸢ ⸤ ⦃ ʊ ᴖ ⦄ ⸥ ⦑ ʊ ⦗ ⸤⸥ ⦑ ⸤ ʊ ⸥ ɷ ᴖ ⦒ ⦘ ⸤ LH(f) ᴖ ⸥ ɷ ⦒ • ɷ ⸣ ⟫, avec $w_{10}$ = f.
binôme(f) = ⦗⦘ et H(f) = (•, ɷ), donc LH(f) = ⦗ • ɷ ⦘ .
D'où $MJ_{10}$ = ⟪ ⸢ ⸤ ⦃ ʊ ᴖ ⦄ ⸥ ⦑ ʊ ⦗ ⸤⸥ ⦑ ⸤ ʊ ⸥ ɷ ᴖ ⦒ ⦘ ⸤ ⦗ • ɷ ⦘ ᴖ ⸥ ɷ ⦒ • ɷ ⸣ ⟫.

$MJ_{10} \in$ AJ*, donc MJ(Ω) = $MJ_{10}$ = ⟪⸢⸤⦃ʊᴖ⦄⸥ ⦑ʊ⦗⸤⸥ ⦑⸤ʊ⸥ ɷᴖ⦒⦘ ⸤⦗•ɷ⦘ᴖ⸥ ɷ⦒•ɷ⸣⟫.

# 5
# Les *Mots Jassologiques Complexes*

## 5.1. L'alphabet *Jassologique* de *Catherine Mouet*

On a défini l'Alphabet *Jassologique* AJ = {[zoucarei], [zouccheirio], [mirounda], [miroundo], [trougouna], [trougouno], [trecouna], [trecouno], [stouna], [stoun], [louna], [loun], [couno], [chio], [baou], [caouly]} :
Les 12 caractères ouvrants et fermants des 6 binômes : [zoucarei][zouccheirio], [mirounda][miroundo], [trougouna][trougouno], [trecouna][trecouno], [stouna][stoun], [louna][loun],
Les 4 monômes : [couno], [chio], [baou], [caouly].
D'autre part, on distinguera :
- les binômes et monômes de type zouc : [zoucarei][zouccheirio], [mirounda][miroundo], [trougouna][trougouno], [trecouna][trecouno], [couno],
- les binômes et monômes de type stoun : [stouna][stoun], [chio],
- les binômes et monômes de type loun : [louna][loun], [baou], [caouly].

| | caractères ouvrants | | caractères fermants | | monômes | |
|---|---|---|---|---|---|---|
| zouc | [zoucarei] | zoucarei | [zouccheirio] | zouccheirio | | |
| | [mirounda] | mirounda | [miroundo] | miroundo | | |
| | [trougouna] | trougouna | [trougouno] | trougouno | | |
| | [trecouna] | trecouna | [trecouno] | trecouno | [couno] | couno |
| | | | | | | |
| stoun | [stouna] | stouna | [stoun] | stoun | [chio] | chio |
| | | | | | | |
| loun | [louna] | louna | [loun] | loun | [baou] | baou |
| | | | | | [caouly] | caouly |

## 5.2. Regroupements

### 5.2.1. *Trouglyre* :
Une *Trouglyre* G est un mot sur AJ qui est une puissance de « [trougouna][trougouno] »: il existe un entier r tel que
G = $([trougouna][trougouno])^r$, avec T = ε si r = 0.

### 5.2.2. *Stougammon* :
Un *Stougammon* Γ est un mot sur AJ qui est une puissance de « [stouna][stoun] »: il existe un entier r tel que
Γ = $([stouna][stoun])^r$, avec Γ = ε si r = 0.

### 5.2.3. *Trénagatte* :
Une *Trénagatte* T est un mot sur AJ qui est une puissance de « [trecouna][trecouno] [couno] »: il existe un entier r tel que
G = $([trecouna][trecouno] [couno])^r$, avec T = ε si r = 0.

### 5.2.4. *Lounafan* :
Un *Lounafan* π est un mot sur AJ tel que π = [baou] $([louna][loun])^r$ [caouly], avec π = [baou][caouly] si r = 0. Dans ce cas, on dira que π est un *Lounafan* simple.

### 5.2.5. *Lounagatte* :
Une *Lounagatte* L est un mot sur AJ tel que L = $\pi_1 \pi_2 \ldots \pi_n$, avec n ≥ 1, où pour tout k, 1 ≤ k ≤ n, $\pi_k$ est un *Lounafan*.
Et si pour tout k, 1 ≤ k ≤ n, $\pi_k$ est un *Lounafan* simple, on dira que L est une *Lounagatte* simple.

### 5.2.6. *Ramajo* , d-*Ramajo* et g-*Ramajo* :
Un *Ramajo* R est un mot sur AJ tel que R = [stouna][stoun] [chio] ou bien R = [stouna][stoun] $\rho_1 \rho_2 \ldots \rho_n$ [chio], avec n ≥ 1, où pour tout k, 1 ≤ k ≤ n, $\rho_k$ = [stouna][stoun] ou [trecouna][trecouno] [couno].
Ce qui induit que R se décompose de façon unique R = $\Gamma_0\ T_1\ \Gamma_1\ T_2\ \Gamma_2 \ldots T_p\ \Gamma_p$ [chio], avec p ≥ 0, où $\Gamma_0$ est est un *Stougammon* non vide et si p ≥ 1, pour tout 1 ≤ i ≤ p, $T_i$ est une *Trénagatte* non vide et $\Gamma_i$ est un *Stougammon* non vide, sauf pour $\Gamma_p$ qui peut éventuellement être vide.

Un d-*Ramajo* $R_d$ est un mot sur AJ tel que $R_d$ = [glyphe] $\rho_1\ \rho_2\ ...\ \rho_n$, avec $n \geq 1$, où $\rho_n = \varepsilon$ ou [glyphe], et si $n \geq 2$, pour tout k, $1 \leq k \leq n - 1$, $\rho_k$ = [glyphe] ou [glyphe] [glyphe].
Ce qui induit que $R_d$ se décompose de façon unique $R_d = \Gamma_0\ T_1\ \Gamma_1\ T_2\ \Gamma_2 \ldots T_p\ \Gamma_p\ \rho_n$, avec $p \geq 0$, où $\Gamma_0$ est est un *Stougammon* non vide et si $p \geq 1$, pour tout $1 \leq i \leq p$, $T_i$ est une *Trénagatte* non vide et $\Gamma_i$ est un *Stougammon* non vide, sauf pour $\Gamma_p$ qui peut éventuellement être vide.

Un g-*Ramajo* $R_g$ est un mot sur AJ tel que $R_g = \theta_0\ \theta_1\ \theta_2\ ...\ \theta_m$ [glyphe], avec $m \geq 0$, où $\theta_0 = \varepsilon$ ou [glyphe], et si $n \geq 1$, pour tout k, $1 \leq k \leq m$, $\rho_k$ = [glyphe] ou [glyphe] [glyphe].
Ce qui induit que $R_g$ se décompose de façon unique $R_g = \rho_0\ \Gamma_0\ T_1\ \Gamma_1\ T_2\ \Gamma_2 \ldots T_p\ \Gamma_p$ [glyphe], avec $p \geq 0$, où pour tout i, $0 \leq i \leq p$, $\Gamma_i$ est un *Stougammon* non vide, sauf pour $\Gamma_0$ et $\Gamma_p$ qui peuvent éventuellement être vides. Et si $p \geq 1$, pour tout i, $1 \leq i \leq p$, $T_i$ est une *Trénagatte* non vide.

Alors pour tout d-*Ramajo* $R_d$ et tout g-*Ramajo* $R_g$, $R = R_d\ R_g$ est un *Ramajo* si et seulement si $\rho_n = \theta_0 = \varepsilon$ ou bien $\rho_n$ = [glyphe] et $\theta_0$ = [glyphe].

**5.2.7. *Stratajo* , d-*Stratajo* et g-*Stratajo* :**
Un *Stratajo* S est un mot sur AJ tel que S = [glyphe] [glyphe] ou bien S = [glyphe] $\rho_1\ \rho_2\ ...\ \rho_n$ [glyphe], avec $n \geq 1$, où pour tout k, $1 \leq k \leq n$, $\rho_k$ = [glyphe] ou [glyphe] ([glyphe] )$^r$ [glyphe], c'est-à-dire un Lounafan.
Ce qui induit que S se décompose de façon unique $S = \Gamma_0\ L_1\ \Gamma_1\ L_2\ \Gamma_2 \ldots L_p\ \Gamma_p$ [glyphe], avec $p \geq 0$, où $\Gamma_0$ est est un *Stougammon* non vide et si $p \geq 1$, pour tout $1 \leq i \leq p$, $L_i$ est une *Lounagatte* non vide et $\Gamma_i$ est un *Stougammon* non vide, sauf pour $\Gamma_p$ qui peut éventuellement être vide.

Un d-*Stratajo* $S_d$ est un mot sur AJ tel que $S_d$ = [glyphe] ou bien $S_d$ = [glyphe] $\rho_1\ \rho_2\ ...\ \rho_n$, avec $n \geq 1$, où pour tout k, $1 \leq k \leq n$, $\rho_k$ = [glyphe] ou [glyphe] ([glyphe] )$^r$ [glyphe].
Ce qui induit que $S_d$ se décompose de façon unique $S_d = \Gamma_0\ L_1\ \Gamma_1\ L_2\ \Gamma_2 \ldots L_p\ \Gamma_p$, avec $p \geq 0$, où $\Gamma_0$ est est un *Stougammon* non vide et si $p \geq 1$, pour tout $1 \leq i \leq p$, $L_i$ est une *Lounagatte* non vide et $\Gamma_i$ est un *Stougammon* non vide, sauf pour $\Gamma_p$ qui peut éventuellement être vide.

Un g-*Stratajo* $S_g$ est un mot sur AJ tel que $S_g$ = [glyphe] ou bien $S_g = \rho_1\ \rho_2\ ...\ \rho_n$ [glyphe], avec $n \geq 1$, où pour tout k, $1 \leq k \leq n$, $\rho_k$ = [glyphe] ou [glyphe] ([glyphe] )$^r$ [glyphe].
Ce qui induit que $S_g$ se décompose de façon unique $S_g = \Gamma_0\ L_1\ \Gamma_1\ L_2\ \Gamma_2 \ldots L_p\ \Gamma_p$ [glyphe], avec $p \geq 0$, où pour tout $0 \leq i \leq p$, $\Gamma_i$ est un *Stougammon* non vide, sauf pour $\Gamma_0$ et $\Gamma_p$ qui peuvent éventuellement être vides, et si $p \geq 1$, pour tout $1 \leq i \leq p$, $L_i$ est une *Lounagatte* non vide.

Alors pour tout d-*Stratajo* $S_d$ et tout g-*Stratajo* $S_g$, $S = S_d\ S_g$ est un *Stratajo*.

**5.3. Règle n°1 : la *Règle d'Emboîtement***

**5.3.1. La *Règle d'Emboîtement* :**
Soit $M = x_1\ x_2 \ldots x_\eta$ un mot sur AJ, avec $\eta \geq 1$ sa longueur.
Alors M vérifie la Règle d'Emboîtement si et seulement si
1) M contient au-moins 3 caractères ouvrants : notons $\alpha_0, \alpha_1, \ldots, \alpha_N$, les entiers tels que $x_{\alpha p}$ soit un caractère ouvrant pour tout p, $0 \leq p \leq N$, avec $N \geq 2$ et $1 \leq \alpha_0 < \alpha_1 < \ldots < \alpha_N \leq \eta$.
2) le mot sur J = {[glyphe], [glyphe]}, obtenu à partir de M en supprimant tous les monômes et en identifiant les caractères ouvrants à « [glyphe] » et les caractères fermants à « [glyphe] », est un *Mot Jassologique Simple* : pour tout p, $1 \leq p \leq N$, il existe donc un et un seul entier $\beta_p$, tel que $x_{\beta p}$ soit le caractère fermant correspondant à $x_{\alpha p}$ dans M. On notera $J_p = (\alpha_p, \beta_p)$.
3) $\alpha_0 = 1$, $\beta_0 = \eta$, et $x_{\alpha 0}\ x_{\beta 0}$ = [glyphe].
4) $\alpha_1 = 2$, $\beta_1 = \eta - 1$, et $x_{\alpha 1}\ x_{\beta 1}$ = [glyphe].
5) pour tout p, $2 \leq p \leq N$, $x_{\alpha p}\ x_{\beta p}$ = [glyphe] ou [glyphe] ou [glyphe] ou [glyphe] .

**5.3.2. L'ensemble $\overline{AJ^*}$ :**
Par définition, $\overline{AJ^*}$ sera le sous-ensemble de AJ* des Mots sur AJ qui vérifient la Règle d'emboîtement.

### 5.4. Le *Dallajascar* relatif à un mot $M \in \overline{AJ^*}$

Soit $M = x_1 x_2 \dots x_\eta \in \overline{AJ^*}$, avec $\eta \geq 1$ sa longueur.
En reprenant les notations précédentes, N+1 est le nombre de caractères ouvrants, et pour tout p, $0 \leq p \leq N$, on a défini $J_p = (\alpha_p, \beta_p)$, avec $x_{\alpha p}\, x_{\beta p}$ = ℭ𝔇 ou ⌜⌝ ou ❛❜ ou ✗❟ ou ⟅⟆ ou ⦅⦆ .

#### 5.4.1. L'ensemble E(M) :

Notons E(M) l'ensemble des entiers $\gamma$, $1 \leq \gamma \leq \eta$, tels que $x_\gamma$ soit un caractère ouvrant ou un monôme, mais pas un caractère fermant. Son cardinal est donc exactement $\eta - (N+1)$. C'est un ensemble totalement ordonné par la relation d'ordre sur les entiers.

#### 5.4.2. La fonction $L : E(M)^* \to AJ^*$ :

Pour tout $\gamma \in E(M)$, si $x_\gamma$ est un monôme, on posera $L(\gamma) = x_\gamma$, sinon il existe p, $0 \leq p \leq N$, tel que $\gamma = \alpha_p$ et donc on posera $L(\gamma) = x_{\alpha p}\, x_{\beta p}$.
Pour tout nesile $s = (\gamma_1, \gamma_2, \dots, \gamma_r)$ de E(M), on posera $L(s) = L(\gamma_1)\, L(\gamma_2) \dots L(\gamma_r)$.

#### 5.4.3. Le *Dallajascar* $(\varphi,H)$ sur E(M) :

Soit $\varphi : E(M) \to E(M) \cup \{\varepsilon\}$ et $H : E(M) \to E(M)^*$ définies par :
i) $\varphi(\alpha_0) = \varepsilon$,
ii) pour tout p, $1 \leq p \leq N$, $\varphi(\alpha_p) = \alpha_q$, où q est le plus grand entier, $0 \leq q \leq p - 1$, tel que
$1 \leq \alpha_q < \alpha_p < \beta_p < \beta_q \leq \eta$. Comme $\alpha_0 = 1$ et $\beta_0 = \eta$, l'entier q est bien défini.
iii) pour tout $\gamma$, $1 \leq \gamma \leq \eta$, tel que $x_\gamma$ soit un monôme, $\varphi(\gamma) = \alpha_q$, où q est le plus grand entier tel que
$1 \leq \alpha_q < \gamma < \beta_q \leq \eta$. Comme $\alpha_0 = 1$ et $\beta_0 = \eta$, l'entier q est bien défini.
iv) pour tout $\gamma \in E(M)$, si $\gamma$ n'a pas d'antécédent par $\varphi$, alors $H(\gamma) = \varepsilon$, sinon $H(\gamma) = (\gamma_1, \gamma_2, \dots, \gamma_r)$, où $\gamma_1, \gamma_2, \dots, \gamma_r$ constituent l'ensemble des antécédents de $\gamma$ par $\varphi$, avec $\gamma_1 < \gamma_2 < \dots < \gamma_r$ selon la relation d'ordre définie sur E(M).
On vérifie que $(\varphi,H)$ constitue bien un *Dallajascar* sur E(M), $\alpha_0$ étant la racine.
De plus, la relation d'ordre définie par $(\varphi,H)$ est exactement la même que celle définie ci-dessus sur E(M) : pour tous éléments $\gamma$ et $\gamma'$ de E(M), $\gamma < \gamma' \Leftrightarrow \gamma <_{(\varphi,H)} \gamma'$.
En particulier, comme $\alpha_0 = 1$, $\beta_0 = \eta$, avec $x_{\alpha 0}\, x_{\beta 0}$ = ℭ𝔇, et $\alpha_1 = 2$, $\beta_1 = \eta - 1$, avec $x_{\alpha 1}\, x_{\beta 1}$ = ⌜⌝, on en déduit que $\varphi(\alpha_1) = \alpha_0$ et $H(\alpha_0) = (\alpha_1)$.

### 5.5. Règle n°2 : la *Règle des Stratojasses*

Soit $M = x_1 x_2 \dots x_\eta \in \overline{AJ^*}$, avec $\eta \geq 1$ sa longueur.
En reprenant les notations précédentes, N+1 est le nombre de caractères ouvrants, et pour tout p, $0 \leq p \leq N$, on a défini $J_p = (\alpha_p, \beta_p)$, avec $x_{\alpha p}\, x_{\beta p}$ = ℭ𝔇 ou ⌜⌝ ou ❛❜ ou ✗❟ ou ⟅⟆ ou ⦅⦆ .
Puis notons $E(M) = \{\gamma_1, \gamma_2, \dots, \gamma_{\eta-N-1}\}$, avec $1 = \gamma_1, < \gamma_2, < \dots, < \gamma_{\eta-N-1} < \eta$.
Rappelons que $\gamma_1 = \alpha_0 = 1$, avec $x_{\alpha 0}\, x_{\beta 0}$ = ℭ𝔇, $\gamma_2 = \alpha_1 = 2$, avec $x_{\alpha 1}\, x_{\beta 1}$ = ⌜⌝, et pour tout p, $2 \leq p \leq N$, $x_{\alpha p}\, x_{\beta p}$ = ❛❜ ou ✗❟ ou ⟅⟆ ou ⦅⦆ .
On a défini les *Stratinos Naturels* et les *Stratinos Décalés*.

#### 5.5.1. La fonction *Stratino*, notée $\sigma$ :

Définissons la fonction *Stratino* $\sigma : E(M) \to NJ^*$ par récurrence sur k, $1 \leq k \leq \eta - N - 1$, avec la propriété que si $x_{\gamma k}$ est un caractère ouvrant, alors $\sigma(\gamma_1)$ est un *Stratino Naturel*.
Pour k = 1 : posons $\sigma(\gamma_1) = \varepsilon$.
Pour k = 2 : posons $\sigma(\gamma_2) = (1)$.
Soit $k \geq 3$ et supposons que $\sigma(\gamma_r)$ soit bien défini pour tout r, $1 \leq r \leq k - 1$, avec la propriété énoncée.
Notons p l'entier tel que $\varphi(\gamma_k) = \alpha_p$, vérifiant par définition $2 \leq \alpha_p < \gamma_k$.
Par hypothèse de récurrence, $\sigma(\alpha_p)$ est bien défini, et comme $x_{\alpha p}$ est un caractère ouvrant, $\sigma(\alpha_p)$ est un *Stratino Naturel*.
Si $\sigma(\alpha_p) = \varepsilon$, alors on posera $\sigma(\gamma_k) = \varepsilon$.
Sinon il existe un et un seul *Stratino* X et un et un seul entier $n > 0$ tel que $\sigma(\alpha_p) = (X, n)$.

<u>1er cas :</u> $x_{\gamma k}$ = [symbole] :

- si n = 1, alors $\sigma(\gamma_k) = \varepsilon$.
- si $n \geq 2$, alors $\sigma(\gamma_k) = (X, n, 1)$.

<u>2ème cas :</u> $x_{\gamma k}$ = [symbole] ou [symbole], alors $\sigma(\gamma_k) = (X, n\#, 1)$.

<u>3ème cas :</u> $x_{\gamma k}$ = [symbole] ou [symbole], alors $\sigma(\gamma_k) = (X, n+1)$.

<u>4ème cas :</u> $x_{\gamma k}$ = [symbole] :

- si X est un *Stratino Naturel*, alors $\sigma(\gamma_k) = \varepsilon$.
- si X est un *Stratino Décalé*, avec X = (Y, m#), alors $\sigma(\gamma_k) = (Y, m+1)$.

<u>5ème cas :</u> $x_{\gamma k}$ = [symbole] ou [symbole] :

- si $X = \varepsilon$ ou bien si X est un *Stratino* unitaire, alors $\sigma(\gamma_k) = \varepsilon$.
- si X est un *Stratino Naturel* non unitaire, avec X = (Y, m), alors $\sigma(\gamma_k) = (Y, m\#)$.
- si X est un *Stratino Décalé*, avec X = (Y, m#), alors $\sigma(\gamma_k) = (Y, m+1)$.

Remarquons que dans tous les cas, si $\sigma(\gamma_k) \neq \varepsilon$, alors $\sigma(\alpha_p) < \sigma(\gamma_k)$ selon la relation d'ordre définie sur les *Stratinos*.

**5.5.2. L'ensemble NJ(M) des *Stratinos* de M :**

$NJ(M) = \{\sigma(\gamma) \mid \gamma \in E(M)\}$ est un ensemble ordonné qui permet de graduer l'*Echelle Jassologique* de M.

**5.5.3. La fonction *Stratojasse*, notée Σ :**

Définissons la fonction *Stratojasse* $\Sigma : NJ^* \rightarrow E(M)^*$ de la façon suivante :
Soit X un *Stratino* quelconque.
Si $X \notin NJ(M)$, alors $\Sigma(X) = \varepsilon$.
Sinon $\Sigma(X) = (\gamma_{k1}, \gamma_{k2}, \ldots, \gamma_{kr})$, avec $r \geq 1$, où $\gamma_{k1}, \gamma_{k2}, \ldots, \gamma_{kr}$ constituent l'ensemble des antécédents de X par σ, avec $\gamma_{k1} < \gamma_{k2} < \ldots < \gamma_{kr}$.
On vérifie alors les résultats suivants :
Si X = (1), alors r = 1 et $\gamma_{k1} = 2$, avec $x_2$ = [symbole],
Si X = (Z, m, 1), m étant un entier, alors $m \geq 2$ et pour tout t, $1 \leq t \leq r$, $x_{\gamma kt}$ = [symbole],
Si X = (Z, m#, 1), m étant un entier, alors pour tout t, $1 \leq t \leq r$, $x_{\gamma kt}$ = [symbole] ou [symbole],
Si X = (Y, n), n étant un entier $\geq 2$, alors pour tout t, $1 \leq t \leq r$, $x_{\gamma kt}$ = [symbole] ou [symbole] ou [symbole] ou [symbole] ou [symbole],
Si X = (Y, n#), n étant un entier $\geq 2$, alors pour tout t, $1 \leq t \leq r$, $x_{\gamma kt}$ = [symbole] ou [symbole],

**5.5.4. La *Règle des Stratojasses* :**

Soit $M \in \overline{AJ^*}$ et reprenons les notations ci-dessus.
On dira que M vérifie la *Règle des Stratojasses* si et seulement si $\Sigma(\varepsilon) = (\gamma_1)$, c'est-à-dire $\Sigma(\varepsilon) = (1)$ puisque $\gamma_1 = \alpha_0 = 1$.
Autrement dit, M vérifie la *Règle des Stratojasses* si et seulement si pour tout $\gamma \in E(M)$ différent de 1, $\sigma(\gamma) \neq \varepsilon$.

**5.6. Règle n°3 : La *Règle du Dallajascar***

Soit $M = x_1 x_2 \ldots x_\eta \in \overline{AJ^*}$, avec $\eta \geq 1$ sa longueur.
En reprenant les notations précédentes, N+1 est le nombre de caractères ouvrants, et pour tout p, $0 \leq p \leq N$, on a défini $J_p = (\alpha_p, \beta_p)$, avec $x_{\alpha p} x_{\beta p}$ = [symbole] ou [symbole] ou [symbole] ou [symbole] ou [symbole] ou [symbole] .
On a défini l'ensemble E(M) et le *Dallajascar* (φ,H) sur E(M), dont $\alpha_0 = 1$ est la racine.
Supposons que M vérifie la *Règle des Stratojasses*.

**5.6.1. Le sous-ensembles A(M) de E(M) et la fonction sΣ :**

Notons $A(M) = \{\alpha_1, \ldots, \alpha_N\}$, qui est un sous-ensemble de E(M).
Définissons la fonction $s\Sigma : NJ^* \rightarrow E(M)^*$ de la façon suivante :
Pour tout *Stratino* X, sΣ(X) sera par définition le *sous-nesile* de Σ(X) constitué de l'ensemble des éléments appartenant à A(M).

**5.6.2. Les sous-ensembles Cg(M) et Cd(M) de A(M) :**

Cg(M) est le sous-ensemble de A(M) défini de la façon suivante :
Soit p, $1 \leq p \leq N$, et notons $X = \sigma(\alpha_p)$. Alors $\alpha_p \in Cg(M)$ si et seulement si :

i) $x_{\alpha p}$ = ɼ ou ɤ ou ϰ, ou bien

ii) $x_{\alpha p}$ = ƛ ou ɮ et $\alpha_p$ est le premier terme de $s\Sigma(X)$, ou bien

iii) $x_{\alpha p}$ = ƛ ou ɮ, il existe $\alpha_q$, $q < p$, qui précède $\alpha_p$ dans $s\Sigma(X)$, et il existe $\gamma \in \Sigma(X)$, tels que $x_\gamma$ = ɷ et $\alpha_q < \gamma < \alpha_p$.

Symétriquement, Cd(M) est le sous-ensemble de A(M) défini de la façon suivante :
Soit p, $1 \leq p \leq N$, et notons $X = \sigma(\alpha_p)$. Alors $\alpha_p \in$ Cd(M) si et seulement si :

i) $x_{\alpha p}$ = ɼ ou ɤ ou ϰ, ou bien

ii) $x_{\alpha p}$ = ƛ ou ɮ et $\alpha_p$ est le dernier terme de $s\Sigma(X)$, ou bien

iii) $x_{\alpha p}$ = ƛ ou ɮ, il existe $\alpha_q$, $q > p$, qui suit $\alpha_p$ dans $s\Sigma(X)$, et il existe $\gamma \in \Sigma(X)$, tels que $x_\gamma$ = ɷ et $\alpha_p < \gamma < \alpha_q$.

**5.6.3. La fonction G : A(M) → E(M)* :**

Soit p, $1 \leq p \leq N$.
$G(\alpha_p)$ sera par définition le *sous-nesile* de $H(\alpha_p)$ constitué des éléments $\gamma$ tels que $x_\gamma$ = ɤ, avec $G(\alpha_p) = \varepsilon$ si $H(\alpha_p)$ ne contient pas d'élément $\gamma$ tels que $x_\gamma$ = ɤ.

**5.6.4. La fonction R : A(M) → E(M)* :**

Soit p, $1 \leq p \leq N$.
$R(\alpha_p)$ sera par définition le *sous-nesile* de $H(\alpha_p)$ constitué des éléments $\gamma$ tels que $x_\gamma$ = ϰ ou • ou ƛ ou ɷ, avec $R(\alpha_p) = \varepsilon$ si $H(\alpha_p)$ ne contient pas d'élément $\gamma$ tel que $x_\gamma$ = ϰ ou • ou ƛ ou ɷ.

**5.6.5. Les fonctions Rd et R' : A(M) → E(M)* :**

Soit p, $1 \leq p \leq N$.
Si $R(\alpha_p)$ ne contient aucun élément $\gamma$ tel que $x_\gamma$ = ɷ, alors $Rd(\alpha_p) = R(\alpha_p)$ et $R'(\alpha_p) = \varepsilon$.
Sinon notons $\gamma$ le plus grand entier appartenant à $R(\alpha_p)$ tel que $x_\gamma$ = ɷ.
Alors $Rd(\alpha_p)$ sera le suffixe de $R(\alpha_p)$ situé après $\gamma$ dans $R(\alpha_p)$, et $R'(\alpha_p)$ sera le préfixe de $R'(\alpha_p)$ dont $\gamma$ est le dernier terme.
Ce qui induit que $R(\alpha_p) = (R'(\alpha_p), Rd(\alpha_p))$.

**5.6.6. La fonction ν : A(M) → {0, 1, 2, 3, ...} et les fonctions Rg et $R_k$ : A(M) → E(M)* pour tout $k \geq 1$ :**

Soit p, $1 \leq p \leq N$, et notons $X = \sigma(\alpha_p)$.
Si $R'(\alpha_p) = \varepsilon$, alors on posera $\nu(\alpha_p) = 0$, $Rg(\alpha_p) = \varepsilon$ et $R_k(\alpha_p) = \varepsilon$ pour tout $k \geq 1$.
Supposons que $R'(\alpha_p) \neq \varepsilon$.
Alors $R'(\alpha_p)$ se décompose de façon unique $R'(\alpha_p) = (R'_0, \gamma_0, R'_1, \gamma_1, \ldots, R'_u, \gamma_u)$, où $\gamma_0, \gamma_1, \ldots, \gamma_u$, constituent l'ensemble des éléments de $R'(\alpha_p)$ tels que $x_{\gamma t}$ = ɷ pour tout t, avec $u \geq 0$. Ce qui induit que pour tout t, $0 \leq t \leq u$ et tout $\gamma \in R'_t$, $x_\gamma$ = ϰ ou • ou ƛ.

<u>1er cas :</u> $\alpha_p \in$ Cg(M) :
Alors $Rg(\alpha_p) = \varepsilon$, $\nu(\alpha_p) = u + 1$, et :

- pour tout k, $1 \leq k \leq \nu(\alpha_p)$, $R_k(\alpha_p) = (R'_{k-1}, \gamma_{k-1})$.
- pour tout $k \geq \nu(\alpha_p) + 1$, $R_k(\alpha_p) = \varepsilon$.

<u>2ème cas :</u> $\alpha_p \notin$ Cg(M) :
Alors $Rg(\alpha_p) = (R'_0, \gamma_0)$, $\nu(\alpha_p) = u$, et :

- si $\nu(\alpha_p) \geq 1$, pour tout k, $1 \leq k \leq \nu(\alpha_p)$, $R_k(\alpha_p) = (R'_k, \gamma_k)$.
- pour tout $k \geq \nu(\alpha_p) + 1$, $R_k(\alpha_p) = \varepsilon$.

**5.6.7. La *Règle du Dallajascar* :**

Soit M ∈ $\overline{AJ}$* qui vérifie la *Règle des Stratojasses* et reprenons les notations ci-dessus.
On a défini la fonction L : E(M)* → AJ*.
On dira que M vérifie la *Règle du Dallajascar* si et seulement si pour tout p, $1 \leq p \leq N$ :
1) si $\alpha_p \in$ Cd(M), alors $G(\alpha_p) = Rd(\alpha_p) = \varepsilon$.
2) si $Rd(\alpha_p) \neq \varepsilon$, alors $L(Rd(\alpha_p))$ est un *d-Ramajo*.
3) si $Rg(\alpha_p) \neq \varepsilon$, alors $L(Rg(\alpha_p))$ est un *g-Ramajo*,

4) s'il existe $\alpha_q$ qui précède $\alpha_p$ dans $s\Sigma(\sigma(\alpha_p)$, alors $R = L(Rd(\alpha_q))\ L(Rg(\alpha_p)) = \varepsilon$ ou bien R est un *Ramajo*.
5) si $\nu(\alpha_p) \geq 1$, alors pour tout k, $1 \leq k \leq \nu(\alpha_p)$, $L(R_k(\alpha_p))$ est un *Ramajo*.
6) $(Rd(\alpha_p), G(\alpha_p))$ est un suffixe de $H(\alpha_p)$.
7) $Rg(\alpha_p)$ est un préfixe de $H(\alpha_p)$.
8) pour tout $k \geq 1$, $R_k(\alpha_p)$ est un facteur de $H(\alpha_p)$.

**5.6.8. Les fonctions $\Delta_k$ : A(M) → E(M)* pour tout $k \geq 0$ :**
Soit $M \in \overline{AJ^*}$ qui vérifie la *Règle des Stratojasses* et la *Règle du Dallajascar*.
Et reprenons les notations ci-dessus.
Soit p, $1 \leq p \leq N$.
Comme $(Rd(\alpha_p), G(\alpha_p))$ est un suffixe de $H(\alpha_p)$, que $Rg(\alpha_p)$ est un préfixe de $H(\alpha_p)$, et que pour tout $k \geq 1$, $R_k(\alpha_p)$ est un facteur de $H(\alpha_p)$, on en déduit que $H(\alpha_p)$ se décompose de la façon suivante :
$H(\alpha_p) = (Rg(\alpha_p), d_0, R_1(\alpha_p), d_1, R_2(\alpha_p), d_2, \ldots, R_{\nu(\alpha p)}(\alpha_p), d_{\nu(\alpha p)}, Rd(\alpha_p), G(\alpha_p))$, avec
$H(\alpha_p) = (Rg(\alpha_p), d_0, Rd(\alpha_p), G(\alpha_p))$ si $\nu(\alpha_p) = 0$.
Pour tout k, $0 \leq p \leq \nu(\alpha_p)$, on posera $\Delta_k(\alpha_p) = d_k$, et pour tout $k \geq \nu(\alpha_p)+1$, $\Delta_k(\alpha_p) = \varepsilon$.
Comme $G(\alpha_p)$ est constitué de tous les éléments $\gamma$ de $H(\alpha_p)$ tels que $x_\gamma$ = ❦, et que $R(\alpha_p)$ est constitué de tous les éléments $\gamma$ de $H(\alpha_p)$ tels que $x_\gamma$ = ❧ ou • ou ʃ ou ɷ, alors les $\Delta_k(\alpha_p)$ sont constitués de tous les éléments $\gamma$ de $H(\alpha_p)$ tels que $x_\gamma$ = ʊ ou ᴖ ou ɬ.
En conclusion, on a donc :
$H(\alpha_p) = (Rg(\alpha_p), \Delta_0(\alpha_p), R_1(\alpha_p), \Delta_1(\alpha_p), R_2(\alpha_p), \Delta_2(\alpha_p), \ldots, R_{\nu(\alpha p)}(\alpha_p), \Delta_{\nu(\alpha p)}(\alpha_p), Rd(\alpha_p), G(\alpha_p))$, avec
$H(\alpha_p) = (Rg(\alpha_p), \Delta_0(\alpha_p), Rd(\alpha_p), G(\alpha_p))$ si $\nu(\alpha_p) = 0$.

## 5.7. Règle n°4 : La *Règle des Lounafans*

Soit $M = x_1\ x_2 \ldots x_\eta \in \overline{AJ^*}$, avec $\eta \geq 1$ sa longueur.
En reprenant les notations précédentes, N+1 est le nombre de caractères ouvrants, et pour tout p, $0 \leq p \leq N$, on a défini $J_p = (\alpha_p, \beta_p)$, avec $x_{\alpha p}\ x_{\beta p}$ = ℭ𝔇 ou ᒥᒧ ou ❦❦ ou ❧❧ ou ʃʅ ou ɬɭ .
On a défini l'ensemble E(M) et le *Dallajascar* $(\varphi, H)$ sur E(M), dont $\alpha_0 = 1$ est la racine.
Supposons que M vérifie la *Règle des Stratojasses* et la *Règle du Dallajascar*.

**5.7.1. L'ensemble *Zouc*(M) et la fonction *fan* :**
Posons $Zouc(M) = \{\gamma \in E(M) \mid x_\gamma$ = ❦ ou ❧$\}$.
Et définissons la fonction fan : Zouc(M) → E(M)* de la façon suivante :
Soit $\alpha \in Zouc(M)$.
Alors $\sigma(\alpha)$ est un *Stratino* unitaire.
<u>1er cas :</u> $x_\alpha$ = ❦ :
Alors il existe un *Stratino* X et un entier $n \geq 2$ tels que $\sigma(\alpha) = (X, n, 1)$.
De plus, $\sigma(\varphi(\alpha)) = (X, n)$ par construction.
Si $\Sigma(X, n\#) = \varepsilon$, ou bien si pour tout $\gamma \in \Sigma(X, n\#)$ $\alpha$ n'emboîte pas $\gamma$, alors $fan(\alpha) = \varepsilon$.
Sinon $fan(\alpha) = (\gamma_1, \gamma_2, \ldots, \gamma_r)$, avec $r \geq 1$, où $\gamma_1, \gamma_2, \ldots, \gamma_r$ constituent l'ensemble des éléments de $\Sigma(X, n\#)$ tels que $\alpha \ulcorner \gamma_t$ pour tout t, avec $\gamma_1 < \gamma_2 < \ldots < \gamma_r$.
Rappelons qu'alors $x_\alpha$ = ʊ ou ᴖ pour tout t.
<u>2ème cas :</u> $x_\alpha$ = ❧ :
Alors il existe un *Stratino* X et un entier $n \geq 1$ tels que $\sigma(\alpha) = (X, n\#, 1)$.
De plus, $\sigma(\varphi(\alpha)) = (X, n)$ par construction.
Si $\Sigma(X, n+1) = \varepsilon$, ou bien si pour tout $\gamma \in \Sigma(X, n+1)$ $\alpha$ n'emboîte pas $\gamma$, alors $fan(\alpha) = \varepsilon$.
Sinon $fan(\alpha) = (\gamma_1, \gamma_2, \ldots, \gamma_r)$, avec $r \geq 1$, où $\gamma_1, \gamma_2, \ldots, \gamma_r$ constituent l'ensemble des éléments de $\Sigma(X, n+1)$ tels que $\alpha \ulcorner \gamma_t$ pour tout t, avec $\gamma_1 < \gamma_2 < \ldots < \gamma_r$.
On vérifie qu'alors $x_\alpha$ = ʊ ou ɬ ou ᴖ pour tout t.

**5.7.2. Prosition :**
Soit p, $1 \le p \le N$.
Alors il existe un et un seul $\alpha \in$ Zouc(M) tel que pour tout $k \ge 0$, $\Delta_k(\alpha_p)$ soit un facteur de fan($\alpha$). En particulier, si $\alpha_p \in$ Zouc(M), alors $\alpha = \alpha_p$.

**5.7.3. Le sous-ensemble Zm(M) de Zouc(M) :**
Zm(M) est le sous-ensemble de Zouc(M) défini de la façon suivante :
Soit $\alpha_p \in$ Zouc(M) et notons q l'entier tel que $\alpha_q = \varphi(\alpha_p)$ : donc $\alpha_p \in H(\alpha_q)$.
Alors $\alpha_p \in$ Zm(M) si et seulement si $x_{\alpha p}$ = [symbole], $Rd(\alpha_q) \neq \varepsilon$, et $\alpha_p$ est le dernier terme de $Rd(\alpha_q)$.
D'après la *Règle du Dallajascar*, cela induit que si on note $\alpha_r$ le terme qui suit $\alpha_q$ dans $s\Sigma(\sigma(\alpha_q)$, alors $R = L(Rd(\alpha_q))$ $L(Rg(\alpha_r))$ est un *Ramajo*. Donc nécessairement, le premier terme de $L(Rg(\alpha_r))$ est $\bullet$.

**5.7.4. La *Règle des Lounafans* :**
Soit $M \in \overline{AJ}^*$ qui vérifie la *Règle des Stratojasses* et la *Règle du Dallajascar* et reprenons les notations ci-dessus.
On dira que M vérifie la *Règle des Lounafans* si et seulement si pour tout $\alpha \in$ Zouc(M) :
1) si $x_\alpha$ = [symbole], alors L(fan($\alpha$)) est un *Lounafan Simple*, c'est-à-dire L(fan($\alpha$)) = [symbole][symbole].
2) si $x_\alpha$ = [symbole], alors L(fan($\alpha$)) est un *Lounafan*, c'est-à-dire L(fan($\alpha$)) = [symbole] ([symbole][symbole] )$^u$ [symbole], avec $u \ge 0$.
3) si $\alpha \notin$ Zm(M), alors [symbole] est le premier terme de $L(\Delta_0(\alpha_p))$.

**5.8. Les *Mots Jassologiques Complexes***

Un *Mot Jassologique Complexe* sera par définition un mot su AJ qui vérifie les 4 règles : la *Règle d'Emboîtement*, la *Règle des Stratojasses*, la *Règle du Dallajascar*, et la *Règle des Lounafans*.

Soit $M = x_1\, x_2 \ldots x_\eta$ un *Mot Jassologique Complexe*, avec $\eta \ge 1$ sa longueur.
En reprenant les notations précédentes, N+1 est le nombre de caractères ouvrants, et pour tout p, $0 \le p \le N$, on a défini $J_p = (\alpha_p, \beta_p)$, avec $x_{\alpha p}\, x_{\beta p}$ = [symbole][symbole] ou [symbole][symbole] ou [symbole][symbole] ou [symbole][symbole] ou [symbole][symbole] ou [symbole][symbole] , l'ensemble E(M) et le *Dallajascar* ($\varphi$,H) sur E(M), dont $\alpha_0 = 1$ est la racine.
Nous allons définir les fonctions $\Delta G$, T et S : $A(M) \to E(M)^*$.

**5.8.1. La fonction $\Delta G : A(M) \to E(M)^*$ :**
Soit p, $1 \le p \le N$.
Si $G(\alpha_p) = \varepsilon$, alors $\Delta G(\alpha_p) = \varepsilon$.
Sinon notons $G(\alpha_p) = (\gamma_1, \gamma_2, \ldots, \gamma_r)$, avec $r \ge 1$, où pour tout t, $1 \le t \le r$, $x_{\gamma t}$ = [symbole].
Donc $\gamma_t \in$ Zouc(M) et $L(\text{fan}(\gamma_t))$ = [symbole][symbole].
On posera alors $\Delta G(\alpha_p) = (\text{fan}(\gamma_1), \text{fan}(\gamma_2), \ldots, \text{fan}(\gamma_r))$.
Ce qui induit que $L(\Delta G(\alpha_p))$ = ([symbole][symbole])$^r$.

**5.8.2. La fonction T : $A(M) \to E(M)^*$ :**
Pour cela, définissons les fonctions Td, Tg et $T_k$ : $A(M) \to E(M)^*$ pour tout $k \ge 1$.
Soit p, $1 \le p \le N$.
$Tg(\alpha_p)$, $Td(\alpha_p)$ et $T_k(\alpha_p)$ seront respectivement les sous-nesils de $Rg(\alpha_p)$, $Rd(\alpha_p)$ et $R_k(\alpha_p)$ constitués des éléments $\gamma$ tels que $x_\gamma$ = [symbole] ou $\bullet$.
On posera alors $T(\alpha_p) = (Tg(\alpha_p), T_1(\alpha_p), T_2(\alpha_p), \ldots, T_{\nu(\alpha p)}(\alpha_p), Td(\alpha_p))$, avec
$T(\alpha_p) = (Tg(\alpha_p), Td(\alpha_p))$ si $\nu(\alpha_p) = 0$.

**5.8.3. La fonction S : $A(M) \to E(M)^*$ :**
Pour cela, définissons les fonctions Sd, Sg et $S_k$ : $A(M) \to E(M)^*$ pour tout $k \ge 1$.
Soit p, $1 \le p \le N$.
Définissons d'abord $S'g(\alpha_p)$, $S'd(\alpha_p)$ et $S'_k(\alpha_p)$ respectivement comme les *sous-nesiles* de $Rg(\alpha_p)$, $Rd(\alpha_p)$ et $R_k(\alpha_p)$ constitués des éléments $\gamma$ tels que $x_\gamma$ = [symbole] ou [symbole] ou [symbole].

Alors $Sg(\alpha_p)$ (resp. $Sd(\alpha_p)$ et $S_k(\alpha_p)$) se définit à partir de $S'g(\alpha_p)$ (resp. $S'd(\alpha_p)$ et $S'_k(\alpha_p)$) en remplaçant les éléments $\gamma$ tels que $x_\gamma$ = [symbol] par fan($\gamma$).
Ainsi $Sg(\alpha_p)$, $Sd(\alpha_p)$ et $S_k(\alpha_p)$ sont constitués d'éléments $\gamma$ tels que $x_\gamma$ = [symbol] ou [symbol] ou [symbol] ou [symbol] ou [symbol].
On vérifie les propriétés suivantes :

1) si $Sg(\alpha_p) \neq \varepsilon$, alors $L(Sg(\alpha_p))$ est un g-*Stratajo*.
2) si $Sd(\alpha_p) \neq \varepsilon$, alors $L(Sd(\alpha_p))$ est un d-*Stratajo*.
3) si $S_k(\alpha_p) \neq \varepsilon$, alors $L(S_k(\alpha_p))$ est un *Stratajo*.

En particulier, si $Sg(\alpha_p) \neq \varepsilon$, alors d'après la *Règle du Dallajascar*, il existe $\alpha_q$ qui précède $\alpha_p$ dans $s\Sigma(\sigma(\alpha_p))$, et $R = L(Rd(\alpha_q))\ L(Rg(\alpha_p))$ est un *Ramajo*.
Ce qui induit que $S = L(Sd(\alpha_q))\ L(Sg(\alpha_p))$ est un *Stratajo*.

**5.8.4. Proposition :**
Soit (X,n) un *Stratino Naturel* appartenant à NJ(M), avec $n \geq 1$.
Donc $\Sigma(X, n) \neq \varepsilon$.
On vérifie que $s\Sigma(X, n) \neq \varepsilon$ et notons $s\Sigma(X, n) = (\alpha_{p1}, \alpha_{p2}, \ldots, \alpha_{pu})$, avec $u \geq 1$.
Alors on vérifie les résultats suivants :
$\Sigma(X, n, 1) = (G(\alpha_{p1}), G(\alpha_{p2}), \ldots, G(\alpha_{pu}))$.
$\Sigma(X, n\#) = (\Delta G(\alpha_{p1}), \Delta G(\alpha_{p2}), \ldots, \Delta G(\alpha_{pu}))$.
$\Sigma(X, n\#, 1) = (T(\alpha_{p1}), T(\alpha_{p2}), \ldots, T(\alpha_{pu}))$.
$\Sigma(X, n+1) = (S(\alpha_{p1}), S(\alpha_{p2}), \ldots, S(\alpha_{pu}))$.
De plus, on vérifie que $L(\Sigma(X, n, 1))$ est une *Trouglyre*, que $L(\Sigma(X, n\#))$ est une *Lounagatte Simple*, que $L(\Sigma(X, n\#, 1))$ est une *Trénagatte*, et si $L(\Sigma(X, n+1)) \neq \varepsilon$, que $L(\Sigma(X, n+1))$ se décompose de façon unique en facteurs de Stratajos : $L(\Sigma(X, n+1)) = S_1\ S_2 \ldots S_r$, avec $r \geq 1$, où pour tout k, $1 \leq k \leq r$, $S_k$ est un *Stratjo*.

**5.8.5. Théorème fondamental :**
Soit $\Omega$ une carte planaire cubique enracinée.
Alors $MJ(\Omega)$ est un *Mot Jassologique Complexe*.
Notons $(\varphi_\Omega, H_\Omega)$ le *Dallajascar* correspondant à $\Omega$ et reprenons les notations précédentes pour $M = MJ(\Omega)$.
On a vu que $\Omega = \{w_0, w_1, \ldots, w_N\}$, avec $w_p <_{(\varphi\Omega,H\Omega)} w_{p+1}$, pour tout p, $0 \leq p \leq N-1$.
Alors le nombre de caractères ouvrants dans $MJ(\Omega)$ est exactement N+1.
Soit $f : E(M) \rightarrow \Omega \cup \{$[symbol], [symbol], [symbol], [symbol]$\}$ définie par $f(\alpha_p) = w_p$ pour tout p, $0 \leq p \leq N$, et $f(\gamma) = x_\gamma$ pour tout $\gamma \in E(M)$ tel que $x_\gamma$ = [symbol] ou [symbol] ou [symbol] ou [symbol].
Alors pour tout p, $0 \leq p \leq N$, $H_\Omega(w_p)$ est le *Nesile* image de $H(\alpha_p)$ par f.
De plus, $NJ(\Omega) = NJ(MJ(\Omega))$.
Notons $\Sigma_\Omega : NJ^* \rightarrow \Omega AJ^*$ et $\Sigma : NJ^* \rightarrow E(MJ(\Omega))^*$ les fonctions *Stratojasse* correspondantes.
Alors pour tout *Stratino* $X \in NJ(\Omega)$, $\Sigma_\Omega(X)$ est le *Nesile* image de $\Sigma(X)$ par f.

**5.9. Exemple**
Reprenons l'exemple de la carte $\Omega = \{w_0, w_1, \ldots, w_{10}\}$ étudiée au chapitre précédent, avec $w_0$ = a, $w_1$ = b, $w_2$ = c, $w_3$ = k, $w_4$ = d, $w_5$ = g, $w_6$ = j, $w_7$ = i, $w_8$ = h, $w_9$ = e, $w_{10}$ = f.
On a vu que $M = MJ(\Omega)$ = [symbol sequence].
Donc $M = x_1\ x_2 \ldots x_{34}$ :

| i | 1 | 2 | 3 | 4 | 5 | 6 | 7 | 8 | 9 | 10 | 11 | 12 | 13 | 14 | 15 | 16 | 17 | 18 | 19 | 20 | 21 | 22 | 23 | 24 | 25 | 26 | 27 | 28 | 29 | 30 | 31 | 32 | 33 | 34 |
|---|---|---|---|---|---|---|---|---|---|---|---|---|---|---|---|---|---|---|---|---|---|---|---|---|---|---|---|---|---|---|---|---|---|---|
| $x_i$ | [symbol] | [symbol] | [symbol] | [symbol] | [symbol] | [symbol] | [symbol] | [symbol] | [symbol] | [symbol] | [symbol] | [symbol] | [symbol] | [symbol] | [symbol] | [symbol] | [symbol] | [symbol] | [symbol] | [symbol] | [symbol] | [symbol] | [symbol] | [symbol] | [symbol] | [symbol] | [symbol] | [symbol] | [symbol] | [symbol] | [symbol] | [symbol] | [symbol] | [symbol] |

**5.9.1. La *Règle d'Emboîtement* :**
Alors $\alpha_0 = 1$, $\alpha_1 = 2$, $\alpha_2 = 3$, $\alpha_3 = 4$, $\alpha_4 = 9$, $\alpha_5 = 11$, $\alpha_6 = 12$, $\alpha_7 = 14$, $\alpha_8 = 15$, $\alpha_9 = 22$, $\alpha_{10} = 23$.
Le premier caractère fermant qui suit $x_{23}$ est $x_{26}$. Donc $\beta_{10} = 26$, avec $x_{\alpha10}\ x_{\beta10}$ = [symbol] .
Après avoir supprimé $x_{23}$ et $x_{26}$, le premier caractère fermant qui suit $x_{22}$ est $x_{28}$. Donc $\beta_9 = 28$, avec $x_{\alpha9}\ x_{\beta9}$ = [symbol] .
Après avoir supprimé $x_{22}$ et $x_{28}$, le premier caractère fermant qui suit $x_{15}$ est $x_{17}$. Donc $\beta_8 = 17$, avec $x_{\alpha8}\ x_{\beta8}$ = [symbol] .
Après avoir supprimé $x_{15}$ et $x_{17}$, le premier caractère fermant qui suit $x_{14}$ est $x_{20}$. Donc $\beta_7 = 20$, avec $x_{\alpha7}\ x_{\beta7}$ = [symbol].

Après avoir supprimé $x_{14}$ et $x_{20}$, le premier caractère fermant qui suit $x_{12}$ est $x_{13}$. Donc $\beta_6 = 13$, avec $x_{\alpha 6}\, x_{\beta 6}$ = [symbol] .
Après avoir supprimé $x_{12}$ et $x_{13}$, le premier caractère fermant qui suit $x_{11}$ est $x_{21}$. Donc $\beta_5 = 21$, avec $x_{\alpha 5}\, x_{\beta 5}$ = [symbol] .
Après avoir supprimé $x_{11}$ et $x_{21}$, le premier caractère fermant qui suit $x_9$ est $x_{30}$. Donc $\beta_4 = 30$, avec $x_{\alpha 4}\, x_{\beta 4}$ = [symbol].
Après avoir supprimé $x_9$ et $x_{30}$, le premier caractère fermant qui suit $x_4$ est $x_7$. Donc $\beta_3 = 7$, avec $x_{\alpha 3}\, x_{\beta 3}$ = [symbol].
Après avoir supprimé $x_4$ et $x_7$, le premier caractère fermant qui suit $x_3$ est $x_8$. Donc $\beta_2 = 8$, avec $x_{\alpha 2}\, x_{\beta 2}$ = [symbol] .
Après avoir supprimé $x_3$ et $x_8$, le premier caractère fermant qui suit $x_2$ est $x_{33}$. Donc $\beta_1 = 33$, avec $x_{\alpha 1}\, x_{\beta 1}$ = [symbol].
Après avoir supprimé $x_2$ et $x_{33}$, le premier caractère fermant qui suit $x_1$ est $x_{34}$. Donc $\beta_0 = 34$, avec $x_{\alpha 1}\, x_{\beta 1}$ = [symbol].
Donc M vérifie bien la *Règle d'Emboîtement*.
On remarquera également que pour tout p, $0 \leq p \leq 10$, $x_{\alpha\pi}\, x_{\beta p}$ = binôme($w_p$) (fonction définie au chapitre 4.2.5.).

**5.9.2. Le *Dallajascar* :**
E(M) = {1, 2, 3, 4, 5, 6, 9, 10, 11, 12, 14, 15, 16, 18, 19, 22, 23, 24, 25, 27, 29, 31, 32}.
Les fonctions L : E(M)* $\rightarrow$ AJ* et $\varphi$ : E(M) $\rightarrow$ E(M)$\cup${ε} sont alors définies par :

| γ | 1 | 2 | 3 | 4 | 5 | 6 | 9 | 10 | 11 | 12 | 14 | 15 | 16 | 18 | 19 | 22 | 23 | 24 | 25 | 27 | 29 | 31 | 32 |
|---|---|---|---|---|---|---|---|---|---|---|---|---|---|---|---|---|---|---|---|---|---|---|---|
| L(γ) | [symbol] | [symbol] | [symbol] | [symbol] | [symbol] | [symbol] | [symbol] | [symbol] | [symbol] | [symbol] | [symbol] | [symbol] | [symbol] | [symbol] | [symbol] | [symbol] | [symbol] | • | [symbol] | [symbol] | [symbol] | • | [symbol] |
| φ(γ) | ε | 1 | 2 | 3 | 4 | 4 | 2 | 9 | 9 | 11 | 11 | 14 | 15 | 14 | 14 | 9 | 22 | 23 | 23 | 22 | 9 | 2 | 2 |

D'où H(1) = (2), H(2) = (3, 9, 31, 32), H(3) = (4), H(4) = (5, 6), H(9) = (10, 11, 22, 29), H(11) = (12, 14),
H(14) = (15, 18, 19), H(15) = (16), H(22) = (23, 27), H(23) = (24, 25), et
H(5) = H(6) = H(10) = H(12) = H(16) = H(18) = H(19) = H(24) = H(25) = H(27) = H(29) = H(31) = H(32) = ε.

**5.9.3. La *Règle des Stratojasses* :**
La fonction *Stratino* σ : E(M) $\rightarrow$ NJ* est alors définie par :
σ(1) = ε et σ(2) = (1).
φ(3) = 2, avec σ(2) = (1) et $x_3$ = [symbol], donc σ(3) = (2).
φ(4) = 3, avec σ(3) = (2) et $x_4$ = [symbol], donc σ(4) = (2, 1).
φ(5) = 4, avec σ(4) = (2, 1) et $x_5$ = [symbol], donc σ(5) = (2#).
φ(6) = 4, avec σ(4) = (2, 1) et $x_6$ = [symbol], donc σ(6) = (2#).
φ(9) = 2, avec σ(2) = (1) et $x_9$ = [symbol], donc σ(9) = (1#, 1).
φ(10) = 9, avec σ(9) = (1#, 1) et $x_{10}$ = [symbol], donc σ(10) = (2).
φ(11) = 9, avec σ(9) = (1#, 1) et $x_{11}$ = [symbol], donc σ(11) = (2).
φ(12) = 11, avec σ(11) = (2) et $x_{12}$ = [symbol], donc σ(12) = (3).
φ(14) = 11, avec σ(11) = (2) et $x_{14}$ = [symbol], donc σ(14) = (2#, 1).
φ(15) = 14, avec σ(14) = (2#, 1) et $x_{15}$ = [symbol], donc σ(15) = (2#, 2).
φ(16) = 15, avec σ(15) = (2#, 2) et $x_{16}$ = [symbol], donc σ(16) = (3).
φ(18) = 14, avec σ(14) = (2#, 1) et $x_{18}$ = [symbol], donc σ(18) = (2#, 2).
φ(19) = 14, avec σ(14) = (2#, 1) et $x_{19}$ = [symbol], donc σ(19) = (3).
φ(22) = 9, avec σ(9) = (1#, 1) et $x_{22}$ = [symbol], donc σ(22) = (1#, 2).
φ(23) = 22, avec σ(22) = (1#, 2) et $x_{23}$ = [symbol], donc σ(23) = (2).
φ(24) = 23, avec σ(23) = (2) et $x_{24}$ = •, donc σ(24) = (2#, 1).
φ(25) = 23, avec σ(23) = (2) et $x_{25}$ = [symbol], donc σ(25) = (3).
φ(27) = 22, avec σ(22) = (1#, 2) et $x_{27}$ = [symbol], donc σ(27) = (2).
φ(29) = 9, avec σ(9) = (1#, 1) et $x_{29}$ = [symbol], donc σ(29) = (1#, 2).
φ(31) = 2, avec σ(2) = (1) et $x_{31}$ = •, donc σ(31) = (1#, 1).
φ(32) = 2, avec σ(2) = (1) et $x_{32}$ = [symbol], donc σ(32) = (2).

D'où NJ(M) = {ε, (1), (1#, 1), (1#, 2), (2), (2, 1), (2#), (2#, 1), (2#, 2), (3)}.
On remarquera que NJ(M) = NJ(Ω).
On en déduit la fonction *Stratojasse* Σ : NJ* $\rightarrow$ E(M)* :
Σ(ε) = (1).
Σ(1) = (2).

$\Sigma$(1#, 1) = (9, 31).
$\Sigma$(1#, 2) = (22, 29).
$\Sigma$(2) = (3, 10, 11, 23, 27, 32).
$\Sigma$(2, 1) = (4).
$\Sigma$(2#) = (5, 6).
$\Sigma$(2#, 1) = (14, 24).
$\Sigma$(2#, 2) = (15, 18).
$\Sigma$(3) = (12, 16, 19, 25).
Et pour tout *Stratino* X $\notin$ NJ(M), $\Sigma$(X) = $\varepsilon$.
En particulier, $\Sigma(\varepsilon)$ = (1), donc M vérifie la *Règle des Stratojasses*.

**5.9.4. La *Règle du Dallajascar* :**
Par définition, A(M) = $\{\alpha_1, \alpha_2, \ldots, \alpha_{10}\}$ = {2, 3, 4, 9, 11, 12, 14, 15, 22, 23}.
s$\Sigma$(1) = (2), s$\Sigma$(1#, 1) = (9), s$\Sigma$(1#, 2) = (22), s$\Sigma$(2) = (3, 11, 23), s$\Sigma$(2, 1) = (4), s$\Sigma$(2#) = $\varepsilon$, s$\Sigma$(2#, 1) = (14), s$\Sigma$(2#, 2) = (15), s$\Sigma$(3) = (12).
D'où Cg(M) = {2, 3, 4, 9, 12, 14, 15, 22} et Cd(M) = {2, 4, 9, 12, 14, 15, 22, 23}.

La fonction G : A(M) $\rightarrow$ E(M)* : seul $x_4$ = et 4 $\in$ H(3), avec 3 = $\alpha_2$, donc G($\alpha_2$) = (4) et pour tout p différent de 2, G($\alpha_p$) = $\varepsilon$.
Comme 3 $\notin$ Cd(M), la propriété 1) de la *Règle du Dallajascar* est vérifiée.

Déterminons maintenant les fonctions R, Rd, R', Rg et $R_k$ : A(M) $\rightarrow$ E(M)* :
H(2) = (3, 9, 31, 32), avec L(H(2)) = L(3) L(9) L(31) L(32) = .
Donc R(2) = H(2), avec L(R(2)) = .
D'où R'(2) = R(2) et Rd(2) = $\varepsilon$.
D'autre part, R'(2) = ($R'_0$, 32), avec L($R'_0$) = et L(32) = .
Or 2 $\in$ Cg(M) donc Rg(2) = $\varepsilon$ et $\nu$(2) = 1, avec $R_1$(2) = ($R'_0$, 32) = H(2), avec L($R_1$(2)) = , qui est un *Ramajo*.
De plus, $R_1$(2) est un facteur de H(2). Donc les propriétés 5) et 8) de la *Règle du Dallajascar* sont vérifiées pour 2.

H(3) = (4), avec LH(3) = L(4) = , donc R(3) = $\varepsilon$.
D'où R'(3) = Rd(3) = Rg(3) = $\varepsilon$ et $\nu$(3) = 0.
En particulier, H(3) = G(3) donc la propriété 6) de la *Règle du Dallajascar* est vérifiée pour 3.
H(4) = (5, 6), avec L(H(4)) = L(5) L(6) = , donc R(4) = $\varepsilon$.
D'où R'(4) = Rd(4) = Rg(4) = $\varepsilon$ et $\nu$(4) = 0.

H(9) = (10, 11, 22, 29), avec L(H(9)) = L(10) L(11) L(22) L(29) = .
Donc R(9) = (22, 29), avec L(R(9)) = .
D'où R'(9) = R(9) et Rd(9) = $\varepsilon$.
D'autre part, R'(9) = ($R'_0$, 29), avec L($R'_0$) = et L(29) = .
Or 9 $\in$ Cg(M) donc Rg(9) = $\varepsilon$ et $\nu$(9) = 1, avec $R_1$(9) = ($R'_0$, 29) = R(9), avec L($R_1$(9)) = , qui est un *Ramajo*.
De plus, $R_1$(9) est un facteur de H(9). Donc les propriétés 5) et 8) de la *Règle du Dallajascar* sont vérifiées pour 9.

H(11) = (12, 14), avec LH(11) = L(12) L(14) = .
Donc R(11) = H(11), avec L(R(11)) = .
D'où R'(11) = $\varepsilon$ et Rd(11) = R(11), avec L(Rd(11)) = , qui est un *d-Ramajo*. De plus, (Rd(11), G(11)) est un suffixe de H(11), donc les propriétés 2) et 6) de la *Règle du Dallajascar* sont vérifiées pour 11.
D'autre part, Rg(11) = $\varepsilon$ et $\nu$(11) = 0.

H(12) = $\varepsilon$ donc R(12) = $\varepsilon$.
D'où R'(12) = Rd(12) = Rg(12) = $\varepsilon$ et $\nu$(12) = 0.

H(14) = (15, 18, 19), avec L(H(14)) = L(15) L(18) L(19) = .
Donc R(14) = (15, 18), avec L(R(14)) = .
D'où R'(14) = R(14) et Rd(14) = $\varepsilon$.

D'autre part, R'(14) = (R'$_0$, 18), avec L(R'$_0$) = [illegible] et L(18) = [illegible].
Or 14 $\in$ Cg(M) donc Rg(14) = $\varepsilon$ et $\nu$(14) = 1, avec $R_1$(14) = (R'$_0$, 18) = R(14), avec L($R_1$(14)) = [illegible] [illegible], qui est un *Ramajo*.
De plus, $R_1$(14) est un facteur de H(14). Donc les propriétés 5) et 8) de la *Règle du Dallajascar* sont vérifiées pour 14.

H(15) = (16), avec L(H(15)) = L(16) = [illegible], donc R(15) = $\varepsilon$.
D'où R'(15) = Rd(15) = Rg(15) = $\varepsilon$ et $\nu$(15) = 0.
H(22) = (23, 27), avec L(H(22)) = L(23) L(27) = [illegible] [illegible], donc R(22) = $\varepsilon$.
D'où R'(22) = Rd(22) = Rg(22) = $\varepsilon$ et $\nu$(22) = 0.

H(23) = (24,25), avec L(H(23)) = L(24) L(25) = • [illegible].
Donc R(23) = H(23), avec L(R(23)) = • [illegible].
D'où R'(23) = R(23) et Rd(23) = $\varepsilon$.
D'autre part, R'(23) = (R'$_0$, 25), avec L(R'$_0$) = • et L(25) = [illegible].
Or 23 $\notin$ Cg(M) donc Rg(23) = (R'$_0$, 25) = H(23), avec L(Rg(23)) = • [illegible], et $\nu$(23) = 0.
LRg(23) est un *g-Ramajo* et Rg(23) est un préfixe de H(23) donc les propriétés 3) et 7) de la *Règle du Dallajascar* sont vérifiées pour 23.

Enfin, pour vérifier la propriété 4), il suffit de vérifier si R = L(Rd(11)) L(Rg(23)) est un *Ramajo*.
En effet, R = [illegible] [illegible] • [illegible] est bien un *Ramajo* par définition.

On en déduit que M vérifie la *Règle du Dallajascar*.

### 5.9.5. Les fonctions $\Delta_k$ : A(M) $\rightarrow$ E(M)* :

H(2) = ($\Delta_0$(2), $R_1$(2), $\Delta_1$(2)), avec $\Delta_0$(2) = $\Delta_1$(2) = $\varepsilon$.
H(3) = ($\Delta_0$(3), G(3)), avec $\Delta_0$(3) = $\varepsilon$.
H(4) = $\Delta_0$(4), avec L($\Delta_0$(4)) = L(5) L(6) = [illegible],
H(9) = ($\Delta_0$(9), $R_1$(9), $\Delta_1$(9)), avec L($\Delta_0$(9)) = L(10) L(11) = [illegible] et $\Delta_1$(9) = $\varepsilon$.
H(11) = ($\Delta_0$(11), Rd(11)), avec $\Delta_0$(11) = $\varepsilon$.
H(12) = $\Delta_0$(12) = $\varepsilon$.
H(14) = ($\Delta_0$(14), $R_1$(14), $\Delta_1$(14)), avec $\Delta_0$(14) = $\varepsilon$ et L($\Delta_1$(14)) = L(19) = [illegible].
H(15) = $\Delta_0$(15), avec L($\Delta_0$(15)) = [illegible].
H(22) = $\Delta_0$(22), avec L($\Delta_0$(22)) = [illegible].
H(23) = (Rg(23), $\Delta_0$(23)), avec $\Delta_0$(23) = $\varepsilon$.

### 5.9.6. La *Règle des Lounafans* :

Zouc(M) = {4, 9, 14} et Zm(M) = {14} puisque 14 est le dernier terme de Rd(11), où 11 = $\varphi$(14).
$x_4$ = [illegible], avec $\sigma$(4) = (2, 1) et $\Sigma$(2#) = (5, 6).
Or 4 ʃ 5 et 4 ʃ 6 puisque $\varphi$(5) = $\varphi$(6) = 4.
Donc fan(4) = (5, 6), avec L(fan(4)) = L(5) L(6) = [illegible].
L(fan(4)) est un *Lounafan Simple* donc la propriété 1) de la *Règle des Lounafans* est vérifiée pour 4.
De plus, [illegible] est le premier terme de L($\Delta_0$(4)), donc la propriété 3) de la *Règle des Lounafans* est vérifiée pour 4.

$x_9$ = [illegible], avec $\sigma$(9) = (1#, 1) et $\Sigma$(2) = (3, 10, 11, 23, 27, 32).
Or 9 ʃ 10 puisque $\varphi$(10) = 9, 9 ʃ 11 puisque $\varphi$(11) = 9, 9 ʃ 23 puisque $\varphi^2$(23) = $\varphi$(22) = 9, 9 ʃ 27 puisque $\varphi^2$(27) = $\varphi$(22) = 9. Mais 9 n'emboîte pas 3 ni 32 puisque $\varphi$(3) = $\varphi$(32) = $\varphi$(9) = 2.
Donc fan(9) = (10, 11, 23, 27), avec L(fan(9)) = L(10) L(11) L(23) L(27) = [illegible].
L(fan(9)) est un *Lounafan* donc la propriété 2) de la *Règle des Lounafans* est vérifiée pour 9.
De plus, [illegible] est le premier terme de L($\Delta_0$(9)), donc la propriété 3) de la *Règle des Lounafans* est vérifiée pour 9.

$x_{14}$ = [illegible], avec $\sigma$(14) = (2#, 1) et $\Sigma$(3) = (12, 16, 19, 25).
Or 14 ʃ 16 puisque $\varphi^2$(16) = $\varphi$(15) = 14, 14 ʃ 19 puisque $\varphi$(19) = 14. Mais 14 n'emboîte pas 12 puisque $\varphi$(12) = $\varphi$(14) = 11, et 14 n'emboîte pas 25 puisque $\varphi^3$(25) = $\varphi^2$(23) = $\varphi$(22) = 9 et $\varphi^2$(14) = $\varphi$(11) = 9.
Donc fan(14) = (16, 19), avec L(fan(14)) = L(16) L(19) = [illegible].

L(fan(14)) est un *Lounafan* donc la propriété 2) de la *Règle des Lounafans* est vérifiée pour 14.
Remarquons que [illegible] n'est pas le premier terme de L($\Delta_0$(14)), ce qui ne contredit pas la *Règle des Lounafans* puisque puisque 14 $\in$ Zm(M).

On en déduit que M vérifie la *Règle des Lounafans*.

**5.9.7. Les fonctions ΔG, T, S : A(M) → E(M)* :**
M est bien un *Mot Jassologique Comlexe* puisqu'il vérifie les quatre règles.
Déterminons maintenant les fonctions ΔG, T, S : A(M) → E(M)*.
H(2) = ($\Delta_0$(2), $R_1$(2), $\Delta_1$(2)), avec $R_1$(2) = (3, 9, 31, 32) et L($R_1$(2)) = [illegible].
Donc ΔG(2) = ε, T(2) = (9, 31), avec L(T(2)) = [illegible], et
S(2) = (3, fan(9), 32) = (3, 10, 11, 23, 27 , 32), avec L(S(2)) = [illegible].

H(3) = ($\Delta_0$(3), G(3)), avec G(3) = (4) et L(G(3)) = [illegible].
Donc ΔG(3) = fan(4) = (5, 6) avec L(ΔG(3)) = [illegible].
De plus, T(3) = S(3) = ε.

H(4) = $\Delta_0$(4), donc ΔG(4) = T(4) = S(4) = ε.

H(9) = ($\Delta_0$(9), $R_1$(9), $\Delta_1$(9)), avec $R_1$(9) = (22, 29) et L($R_1$(9)) = [illegible].
Donc ΔG(9) = T(9) = ε, et S(9) = (22, 29) avec L(S(9)) = [illegible].

H(11) = ($\Delta_0$(11), Rd(11)), avec Rd(11) = (12, 14) et L(Rd(11)) = [illegible].
Donc ΔG(11) = ε, T(11) = Td(11) = (14) avec L(T(11)) = [illegible], et
S(11) = Sd(11) = (12, fan(14)) = (12, 16, 19) avec L(S(11)) = [illegible].

H(12) = $\Delta_0$(12), donc ΔG(12) = T(12) = S(12) = ε.

H(14) = ($\Delta_0$(14), $R_1$(14), $\Delta_1$(14)), avec $R_1$(14) = (15, 18) et L($R_1$(14)) = [illegible].
Donc ΔG(14) = T(14) = ε, et S(14) = (15, 18) avec L(S(14)) = [illegible].

H(15) = $\Delta_0$(15), donc ΔG(15) = T(15) = S(15) = ε.

H(22) = $\Delta_0$(22), donc ΔG(22) = T(22) = S(22) = ε.

H(23) = (Rg(23), $\Delta_0$(23)), avec Rg(23) = (24, 25) et L(Rg(23)) = [illegible].
Donc ΔG(23) = ε, T(23) = Tg(23) = (24) avec L(T(23)) = [illegible], et S(23) = (25) avec L(S(23)) = [illegible].

Vérifions maintenant la proposition 5.8.4. :
sΣ(1) = (2). Alors
Σ(1, 1) = (G(2)) = ε, Σ(1#) = (ΔG(2)) = ε,
Σ(1#, 1) = (T(2)) = (9, 31), avec L(Σ(1#, 1)) = [illegible], qui est bien une *Trénagatte*.
Σ(2) = (S(2)) = (3, 10, 11, 23, 27, 32), avec L(Σ(2)) = [illegible], qui est composé d'un seul *Stratajo*.

sΣ(2) = (3, 11, 23). Alors
Σ(2, 1) = (G(3), G(11), G(23)) = (4), avec L(Σ(2, 1)) = [illegible], qui est bien une *Trouglyre*.
Σ(2#) = (ΔG(3), ΔG(11), ΔG(23)) = (5, 6), avec L(Σ(2#)) = [illegible], qui est bien une *Lounagatte Simple*.
Σ(2#, 1) = (T(3), T(11), T(23)) = (14, 24), avec L(Σ(2#, 1)) = [illegible], qui est bien une *Trénagatte*.
Σ(3) = (S(3), S(11), S(23)) = (12, 16, 19, 25), avec L(Σ(3)) = [illegible], qui est composé d'un seul *Stratajo*.

sΣ(3) = (12). Alors Σ(3, 1) = (G(12)) = ε, Σ(3#) = (ΔG(12)) = ε, Σ(3#, 1) = (T(12)) = ε, Σ(4) = (S(12)) = ε.

sΣ(1#, 1) = (9). Alors
Σ(1#, 1, 1) = (G(9)) = ε, Σ(1#, 1#) = (ΔG(9)) = ε, Σ(1#, 1#, 1) = (T(9)) = ε,
Σ(1#, 2) = (S(9)) = (22, 29), avec L(Σ(1#, 2)) = ʎʅ ɷ, qui est composé d'un seul *Stratajo*.

sΣ(1#, 2) = (22).
Alors Σ(1#, 2, 1) = (G(22)) = ε, Σ(1#, 2#) = (ΔG(22)) = ε, Σ(1#, 2#, 1) = (T(22)) = ε, Σ(1#, 3) = (S(22)) = ε.

sΣ(2, 1) = (4). Alors
Σ(2, 1, 1) = (G(4)) = ε, Σ(2, 1#) = (ΔG(4)) = ε, Σ(2, 1#, 1) = (T(4)) = ε, Σ(2, 2) = (S(4)) = ε.

sΣ(2#, 1) = (14). Alors
Σ(2#, 1, 1) = (G(14)) = ε, Σ(2#, 1#) = (ΔG(14)) = ε, Σ(2#, 1#, 1) = (T(14)) = ε,
Σ(2#, 2) = (S(14)) = (15, 18), avec L(Σ(2#, 2)) = ʎʅ ɷ, qui est composé d'un seul *Stratajo*.

sΣ(2#, 2) = (15). Alors
Σ(2#, 2, 1) = (G(15)) = ε, Σ(2#, 2#) = (ΔG(15)) = ε, Σ(2#, 2#, 1) = (T(15)) = ε, Σ(2#, 3) = (S(15)) = ε.

Il y a donc au final quatre *Stratajos* :
$S_1$ = L(3, 10, 11, 23, 27, 32), $S_2$ = L(12, 16, 19, 25), $S_3$ = L(22, 29), et $S_4$ = L(15, 18).

# 6
# Construction de la Carte Géométrique correspondant à un *Mot Jassologique Complexe*

Soit $M = x_1\, x_2 \ldots x_\eta$ un *Mot Jassologique Complexe*, avec $\eta \geq 1$ sa longueur.
En reprenant les notations précédentes, N+1 est le nombre de caractères ouvrants, et pour tout p, $0 \leq p \leq N$, on a défini $J_p = (\alpha_p, \beta_p)$, avec $x_{\alpha p}\, x_{\beta p}$ = 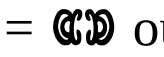 ou 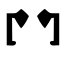 ou 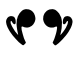 ou 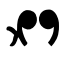 ou 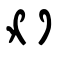 ou 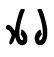 , l'ensemble E(M) et le *Dallajascar* (φ,H) sur E(M), dont $\alpha_0 = 1$ est la racine.

## 6.1. Construction de l'*Echelle Jassologique*

Elle se construit à partir de l'ensemble NJ(M) des *Stratinos* de M.
Notons $NJ(M) = \{X_0, X_1, \ldots, X_U\}$, avec $X_0 < X_1 < \ldots < X_U$ selon la relation d'ordre définie sur l'ensemble des *Stratinos*.
Comme ε et (1) appartiennent à NJ(M), alors nécessairement $X_0 = \varepsilon$ et $X_1 = (1)$.
Pour tout k, $0 \leq k \leq U$, notons $\tau_k$ la longueur du *Stratino* $X_k$.

### 6.1.1. 1ère étape :

Commençons par construire un tableau de U+1 lignes, numérotées de 0 à U en partant du haut, et de T+1 colonnes, numérotées de 0 à T en partant de la gauche.
Sur la ligne n°0, $X_0 = \varepsilon$ occupera la case de la colonne n°0, et les autres cases seront vides.
Pour tout k, $1 \leq k \leq U$, sur la ligne n°k, correspondant au *Stratino* $X_k$, le $p^{ème}$ terme de $X_k$ occupera la case de la colonne n°p pour tout p, $1 \leq p \leq \tau_k$, et les autres cases seront vides. En particulier, la case de la colonne n°0 est vide et si $\tau_k < T$, pour tout p, $\tau_k+1 \leq p \leq T$, la case de la colonne n°p est vide.

### 6.1.2. 2ème étape :

Notons $X_{u1}, X_{u2}, \ldots, X_{uc}$ l'ensemble des *Stratinos* unitaires de NJ(M), avec $c \geq 1$ et $1 = u_1 < u_2 < \ldots < u_c \leq U$. Pour tout k, $1 \leq k \leq c$, notons $Z_{uk}$ le *Stratino* tel que $X_{uk} = (Z_{uk}, 1)$.
Il existe alors $a_k$ le plus grand entier tel que $(Z_{uk}, a_k) \in NJ(M)$. Notons $r_k$, $u_k \leq r_k \leq U$, l'entier tel que $X_{rk} = (Z_{uk}, a_k)$.
Ajoutons alors une ligne vide en-dessous de la ligne correspondant à $X_{rk}$.
On notera $X'_{rk}$ cette ligne.
On vérifie que si $a_k \geq 2$, alors pour tout entier a, $2 \leq a \leq a_k$, le *Stratino* $(Z_{uk}, a)$ appartient à NJ(M).
En revanche, le *Stratino* $(Z_{uk}, a\#)$ n'appartient pas nécessairement à NJ(M). En particulier, $(Z_{uk}, 1\#)$ et $(Z_{uk}, a_k\#)$ n'appartiennent pas à NJ(M).
L'*Echelle Jassologique* correspondant à M est donc ce tableau composé de U+1+c lignes et de T+1 colonnes.
En particulier, $u_1 = 1$ et on vérifie que nécessairement $r_1 = U$.

## 6.2. Construction du *Tableau Jassologique*

Il se construit à partir de l'*Echelle Jassologique* correspondant à M.

### 6.2.1. 1ère étape :

Ajoutons η colonnes à l'*Echelle Jassologique*, η étant la longueur de M.
Numérotons-les de 1 à η en partant de la gauche.
Pour tout γ, $1 \leq \gamma \leq \eta$, le caractère $x_\gamma$ occupera la case de la ligne correspondant au *Stratino* $X_k$, tel que $X_k = \sigma(\gamma)$ si $\gamma \in E(M)$, ou $X_k = \sigma(\alpha_p)$ si $\gamma = \beta_p$ avec $0 \leq p \leq N$.
Toutes les autres cases de la colonne seront vide.
En particulier, pour tout k, $1 \leq k \leq c$, la ligne $X'_{rk}$ ne contient aucun caractère.
Mais pour tout v, $1 \leq v \leq U$, la ligne $X_v$ contient exactement tous les caractères de $L(\Sigma(X_v))$, ordonnés de gauche à droite dans leur colonne respective.

**6.2.2. 2ème étape :**
Le *Tableau Jassologique* correspondant à M s'obtient alors en supprimant toutes les colonnes correspondant aux caractères $x_{\beta p}$, $0 \le p \le N$, tels que $x_{\beta p}$ = ❟ ou ɔ ou ʖ .
Les colonnes restantes correspondent aux caractères $x_\gamma$ pour tous les entiers $\gamma \in$ E(M) ainsi qu'aux caractères $x_{\beta p}$, $0 \le p \le$ N, tels que $x_{\beta p}$ = Ɗ ou ⁊ ou ꝩ.
La ligne $X_0$ contient donc 2 caractères : le caractère $x_1$ = Ͼ dans la colonne 1 et le caractère $x_\eta$ = Ɗ dans la colonne $\eta$.
La ligne $X_1$ contient aussi 2 caractères : le caractère $x_2$ = ɼ dans la colonne 2 et le caractère $x_{\eta-1}$ = ⁊ dans la colonne $\eta$-1.
Pour tout k, $1 \le k \le c$, si $Z_{uk}$ est un *Stratino Naturel*, alors la ligne $X_{uk} = (Z_{uk}, 1)$ contient uniquement les caractères ɥ et ꝩ, et si $Z_{uk}$ est un *Stratino Décalé*, alors la ligne $X_{uk} = (Z_{uk}, 1)$ contient uniquement les caractères ɻ et •, c'est-à-dire ceux de $\Sigma(X_{uk})$ ordonnés de gauche à droite dans leur colonne respective.
Et si $a_k \ge 2$, alors pour tout entier a, $2 \le a \le a_k$, la ligne correspondant au *Stratino* $(Z_{uk}, a)$ contient uniquement les caractères ɩ, ɷ, ɤ, ʋ et ȣ, c'est-à-dire ceux de $\Sigma(Z_{uk}, a)$ ordonnés de gauche à droite dans leur colonne respective.
De plus, si $(Z_{uk}, a\#)$ appartient à NJ(M), alors la ligne correspondante contient uniquement les caractères ʋ et ȣ, c'est-à-dire ceux de $\Sigma(Z_{uk}, a\#)$ ordonnés de gauche à droite dans leur colonne respective.

**6.2.3. Repère Cartésien :**
Choisissons une même hauteur unité pour toutes les lignes et une même largeur unité pour toutes les colonnes.
La case qui contient le caractère Ͼ, c'est-à-dire la case qui se trouve sur la ligne $X_0$ et dans la colonne 1, est un rectangle : plaçons l'origine O sur le sommet en haut à gauche du rectangle, le point I sur le sommet en haut à droite du rectangle, et le sommet J en bas à droite du rectangle.
On a ainsi notre Repère Cartésien (O, I, J).
Pour tout $\gamma$, $1 \le \gamma \le \eta$, tel que $x_\gamma \ne$ ❟, ɔ et ʖ , la case qui contient le caractère $x_\gamma$ est un rectangle : on notera $A_\gamma$ le sommet en haut à gauche de ce rectangle.

**6.3. Construction des *Courbes Rovéjasses***

Pour $\alpha_1$ et tout $\alpha_p \in$ Zouc(M), nous allons construire la *Courbe Rovéjasse* $C_{\alpha p}$.
Notons k l'entier tel que $\sigma(\alpha_p) = X_{uk} = (Z_{uk}, 1)$.

**6.3.1. Cas où p = 1 :**
Alors $\alpha_1 = 2$, $\beta_1 = \eta$-1, et $\sigma(\alpha_1) = X_1$.
La ligne $X_1$ contient le caractère $x_{\alpha 1}$ = ɼ dans la colonne 2 et le caractère $x_{\beta 1}$ = ⁊ dans la colonne $\eta$-1. De plus, on a vu que $r_1$ = U. La ligne $X'_U$ est donc la dernière ligne du tableau.
On a défini les sommets $A_{\alpha 1}$ et $A_{\beta 1}$.
Notons $B_{\alpha 1}$ et $B_{\beta 1}$ les sommets situés en haut à gauche des cases situées sur la ligne $X'_U$ et dans les colonnes 2 et $\eta$-1 respectivement.
Alors $[A_{\alpha 1}, A_{\beta 1}]$ et $[B_{\alpha 1}, B_{\beta 1}]$ sont des segments horizontaux, alors que $[A_{\alpha 1}, B_{\alpha 1}]$ et $[A_{\beta 1}, B_{\beta 1}]$ sont des segments verticaux.
Posons $C_{\alpha 1} = [A_{\alpha 1}, B_{\alpha 1}] \cup [B_{\alpha 1}, B_{\beta 1}] \cup [B_{\beta 1}, A_{\beta 1}] \cup [A_{\beta 1}, A_{\alpha 1}]$, qui constitue un rectangle.
On dira que $[A_{\alpha 1}, A_{\beta 1}]$ et $[B_{\alpha 1}, B_{\beta 1}]$ sont respectivement les segments horizontaux supérieur et inférieur de la *Courbe Rovéjasse* $C_{\alpha 1}$.

**6.3.2. Cas où $x_{\alpha p}$ = ɥ :**
La ligne $X_{uk}$ contient le caractère $x_{\alpha p}$ = ɥ dans la colonne $\alpha_p$ et le caractère $x_{\beta p}$ = ꝩ dans la colonne $\beta_p$. On a défini les sommets $A_{\alpha p}$ et $A_{\beta p}$.
Notons $B_{\alpha p}$ et $B_{\beta p}$ les sommets situés en haut à gauche des cases situées sur la ligne $X'_{rk}$ et dans les colonnes $\alpha_p$ et $\beta_p$ respectivement.
D'après la règle des *Lounafans*, L(fan($\alpha_p$)) = ʋȣ, donc il existe deux entiers $\gamma_b$ et $\gamma_c$, tels que fan($\alpha_p$) = $(\gamma_b, \gamma_c)$, avec $x_{\gamma b}$ = ʋ et $x_{\gamma c}$ = ȣ. De plus par définition $\alpha_p$ ʃ $\gamma_b$ et $\alpha_p$ ʃ $\gamma_c$.
Donc $\alpha_p < \gamma_b < \gamma_c < \beta_p$.
Notons v l'entier tel que $\sigma(\gamma_b) = \sigma(\gamma_c) = X_v$, sachant que $X_v$ est un *Stratino Décalé* : si on note Y le *Stratino* et n l'entier $\ge$ 2 tel que $X_v$ = (Y, n#), alors $X_{uk}$ = (Y, n, 1).
On en déduit que le *Stratino* $X_{rk} = (Y, n, a_k) < X_v$.

Donc la ligne $X'_{rk}$ se situe entre la ligne $X_{rk}$ et la ligne $X_v$.
On a défini les sommets $A_{\gamma b}$ et $A_{\gamma c}$.
Notons $B_{\gamma b}$ et $B_{\gamma c}$ les sommets situés en haut à gauche des cases situées sur la ligne $X'_{rk}$ et dans les colonnes $\gamma_b$ et $\gamma_c$ respectivement.
$A_{\gamma b}$ se trouve donc juste en-dessous de $B_{\gamma b}$ sur la même perpendiculaire, et $A_{\gamma c}$ juste en-dessous de $B_{\gamma c}$ sur la même perpendiculaire.
De plus, $B_{\alpha p}$, $B_{\gamma b}$, $B_{\gamma c}$ et $B_{\beta p}$ se suivent dans cet ordre de gauche à droite sur la même horizontale.
Posons alors : $C_{\alpha p} = [A_{\alpha p}, B_{\alpha p}] \cup [B_{\alpha p}, B_{\gamma b}] \cup [B_{\gamma b}, A_{\gamma b}] \cup [A_{\gamma b}, A_{\gamma c}] \cup$
$[A_{\gamma c}, B_{\gamma c}] \cup [B_{\gamma c}, B_{\beta p}] \cup [B_{\beta p}, A_{\beta p}] \cup [A_{\beta p}, A_{\alpha p}]$.
On dira que $[A_{\alpha p}, A_{\beta p}]$ est le segment horizontal supérieur de $C_{\alpha p}$ et que $[B_{\alpha p}, B_{\gamma b}]$, $[A_{\gamma b}, A_{\gamma c}]$ et $[B_{\gamma c}, B_{\beta p}]$ sont les segments horizontaux inférieurs de $C_{\alpha p}$.

**6.3.3. Cas où $x_{\alpha p}$ = ⚲ :**
La ligne $X_{uk}$ contient le caractère $x_{\alpha p}$ = ⚲ dans la colonne $\alpha_p$, mais pas le caractère $x_{\beta p}$ = ৵. Mais on a vu que $L(\Sigma(X_{uk}))$ est une *Trénagatte*, donc le caractère qui suit $x_{\beta p}$ dans $L(\Sigma(X_{uk}))$ est •. Notons $\gamma$ l'entier qui suit $\alpha_p$ dans $\Sigma(X_{uk})$, qui vérifie donc $x_\gamma$ = •.
On a défini les sommets $A_{\alpha p}$ et $A_\gamma$.
Notons $B_{\alpha p}$ et $B_\gamma$ les sommets situés en haut à gauche des cases situées sur la ligne $X'_{rk}$ et dans les colonnes $\alpha_p$ et $\gamma$ respectivement.
D'après la *Règle des Lounafans*, $L(fan(\alpha_p))$ = ʊ (ꝭꝭ )$^u$ ɷ, avec $u \geq 0$, donc si on note $\gamma_b$ et $\gamma_c$ les premier et dernier termes de $fan(\alpha_p)$, $x_{\gamma b}$ = ʊ et $x_{\gamma c}$ = ɷ. De plus par définition $\alpha_p$ ⌠ $\gamma_b$ et $\alpha_p$ ⌠ $\gamma_c$.
Donc $\alpha_p < \gamma_b < \gamma_c < \beta_p < \gamma$.
Notons v l'entier tel que $\sigma(\gamma_b) = \sigma(\gamma_c) = X_v$, sachant que $X_v$ est un *Stratino Naturel* : si on note Y le *Stratino* et n l'entier $\geq$ 2 tel que $X_v = (Y, n)$, alors $X_{uk} = (Y, n-1, 1)$.
On en déduit que le *Stratino* $X_{rk} = (Y, n-1, a_k) < X_v$.
Donc la ligne $X'_{rk}$ se situe entre la ligne $X_{rk}$ et la ligne $X_v$.
On a défini les sommets $A_{\gamma b}$ et $A_{\gamma c}$.
Notons $B_{\gamma b}$ et $B_{\gamma c}$ les sommets situés en haut à gauche des cases situées sur la ligne $X'_{rk}$ et dans les colonnes $\gamma_b$ et $\gamma_c$ respectivement.
$A_{\gamma b}$ se trouve donc juste en-dessous de $B_{\gamma b}$ sur la même perpendiculaire, et $A_{\gamma c}$ juste en-dessous de $B_{\gamma c}$ sur la même perpendiculaire.
De plus, $B_{\alpha p}$, $B_{\gamma b}$, $B_{\gamma c}$ et $B_\gamma$ se suivent dans cet ordre de gauche à droite sur la même horizontale.
Posons alors : $C_{\alpha p} = [A_{\alpha p}, B_{\alpha p}] \cup [B_{\alpha p}, B_{\gamma b}] \cup [B_{\gamma b}, A_{\gamma b}] \cup [A_{\gamma b}, A_{\gamma c}] \cup$
$[A_{\gamma c}, B_{\gamma c}] \cup [B_{\gamma c}, B_\gamma] \cup [B_\gamma, A_\gamma] \cup [A_\gamma, A_{\alpha p}]$.
On dira que $[A_{\alpha p}, A_\gamma]$ est le segment horizontal supérieur de $C_{\alpha p}$ et que $[B_{\alpha p}, B_{\gamma b}]$, $[A_{\gamma b}, A_{\gamma c}]$ et $[B_{\gamma c}, B_\gamma]$ sont les segments horizontaux inférieurs de $C_{\alpha p}$.

**6.3.4. Proposition :**
Les *Courbes Rovéjasses* sont 2-à-2 disjointes.

### 6.4. Construction des *Segments Rectangulaires Stratajos*

Soit t, $1 \leq t \leq U$, tel que $X_t$ soit un *Stratino Naturel* non unitaire.
Donc $X_t = (Z, n)$, où n est un entier $\geq 2$.
Or on vérifie que nécessairement le *Stratino* Unitarie (Z, 1) appartient à NJ(M).
Donc il existe k, $1 \leq k \leq c$, tel que $Z = Z_{uk}$.
D'autre part, on a vu que $L(\Sigma(X_t))$ se décompose de façon unique en facteurs de *Stratajos*.
Alors pour tout *Stratajo* S facteur de $L(\Sigma(X_t))$, nous allons construire le *Segment Rectangulaire Stratajo* R(S).
Soit S un de ces *Stratajos*.

**6.4.1. Construction :**
Par définition, le premier terme de S est ʆ, le dernier est ɷ, et les éventuels autres termes sont ʆ ou ɫ ou ʊ ou ᴖ mais pas ɷ.
Notons p, $1 \leq p \leq N$, l'entier tel que $x_{\alpha p}$ = ʆ soit le premier terme de S, et $\gamma$, $\alpha_p < \gamma \leq \eta$, l'entier tel que $x_\gamma$ = ɷ soit le dernier.
On vérifie alors qu'il existe un et un seul entier q, $1 \leq q \leq N$, tel que $\alpha_q \in \Sigma(X_{uk})$ et tel que $\alpha_q$ *emboîte* $\alpha_p$ et $\gamma$ : $\alpha_q$ ⌠ $\alpha_p$ et $\alpha_q$ ⌠ $\gamma$. Autrement dit : $\alpha_q < \alpha_p < \gamma < \beta_p$.
Ce qui induit que $x_{\alpha p}$ et $x_\gamma$ se trouvent à l'intérieur de la surface délimitée par la *Courbe Rovéjasse* $C_{\alpha q}$. Alors les points $A_{\alpha p}$ et $A_\gamma$ se projettent chacun verticalement sur un et un seul des segments horizontaux inférieurs de $C_{\alpha q}$. : notons $B_{\alpha p}$ et $B_\gamma$ leurs projetés respectifs.
Posons alors : $R(S) = [B_{\alpha p}, A_{\alpha p}] \cup [A_{\alpha p}, A_\gamma] \cup [A_\gamma, B_\gamma]$.
En particulier, $[B_{\alpha p}, A_{\alpha p}]$ et $[A_\gamma, B_\gamma]$ sont verticaux, alors que $[A_{\alpha p}, A_\gamma]$ est horizontal.
On dira que $[B_{\alpha p}, A_{\alpha p}]$ est le segment vertical gauche et $[A_\gamma, B_\gamma]$ le segment vertical droit de R(S).

**6.4.2. Proposition :**
Notons $L(\Sigma(X_t)) = S_1 S_2 \ldots S_r$, avec $r \geq 1$, où pour tout k, $1 \leq k \leq r$, $S_k$ est un *Stratjo*.
Alors $R(S_1)$, $R(S_2)$, …, $R(S_r)$ sont deux à deux disjoints.
En particulier, leurs segments horizontaux sont alignés dans cet ordre et de gauche à droite sur la même verticale.

## 6.5. Construction des *Arêtes Transversales* aux *Segments Rectangulaires Stratajos*

Soit t, $1 \leq t \leq U$, tel que $X_t$ soit un *Stratino Naturel* non unitaire et k l'entier tel que $X_t = (Z_{uk}, n)$, où n est un entier $\geq 2$.
Soit S un *Stratajo* facteur de $L(\Sigma(X_t))$.
Nous allons construire les *Arêtes Transversales* à R(S).

Notons $\Sigma$ le facteur de $\Sigma(X_t)$ tel que $S = L(\Sigma)$, et $s\Sigma = (\alpha_{p0}, \alpha_{p1}, \ldots, \alpha_{pr})$, avec $r \geq 1$, le sous-nesile de $\Sigma$ constitué des éléments appartenant à A(M).
Si r = 0, alors il n'y aura pas d'*Arête Transversale* à R(S).
Supposons que $r \geq 2$.
Il y aura alors u-1 *Arêtes Transversales* à R(S), notées $\phi_1(S)$, $\phi_2(S)$, …, $\phi_r(S)$.
Notons $\gamma$ l'entier tel que $x_\gamma$ = ɷ soit le dernier terme de S.
Par définition, $[A_{\alpha p0}, A_\gamma]$ est le segment horizontal de R(S).
Comme $\alpha_{p0} < \alpha_{p1} < \ldots < \alpha_{pr} < \gamma$, on en déduit que les points $A_{\alpha p0}$, $A_{\alpha p1}$, …, $A_{\alpha pr}$, $A_\gamma$, se suivent dans cet ordre sur $[A_{\alpha p0}, A_\gamma]$.
De plus, on a vu qu'il existe un et un seul entier q, $1 \leq q \leq N$, tel que $\alpha_q \in \Sigma(X_{uk})$ et tel que $\alpha_q$ *emboîte* $\alpha_{p0}$ et $\gamma$. Alors nécessairement, $\alpha_q$ emboîte $\alpha_{pi}$ pour tout i, $1 \leq i \leq r$, et donc $x_{\alpha pi}$ se trouve à l'intérieur de la surface délimitée par la *Courbe Rovéjasse* $C_{\alpha q}$.
Soit i, $1 \leq i \leq r$ et construisons $\phi_i(S)$.

**6.5.1. Cas où $G(\alpha_{p(i-1)}) = \varepsilon$ :**
Alors $A_{\alpha pi}$ se projette verticalement sur un et un seul des segments horizontaux inférieurs de $C_{\alpha q}$ : notons $B_{\alpha pi}$ son projeté.
1er cas : $Rg(\alpha_{pi}) = \varepsilon$.
Alors le segment $[A_{\alpha pi}, B_{\alpha pi}]$ ne rencontre aucune autre *Courbe Rovéjasse* ni aucun autre *Segment Rectangulaire Stratajo*.
On posera $\phi_i(S) = [A_{\alpha pi}, B_{\alpha pi}]$.
2ème cas : $Rg(\alpha_{pi}) \neq \varepsilon$ et le premier terme de $L(Rg(\alpha_{pi}))$ est •.
Alors le segment $[A_{\alpha pi}, B_{\alpha pi}]$ rencontre une *Courbe Rovéjasse* $C_{\alpha q'}$ telle que $\sigma(\alpha_q) = (X_t, 1)$.
Notons $C_{\alpha pi}$ le point d'intersection avec le segment horizontal supérieur de $C_{\alpha q'}$.
Alors le segment $[A_{\alpha pi}, C_{\alpha pi}]$ ne rencontre aucune autre *Courbe Rovéjasse* ni aucun autre *Segment Rectangulaire Stratajo*.
On posera $\phi_i(S) = [A_{\alpha pi}, C_{\alpha pi}]$.

3ème cas : $Rg(\alpha_{pi}) \neq \varepsilon$ et le premier terme de $L(Rg(\alpha_{pi}))$ est différent de •.
Alors le segment $[A_{\alpha pi}, B_{\alpha pi}]$ rencontre un autre *Segment Rectangulaire Stratajo* R(S') tel que S' soit un facteur de $L(\Sigma(Z_{uk}, n+1))$.
Notons $C_{\alpha pi}$ le point d'intersection avec le segment horizontal de R(S').
Alors le segment $[A_{\alpha pi}, C_{\alpha pi}]$ ne rencontre aucune autre *Courbe Rovéjasse* ni aucun autre *Segment Rectangulaire Stratajo*.
On posera $\phi_i(S) = [A_{\alpha pi}, C_{\alpha pi}]$.

**6.5.2. Cas où $G(\alpha_{p(i-1)}) \neq \varepsilon$ :**
Notons $G(\alpha_{p(i-1)}) = (\gamma_1, \gamma_2, \ldots, \gamma_v)$, avec $v \geq 1$, où pour tout j, $1 \leq j \leq v$, $x_{\gamma j}$ = ɤ et $L(fan(\gamma_j))$ = ϑ ɷ.
Donc $\Delta G(\alpha_{p(i-1)}) = (fan(\gamma_1), fan(\gamma_2), \ldots, fan(\gamma_r))$, avec $L(\Delta G(\alpha_p))$ = (ϑ ɷ)$^v$.
On sait que $G(\alpha_{p(i-1)})$ est un facteur de $\Sigma(Z_{uk}, n, 1)$ et $\Delta G(\alpha_{p(i-1)})$ un facteur de $\Sigma(Z_{uk}, n\#)$.
De plus, $G(\alpha_{p(i-1)})$ est un suffixe de $H(\alpha_{p(i-1)})$.
Notons $\gamma_b$ le premier terme de $\Delta G(\alpha_{p(i-1)})$ et $\gamma_c$ le dernier, vérifiant $x_{\gamma b}$ = ϑ et $x_{\gamma c}$ = ɷ.
Donc $\alpha_{p(i-1)} < \gamma_b < \gamma_c < \beta_{p(i-1)} < \alpha_{pi}$.
On sait que $\Delta G(\alpha_{p(i-1)})$ est un facteur de $\Sigma(Z_{uk}, n\#)$, donc $\gamma_b$ et $\gamma_c$ se trouvent sur la ligne $(Z_{uk}, n\#)$.
Notons $A'_{\alpha pi}$ le sommet en haut à gauche du rectangle de la case se trouvant sur la ligne $(Z_{uk}, n\#)$ et dans la colonne $\alpha_{pi}$.
Ce qui induit que $[A_{\alpha pi}, A'_{\alpha pi}]$ est un segment vertical et $[A'_{\alpha pi}, A_{\gamma b}]$ un segment horizontal.
Alors $A_{\gamma b}$ se projette verticalement sur un et un seul des segments horizontaux inférieurs de $C_{\alpha q}$ : notons $B_{\gamma b}$ son projeté.
1er cas : $Rg(\alpha_{pi}) = \varepsilon$.
Alors le segment $[A_{\gamma b}, B_{\gamma b}]$ ne rencontre aucune autre *Courbe Rovéjasse* ni aucun autre *Segment Rectangulaire Stratajo.*
On posera $\phi_i(S) = [A_{\alpha pi}, A'_{\alpha pi}] \cup [A'_{\alpha pi}, A_{\gamma b}] \cup [A_{\gamma b}, B_{\gamma b}]$.
2ème cas : $Rg(\alpha_{pi}) \neq \varepsilon$ et le premier terme de $L(Rg(\alpha_{pi}))$ est •.
Alors le segment $[A_{\gamma b}, B_{\gamma b}]$ rencontre une *Courbe Rovéjasse* $C_{\alpha q'}$ telle que $\sigma(\alpha_q) = (X_t, 1)$.
Notons $C_{\gamma b}$ le point d'intersection avec le segment horizontal supérieur de $C_{\alpha q'}$.
Alors le segment $[A_{\gamma b}, C_{\gamma b}]$ ne rencontre aucune autre *Courbe Rovéjasse* ni aucun autre *Segment Rectangulaire Stratajo.*
On posera $\phi_i(S) = [A_{\alpha pi}, A'_{\alpha pi}] \cup [A'_{\alpha pi}, A_{\gamma b}] \cup [A_{\gamma b}, C_{\gamma b}]$.
3ème cas : $Rg(\alpha_{pi}) \neq \varepsilon$ et le premier terme de $L(Rg(\alpha_{pi}))$ est différent de •.
Alors le segment $[A_{\gamma b}, B_{\gamma b}]$ rencontre un autre *Segment Rectangulaire Stratajo* R(S') tel que S' soit un facteur de $L(\Sigma(Z_{uk}, n+1))$.
Notons $C_{\gamma b}$ le point d'intersection avec le segment horizontal de R(S').
Alors le segment $[A_{\gamma b}, C_{\gamma b}]$ ne rencontre aucune autre *Courbe Rovéjasse* ni aucun autre *Segment Rectangulaire Stratajo.*
On posera $\phi_i(S) = [A_{\alpha pi}, A'_{\alpha pi}] \cup [A'_{\alpha pi}, A_{\gamma b}] \cup [A_{\gamma b}, C_{\gamma b}]$.

**6.5.3. Proposition :**
$\phi_1(S), \phi_2(S), \ldots, \phi_r(S)$ sont deux à deux disjointes.

## 6.6. La Carte Géométrique

**6.6.1. Proposition :**
Après avoir construit les *Courbes Rovéjasses*, les *Segments Rectangulaires Stratajos* et les *Arêtes Transversales* aux *Segments Rectangulaires Stratajos*, on obtient une carte planaire cubique.
En particulier, cette carte contient exactement N+1 cellules qu'on peut noter $w_0, w_1, \ldots, w_N$, où pour tout p, $0 \leq p \leq N$, $w_p$ est l'unique cellule contenant la case du caractère $\alpha_p$.

**6.6.2. Définition de l'arête orientée :**
Rappelons que $C_{\alpha 1} = [A_{\alpha 1}, B_{\alpha 1}] \cup [B_{\alpha 1}, B_{\beta 1}] \cup [B_{\beta 1}, A_{\beta 1}] \cup [A_{\beta 1}, A_{\alpha 1}]$, où $[A_{\alpha 1}, A_{\beta 1}]$ et $[B_{\alpha 1}, B_{\beta 1}]$ sont respectivement les segments horizontaux supérieur et inférieur de la *Courbe Rovéjasse* $C_{\alpha 1}$.
Notons $\gamma_a$ et $\gamma_b$ les premier et dernier termes de $\Sigma(2)$, qui vérifient nécessairement $x_{\gamma a}$ = ɩ et $x_{\gamma b}$ = ɵ.
En particulier, $x_{\gamma a}$ et $x_{\gamma b}$ sont respectivement les premier et dernier termes des *Stratajos* Préfixe et Suffixe $S_a$ et $S_b$ de $L(\Sigma(2))$ (qui sont éventuellement les mêmes).
On a défini $[A_{\gamma a}, B_{\gamma a}]$ et $[A_{\gamma b}, B_{\gamma b}]$ les segments verticaux gauche et droit de $R(S_a)$ et $R(S_b)$ respectivement, où $B_{\gamma a}$ et $B_{\gamma b}$ appartiennent à $[B_{\alpha 1}, B_{\beta 1}]$.
Alors $\lambda = [B_{\gamma b}, B_{\beta 1}] \cup [B_{\beta 1}, A_{\beta 1}] \cup [A_{\beta 1}, A_{\alpha 1}] \cup [A_{\alpha 1}, B_{\alpha 1}] \cup [B_{\alpha 1}, B_{\gamma a}]$ est une arête de la carte qui est à la frontière entre $w_0$ et $w_1$.
Orientons-là de $B_{\gamma b}$ vers $B_{\gamma a}$.
La carte géométrique correspondant à M, qu'on notera $\Omega(M)$, sera donc cette carte planaire cubique enracinée dont $\lambda$ est l'arête distinguée.

**6.6.3. Théorème fondamental :**
$MJ(\Omega(M)) = M$ : le *Mot Jassologique* correspondant à $\Omega(M)$ est exactement M.

**6.6.4. Corollaire :**
Soit $\Omega_0$ une carte planaire cubique enracinée et notons $M = MJ(\Omega_0)$.
Alors $\Omega(M)$ est équivalente à $\Omega_0$ puisque $MJ(\Omega(M)) = MJ(\Omega_0) = M$.

**6.6.5. Corollaire fondamental :**
On obtient une bijection entre l'ensemble des *Mots Jassologiques Complexes* et celui des classes d'équivalence des cartes planaires cubiques enracinées.

**6.7. Exemple**
Reprenons l'exemple de la carte $\Omega = \{w_0, w_1, \dots, w_{10}\}$ étudiée aux chapitres précédents, avec $w_0 = a$, $w_1 = b$, $w_2 = c$, $w_3 = k$, $w_4 = d$, $w_5 = g$, $w_6 = j$, $w_7 = i$, $w_8 = h$, $w_9 = e$, $w_{10} = f$.
On a vu que $M = MJ(\Omega) =$ [illegible] est un *Mot Jassologique Complexe*.
$M = x_1\, x_2 \dots x_{34}$ :

| i | 1 | 2 | 3 | 4 | 5 | 6 | 7 | 8 | 9 | 10 | 11 | 12 | 13 | 14 | 15 | 16 | 17 | 18 | 19 | 20 | 21 | 22 | 23 | 24 | 25 | 26 | 27 | 28 | 29 | 30 | 31 | 32 | 33 | 34 |
|---|---|---|---|---|---|---|---|---|---|---|---|---|---|---|---|---|---|---|---|---|---|---|---|---|---|---|---|---|---|---|---|---|---|---|
| $x_i$ | [illegible] | [illegible] | [illegible] | [illegible] | [illegible] | [illegible] | [illegible] | [illegible] | [illegible] | [illegible] | [illegible] | [illegible] | [illegible] | [illegible] | [illegible] | [illegible] | [illegible] | [illegible] | [illegible] | [illegible] | [illegible] | [illegible] | [illegible] | [illegible] | [illegible] | [illegible] | [illegible] | [illegible] | [illegible] | [illegible] | [illegible] | [illegible] | [illegible] | [illegible] |

**6.7.1. Construction de l'*Echelle Jassologique* et du *Tableau Jassologique* (Figure 4) *:***
On a vu (cf. Chapitre 5.9.3.) que NJ(M) = {ε, (1), (1#, 1), (1#, 2), (2), (2, 1), (2#), (2#, 1), (2#, 2), (3)}, qui est ordonné dans cet ordre selon la relation d'ordre sur les S*tratinos*.
Les *Stratinos* unitaires étant (1), (1#, 1), (2, 1) et (2#, 1), l'*Echelle Jassologique* correspondante est alors la suivante :

| | j0 | j1 | j2 |
|---|---|---|---|
| $X_0$ | ε | | |
| $X_1$ | | 1 | |
| $X_2$ | | 1# | 1 |
| $X_3$ | | 1# | 2 |
| $X'_3$ | | | |
| $X_4$ | | 2 | |
| $X_5$ | | 2 | 1 |
| $X'_5$ | | | |
| $X_6$ | | 2# | |
| $X_7$ | | 2# | 1 |
| $X_8$ | | 2# | 2 |
| $X'_8$ | | | |
| $X_9$ | | 3 | |
| $X'_9$ | | | |

On en déduit le *Tableau Jassologique* correspondant (cf. le chapitre 5.9.3. où on a calculé $\sigma(\gamma)$ pour tout $\gamma \in E(M)$) :

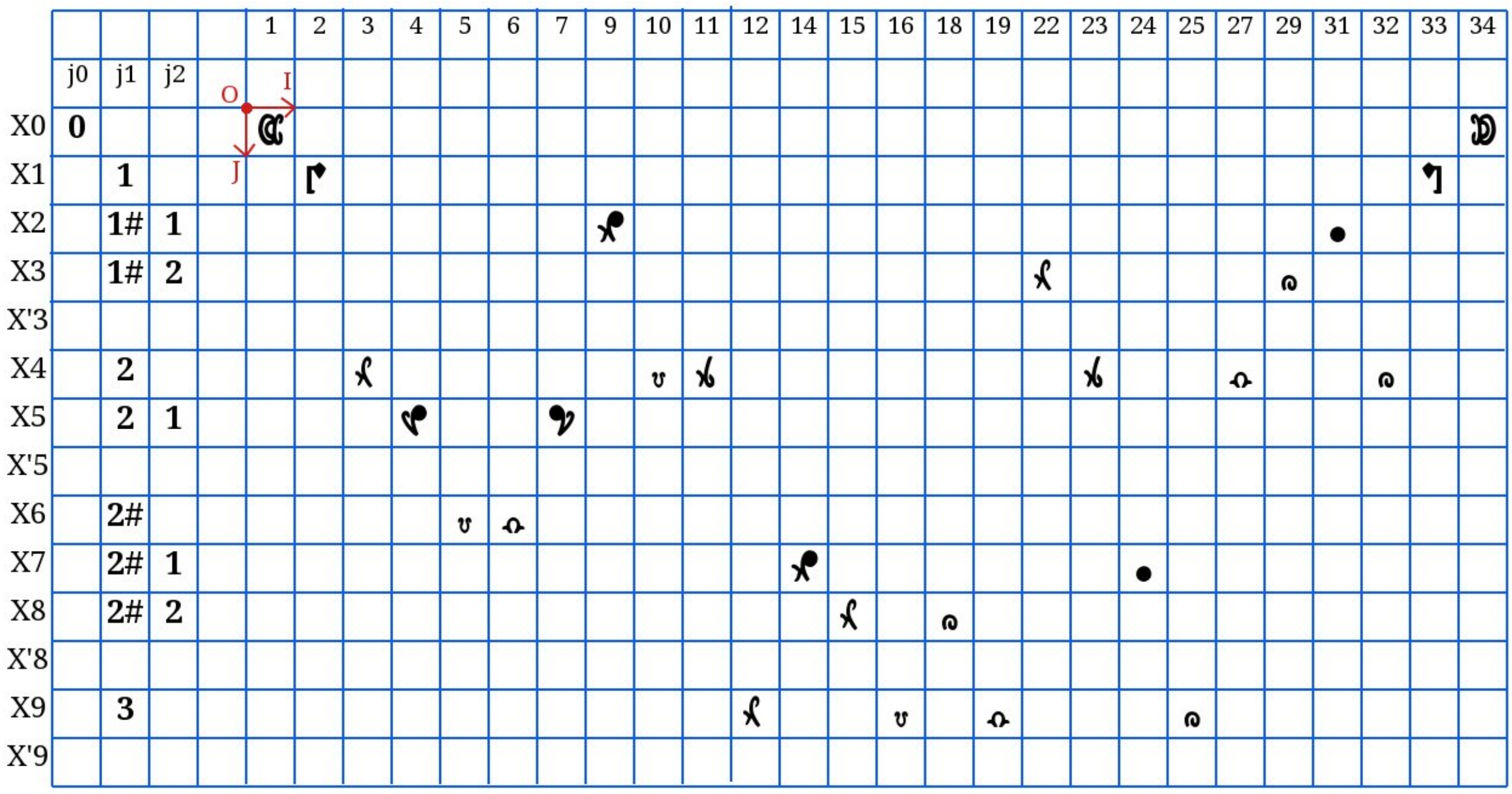

Figure 4. Le *Tableau Jassologique* de

M = MJ(Ω) = 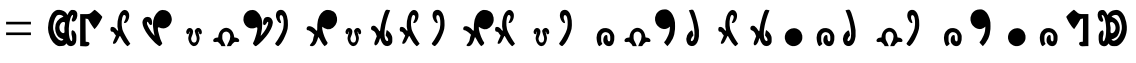

**6.7.2. Construction des *Courbes Rovéjasses* (Figure 5) :**

Pour $\alpha_1 = 2$ : $\beta_1 = 33$ et $\sigma(\alpha_1) = X_1$.

Alors $C_2 = [A_2, B_2] \cup [B_2, B_{33}] \cup [B_{33}, A_{33}] \cup [A_{33}, A_2]$.

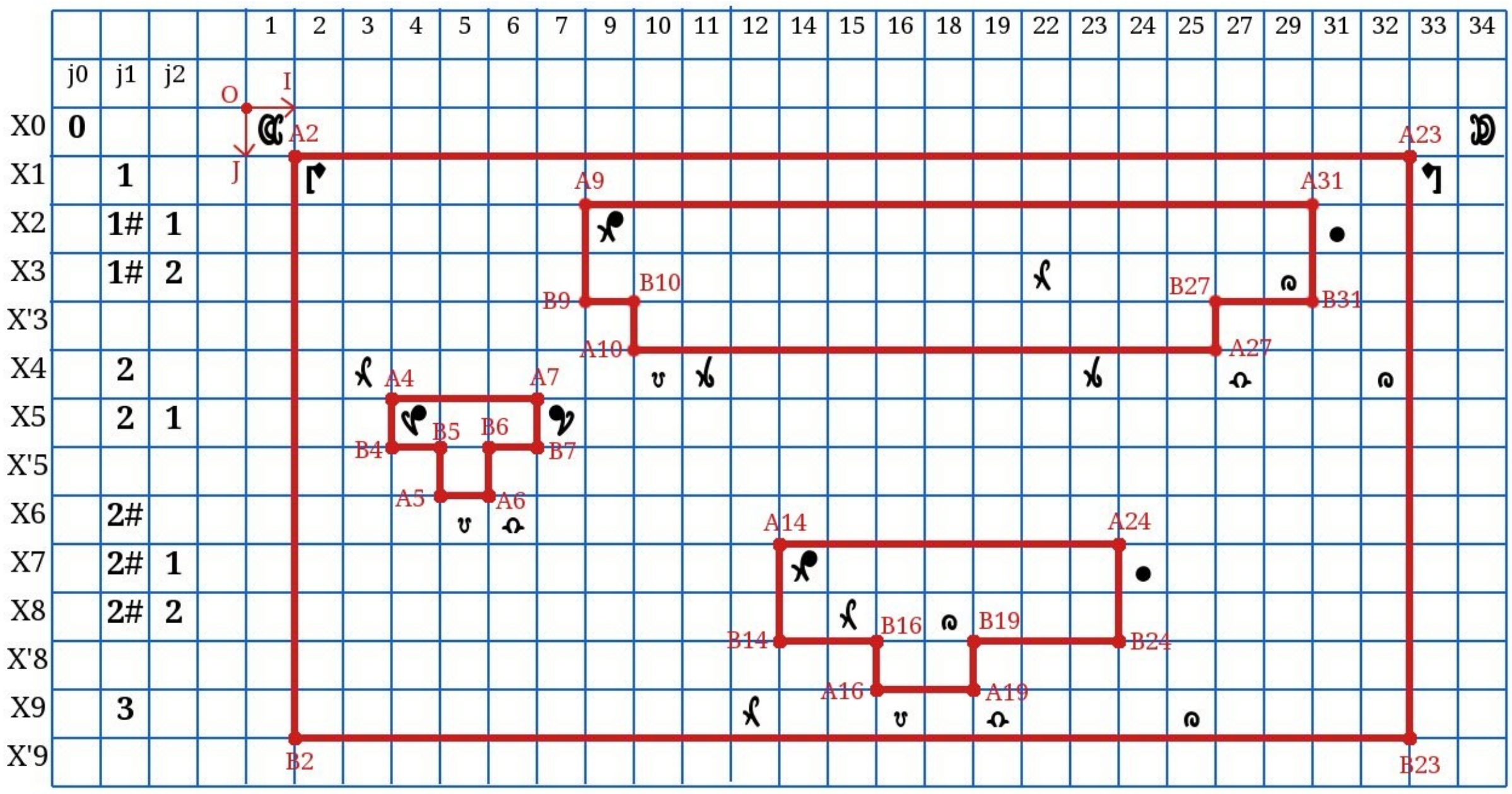

Figure 5. Les *Courbes Rovéjasses* $C_2$, $C_4$, $C_9$, et $C_{14}$ (en rouge).

Rappelons que Zouc(M) = {4, 9, 14} (cf. chapitre 5.9.6).
Pour $\alpha_p = 4$ : $x_4$ = ɤ, $\beta_p = 7$ et $\sigma(4) = (2, 1) = X_5$.
De plus, fan(4) = (5, 6), avec $x_5$ = ʊ et $x_6$ = ɷ.
Alors $C_4 = [A_4, B_4] \cup [B_4, B_5] \cup [B_5, A_5] \cup [A_5, A_6] \cup [A_6, B_6] \cup [B_6, B_7] \cup [B_7, A_7] \cup [A_7, A_4]$.
Pour $\alpha_p = 9$ : $x_9$ = ɤ, $\beta_p = 31$ et $\sigma(9) = (1\#, 1) = X_2$.
De plus, fan(9) = (10, 11, 23, 27), avec $x_{10}$ = ʊ et $x_{27}$ = ɷ.
Alors $C_9 = [A_9, B_9] \cup [B_9, B_{10}] \cup [B_{10}, A_{10}] \cup [A_{10}, A_{27}] \cup [A_{27}, B_{27}] \cup [B_{27}, B_{31}] \cup [B_{31}, A_{31}] \cup [A_{31}, A_9]$.
Pour $\alpha_p = 14$ : $x_{14}$ = ɤ, $\beta_p = 24$ et $\sigma(14) = (2\#, 1) = X_7$.
De plus, fan(14) = (16, 19), avec $x_{16}$ = ʊ et $x_{19}$ = ɷ.
Alors $C_{14} = [A_{14}, B_{14}] \cup [B_{14}, B_{16}] \cup [B_{16}, A_{16}] \cup [A_{16}, A_{19}] \cup [A_{19}, B_{19}] \cup [B_{19}, B_{24}] \cup [B_{24}, A_{24}] \cup [A_{24}, A_{14}]$.

**6.7.3. Construction des Segments *Rectangulaires Stratajos* (Figure 6) :**
On a vu que M admet quatre *Stratajos* (cf. chapitre 5.9.7.) :
$S_1$ = L(3, 10, 11, 23, 27, 32), $S_2$ = L(12, 16, 19, 25), $S_3$ = L(22, 29), et $S_4$ = L(15, 18).
Pour $S_1$ : $x_3$ = ɬ et $x_{32}$ = ɷ se trouvent à l'intérieur de la surface délimitée par la *Courbe Rovéjasse* $C_2$.
Alors $R(S_1) = [B_3, A_3] \cup [A_3, A_{32}] \cup [A_{32}, B_{32}]$, où $B_3$ et $B_{32}$ appartiennent à $[B_2, B_{33}] \subset C_2$.
Pour $S_2$ : $x_{12}$ = ɬ et $x_{25}$ = ɷ se trouvent à l'intérieur de la surface délimitée par la *Courbe Rovéjasse* $C_2$.
Alors $R(S_2) = [B_{12}, A_{12}] \cup [A_{12}, A_{25}] \cup [A_{25}, B_{25}]$, où $B_{12}$ et $B_{25}$ appartiennent à $[B_2, B_{33}] \subset C_2$.
Pour $S_3$ : $x_{22}$ = ɬ et $x_{29}$ = ɷ se trouvent à l'intérieur de la surface délimitée par la *Courbe Rovéjasse* $C_9$.
Alors $R(S_3) = [B_{22}, A_{22}] \cup [A_{22}, A_{29}] \cup [A_{29}, B_{29}]$, où $B_{22}$ et $B_{29}$ appartiennent respectivement à à $[A_{10}, A_{27}]$ et $[B_{27}, B_{31}]$ contenus dans $C_9$.
Pour $S_4$ : $x_{15}$ = ɬ et $x_{18}$ = ɷ se trouvent à l'intérieur de la surface délimitée par la *Courbe Rovéjasse* $C_{14}$.
Alors $R(S_4) = [B_{15}, A_{15}] \cup [A_{15}, A_{18}] \cup [A_{18}, B_{18}]$, où $B_{15}$ et $B_{18}$ appartiennent respectivement à à $[B_{14}, B_{16}]$ et $[A_{16}, A_{19}]$ contenus dans $C_{14}$.

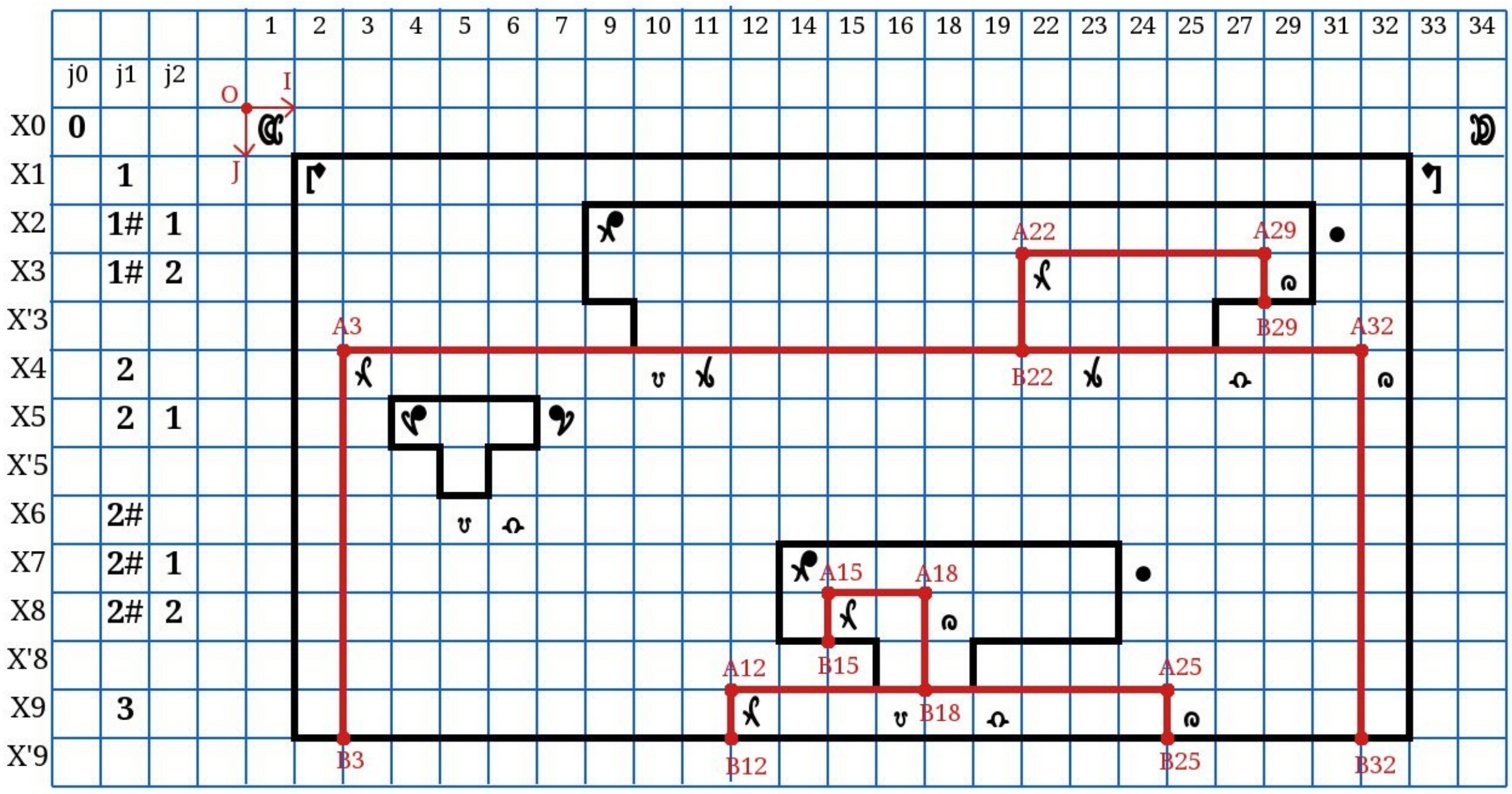

Figure 6. Les *Segments Rectangulaires Stratajos* $R(S_1)$, $R(S_2)$, $R(S_3)$, et $R(S_4)$ (en rouge).

**6.7.4. Construction des *Arêtes Transversales aux Segments Rectangulaires Stratajos* (Figure 7) :**
Seul le *Stratajo* $S_1 = L(\Sigma)$, avec $\Sigma$ = (3, 10, 11, 23, 27, 32), contient plus d'un élément de A(M) puisque $s\Sigma$ = (3, 11, 23) = $(\alpha_{p0}, \alpha_{p1}, \alpha_{p2})$.
Il y aura donc deux Arêtes Transversales $\phi_1(S_1)$ et $\phi_2(S_1)$ correspondant respectivement à $\alpha_{p1} = 11$ et $\alpha_{p2} = 23$.
Notamment, $x_{11}$ = ɬ et $x_{23}$ = ɬ se trouvent à l'intérieur de la surface délimitée par la *Courbe Rovéjasse* $C_2$.
Pour $\alpha_{p1} = 11$ : $G(\alpha_{p0}) = G(3) = (4)$, avec $x_4$ = ɤ, fan(4) = (5, 6), $x_5$ = ʊ et $x_6$ = ɷ.

Rg(11) = ε, donc on est dans le 1[er] cas : $\phi_1(S_1) = [A_{11}, A'_{11}] \cup [A'_{11}, A_5] \cup [A_5, B_5]$.
Pour $\alpha_{p2} = 23$ : $G(\alpha_{p2}) = \varepsilon$, et le premier terme de Rg(23) est •.
Donc on est dans le 2[ème] cas : $\phi_2(S_1) = [A_{23}, C_{23}]$.
Remarquons que $\phi_1(S)$ et $\phi_2(S)$ sont bien disjointes.

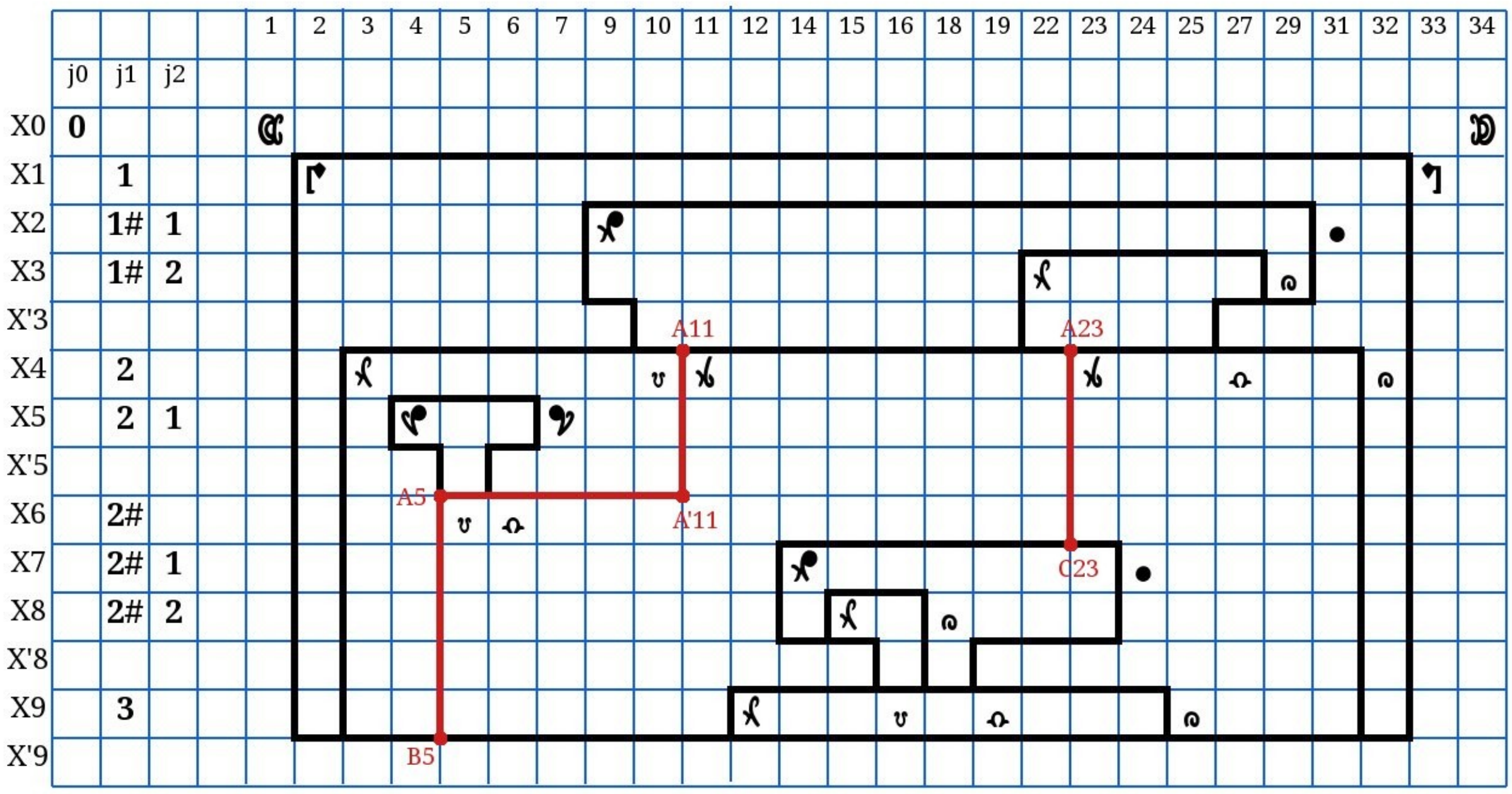

Figure 7. Les *Arêtes Transversales* $\phi_1(S_1)$ et $\phi_2(S_1)$ (en rouge).

**6.7.5. L'arête distinguée de la carte géométrique Ω(M) (Figure 8) :**
On obtient bien une carte planaire cubique contenant 11 cellules notées $w_0, w_1, \ldots, w_{11}$, où pour tout p, $0 \le p \le 10$, $w_p$ est l'unique cellule contenant la case du caractère $\alpha_p$.

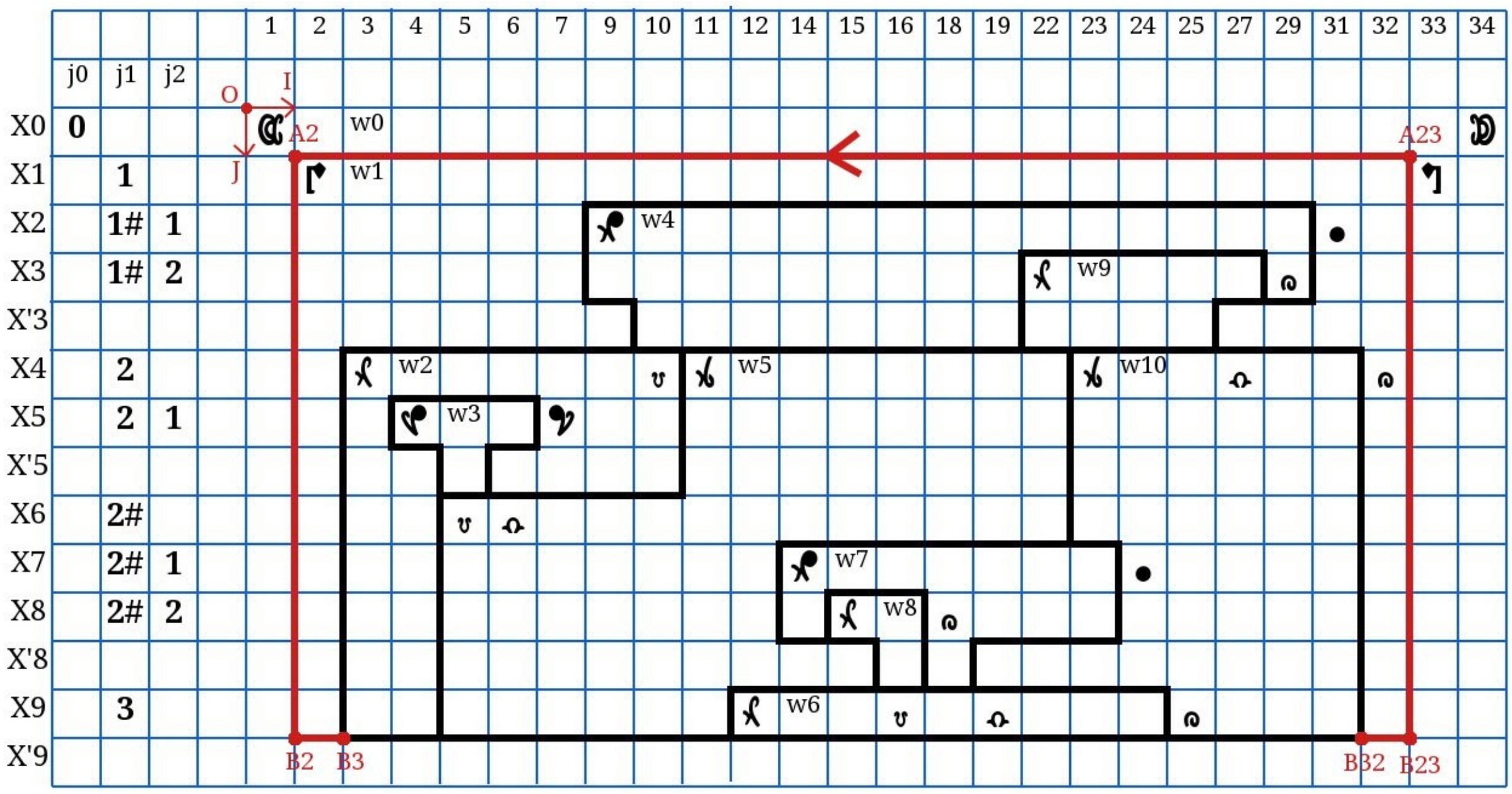

Figure 8. L'arête distinguée (en rouge) de la carte géométrique Ω(M).

Les premier et dernier termes de $\Sigma(2) = (3, 10, 11, 23, 27, 32)$ sont $\gamma_a = 3$ et $\gamma_b = 32$, avec $x_3$ = ℓ et $x_{32}$ = ɷ.
Alors $\lambda = [B_{32}, B_{33}] \cup [B_{33}, A_{33}] \cup [A_{33}, A_2] \cup [A_2, B_2] \cup [B_2, B_3]$ est par définition l'arête distinguée de Ω(M), orientée de $B_{32}$ vers $B_3$. Elle est à la frontière entre $w_0$ et $w_1$, faisant de $w_0$ la face racine (-) et $w_1$ la face racine (+).

**6.7.6. Comparaison avec la carte initiale Ω (Figure 9) :**
Remplaçons $w_0, w_1, \ldots, w_{10}$ par a, b, c, d, e, f, g, h, i, j, et k, avec $w_0$ = a, $w_1$ = b, $w_2$ = c, $w_3$ = k, $w_4$ = d, $w_5$ = g, $w_6$ = j, $w_7$ = i, $w_8$ = h, $w_9$ = e, $w_{10}$ = f.

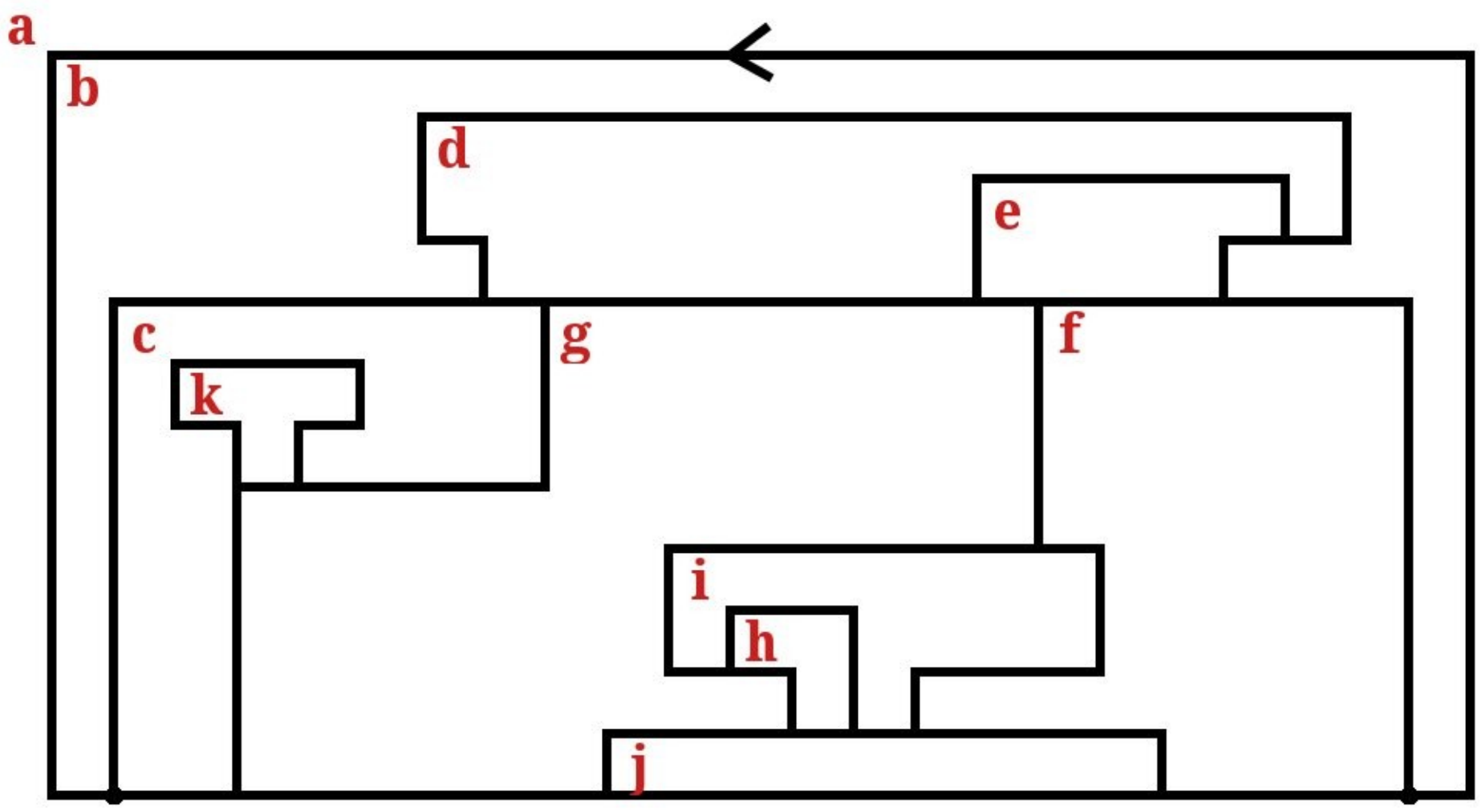

Figure 9. La carte géométrique Ω(M) est équivalente à la carte initiale Ω.

On voit que les bordures étendues sont exactement les mêmes que celles des faces de Ω :
$\overline{B}$(a) = (b, f, j, g, c, b), $\overline{B}$(b) = (a, c, d, e, f, a), $\overline{B}$(c) = (a, g, k, g, d, b, a), $\overline{B}$(d) = (c, g, e, b, c), $\overline{B}$(e) = (g, f, b, d, g), $\overline{B}$(f) = (a, b, e, g, i, j, a), $\overline{B}$(g) = (c, a, j, h, i, f, e, d, c, k, c), $\overline{B}$(h) = (g, j, i, g), $\overline{B}$(i) = (f, g, h, j, f), $\overline{B}$(j) = (i, h, g, a, f, i), et $\overline{B}$(k) = (c, g, c).
Ω(M) est donc équivalente à Ω. cqfd

## Lexique étymologique des mots inventés

**B**
baou : vient du provençal « baus » qui signifie « rocher escarpé, terrasse ».

bordure : viens du français « bord ».

**C**
caouly : vient du mot français « col de montagne ».

chio : vient du mot italien « chiodo », qui signifie « clou » mais aussi « piton d'alipiniste ». Se prononce [kio].

couno : dérivé de « trecouna ».

**D**
dallajascar : vient du mot turc « dallara », qui signifie « ramification », et du mot provençal « jas », qui signifie « dépôt, lie, stratification ».

**J**
jasse : viennent du mot provençal « jas », qui signifie « dépôt, lie, stratification », ou « jaseiran », qui signifie « collier formé de mailles, chaîne, bracelet ».

**L**
loun, louna, lounagatte : extrait du mot français « cellule ».

lounafan : vient de « loun » et du mot français « olifant », qui est un instrument de musique du Moyen Age.

**M**
mirounda, miroundo : viennent du mot provençal « mira », qui signifie « regarder, admirer », ou encore « miroundello », qui signifie « affiche, montre ».

**N**
nesile : viennent du mot turc « nesil », qui signifie « descendance, génération ».

nesilpartition : vient de « nesil » et du mot français « partition ».

**R**
ramajo : vient des mots français « ramifié » et « banjo », qui est un instrument de musique.

rôvéjasse : vient du mot hébreux « rôvéde », qui signifie « couche », et du mot provençal « jas », qui signifie « dépôt, lie, stratification ».

**S**
stougammon : vient de « stoun » et du mot français « gamme musicale ».

stoun, stouna, stratino : viennnent du mot français « strate ».

stratajo : vient des mots français « strate » et « banjo », qui est un instrument de musique.

stratojasse : vient du mot français « strate » et du mot provençal « jas », qui signifie « dépôt, lie, stratification ».

**T**
trecouna, trecouno : viennent du mot provençal « trecou », qui signifie « point culminant d'une montagne ».

trénagatte : dérivé de « trecouna ».

trouglyre : vient de « trougouna » et du mot français « lyre », qui est un instrument de musique.

trougouna, trougouno : dérivés de « trecouna ».

**Z**
zouc, zoucarei, zouccheirio : viennent du provençal « souco », qui signifie « souche ».

## Références

## Remerciements

**Mathématiciens :**
Mon Directeur de Thèse Max Karoubi (Université Paris 7).
Sylvain Cappell (Courant Institute NYU), Gérard Benarous (Courant Institute NYU), Edward Y. Miller (Polytechnics NYU).

**Non mathématiciens :**
Un immense merci à Gilles Trehin et Stéphane De Bock qui m'ont été de précieux interlocuteurs. En particulier, merci à Gilles pour ses connaissances en histoire-géographie et langues étrangères, qui m'ont largement influencée pour l'origine des mots inventés. Et merci à Stéphane, dont les compétences en informatique et programmation ont été extrêmement bienvenues.
Un très grand merci à Paul et Chantal Trehin, Chris et Ralph Rizzuto, Charles et Geneviève Durham, pour leur aide et soutien sans faille. En particulier, merci à Chantal pour m'avoir fait connaître le livre de Paul Hoffman sur la vie d'Erdös et celui de Robin Wilson sur le théorème des 4 couleurs.
Un très grand merci également à Josef Schovanec.